\documentclass[11pt]{amsart}
\usepackage{amssymb}
\usepackage{amscd}
\usepackage{curves}
\usepackage{epsfig}
\usepackage[pass]{geometry}
\usepackage{graphicx}
\usepackage{pstricks}
\usepackage{pstricks-add}
\usepackage{pst-grad}
\usepackage{pst-plot}
\usepackage{verbatim}
\usepackage{latexsym,eucal,amsfonts,amssymb,amsmath,graphicx}
\usepackage{auto-pst-pdf}
\numberwithin{equation}{section}
\newtheorem{theorem}{Theorem}[section]

\newtheorem{lemma}[theorem]{Lemma}
\newtheorem{prop}[theorem]{Proposition}

\newtheorem{definition}[theorem]{Definition}

\def \bpf {\begin{proof}}
\def \epf {\end{proof}}
\def \beq {\begin{equation*}}
\def \eeq {\end{equation*}}
\def \bsp{\begin{split}}
\def \esp{\end{split}}
\def \beqq {\begin{equation}}
\def \eeqq {\end{equation}}

\def \mcd {{\mathcal D}}
\def \mce {{\mathcal E}}
\def \mcf {{\mathcal F}}

\def \mck {{\mathcal K}}
\def \mcl {{\mathcal L}}

\def \mco {{\mathcal O}}
\def \mcp {{\mathcal P}}

\def \mcr {{\mathcal R}}
\def \mcs {{\mathcal S}}

\def \mcu {{\mathcal U}}
\def \mcv {{\mathcal V}}

\def \mcx {{\mathcal X}}
\def \mcy {{\mathcal Y}}

\def \mbn {{\mathbb N}}
\def \mbr {{\mathbb R}}

\def \id {\operatorname{Id}}
\def \comp {\operatorname{comp}}

\def \loc {\operatorname{loc}}

\def \det {\operatorname{det}}

\def \diag{\textrm{Diag}}

\def \supp {\text{supp }}

\def \eps {\epsilon}   
   
\def \La {\Lambda}

\def \p {\partial}

\def \eps {\epsilon}
\def \det {\text{det}}

\def \vol {\text{vol}}
\def \dim {\text{dim}}
\def \ha {\frac{1}{2}}

\def \fnf {\frac{n}{4}}
\def \WF {\operatorname{WF}}
\def \loc {\operatorname{loc}}
\def \comp {\operatorname{comp}}

\def \hml {H_{\text{ml}}}

\def \ba {\begin {eqnarray*} }
\def \ea {\end {eqnarray*} }
\def\R{\mathbb R}

\def\Y{{\mathcal Y}}
\def \empty {\emptyset}

\begin{document}
\title{Inverse problems for semilinear wave equations on Lorentzian manifolds}

\author{Matti Lassas} 
\address{Matti Lassas
\newline
\indent Department of Mathematics, University of Helsinki}
\email{matti.lassas@helsinki.fi}

\author{Gunther Uhlmann} 
\address{Gunther Uhlmann
\newline
\indent Department of Mathematics, University of Washington, Seattle, 
\newline
\indent Institute for Advanced Study, the Hong Kong University of Science and Technology 
\newline
\indent and Department of Mathematics, University of Helsinki}
\email{gunther@math.washington.edu}

\author{Yiran Wang}
\address{Yiran Wang
\newline
\indent Department of Mathematics, University of Washington, Seattle 
\newline
\indent and Institute for Advanced Study, the Hong Kong University of Science and Technology}
\email{wangy257@math.washington.edu}

\begin{abstract} 
We consider inverse problems in space-time $(M, g)$, a $4$-dimensional Lorentzian manifold. For semilinear wave equations $\square_g u + H(x, u) = f$, where $\square_g$ denotes the usual Laplace-Beltrami operator, we prove that the source-to-solution map $L: f \rightarrow u|_V$, where $V$ is a neighborhood of a time-like geodesic $\mu$, determines the topological, differentiable structure and the conformal class of the metric of the space-time in the maximal set where waves can propagate from $\mu$ and return back.  Moreover, on a given space-time $(M, g)$, the source-to-solution map determines some coefficients of the Taylor expansion of $H$ in $u$ .
\end{abstract}

\maketitle
\tableofcontents

\section{Introduction}
\subsection{The inverse problem}
We study inverse problems for semilinear wave equations on a $4$-dimensional Lorentzian manifold. To set up the problem, we briefly recall the background from Lorentzian geometry. The details and references can be found in Section \ref{prel}.

\begin{figure}[htbp]
\centering
\scalebox{0.7} 
{
\begin{pspicture}(0,-4.68)(6.8418946,4.68)
\psellipse[linewidth=0.04,dimen=outer](3.11,-0.13)(0.55,3.69)
\psbezier[linewidth=0.04](3.0557144,-4.66)(3.26,-3.04)(3.2428572,-2.26)(3.1114285,-0.18)(2.98,1.9)(3.0557144,3.58)(3.2414286,4.66)
\psellipse[linewidth=0.04,linestyle=dashed,dash=0.16cm 0.16cm,dimen=outer](3.1,-0.23)(3.1,0.61)
\psdots[dotsize=0.12](3.08,2.84)
\psdots[dotsize=0.12](3.18,-3.18)
\psline[linewidth=0.04cm,linestyle=dashed,dash=0.16cm 0.16cm](3.12,2.84)(6.1,-0.08)
\psline[linewidth=0.04cm,linestyle=dashed,dash=0.16cm 0.16cm](3.04,2.84)(0.08,-0.1)
\psline[linewidth=0.04cm,linestyle=dashed,dash=0.16cm 0.16cm](6.16,-0.3)(3.18,-3.2)
\psline[linewidth=0.04cm,linestyle=dashed,dash=0.16cm 0.16cm](0.04,-0.3)(3.16,-3.2)
\usefont{T1}{ptm}{m}{n}
\rput(4.5,2.965){$p_+ = \hat\mu(s_+)$}
\usefont{T1}{ptm}{m}{n}
\rput(4.5,-3.255){$p_- = \hat\mu(s_-)$}
\usefont{T1}{ptm}{m}{n}
\rput(1.5614551,-0.19){$I(p_-, p_+)$}
\usefont{T1}{ptm}{m}{n}
\rput(2.95,-2.235){$V$}
\usefont{T1}{ptm}{m}{n}
\rput(3.5071455,4.385){$\hat\mu$}
\end{pspicture} 
}
\caption{$\hat\mu$ is a time-like geodesic. The set $V$ is an open neighborhood of $\hat\mu([-1, 1])$. The source $f$ is supported in $V$ and we take measurements $L(f)$ in $V$. The set $I(p_-, p_+)$ is the set where the wave can propagate to from $\hat \mu$ and return back to $\hat \mu$. We study the inverse problem of determining the metric and the nonlinearity in $I(p_-, p_+)$ bounded by the dashed curves.}
\label{figinv}
\end{figure}
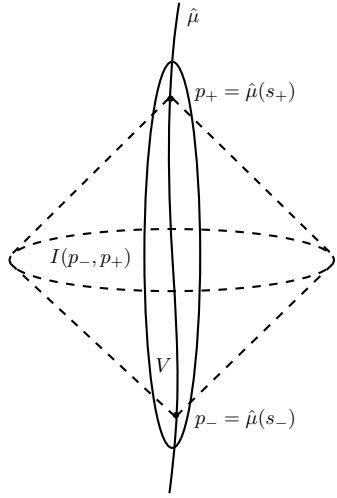

Let $(M, g)$ be an $1+ 3$ dimensional time oriented globally hyperbolic Lorentzian manifold. For $p, q\in M$, we denote by $p\ll q$ if $p\neq q$ and $p$ can be joined to $q$ by a future pointing time-like curve. We denote by $p<q$ if $p\neq q$ and $p$ can be joined to $q$ by a future pointing causal curve. We use $p\leq q$ if $p = q$ or $p<q$. The chronological future of $p\in M$ is denoted by $I^+(p) = \{q\in M: p\ll q\}$. The causal future of $p\in M$ is $J^+(p) = \{q\in M: p\leq q\}$. The chronological past and causal past are denoted by $I^-(p)$ and $J^-(p)$ respectively. For any set $A\subset M$, we denote $J^\pm(A) = \cup_{p\in A}J^\pm(p)$. Also, we denote $J(p, q) = J^+(p)\cap J^-(q)$ and $I(p, q) = I^+(p)\cap I^-(q)$. See Fig. \ref{figinv}. When it becomes necessary, we use subscript to indicate the dependence on $g$.

For globally hyperbolic Lorentzian manifold, $(M, g)$ can be identified with the product manifold 
\beq
\mbr\times N \text{ with metric } g = -\beta(t, y) dt^2 + \kappa(t, y),
\eeq
where $N$ is a $3$-dimensional manifold, $\beta$ is smooth and $\kappa$ is a family of Riemannian metrics on $N$ smoothly depending on $t$, see \cite{BS0} . Let $\hat \mu(t) \subset M$ be a time-like geodesic where $t\in [-1, 1]$. Let $p_\pm = \hat\mu(s_\pm), -1<s_-<s_+<1$ be two points on $\hat \mu$ and $V$ be an open relatively compact neighborhood of $\hat \mu([s_-, s_+])$. Take $M_0 = (-\infty, T_0)\times N, T_0>0$ such that $V\subset M_0$. Let $\square_g$ be the Laplace-Beltrami operator on $(M, g)$, we consider semilinear wave equations with source terms
\beqq\label{eqsem}
\begin{gathered}
\square_{g} u(x) + H(x, u(x))  = f(x), \text{ on } M_0,\\
u = 0 \text{ in } M_0\backslash J^+_{g}(\supp(f)),
\end{gathered}
\eeqq
where $x = (t, y) \in \mbr\times N, \supp(f)\subset V$ and $H$ is smooth. The local well-posedness of this problem has been studied in \cite{KLU1}, see also  \cite[Appendix III]{Cb}, \cite{So, Tay}.  Roughly speaking, there is a unique solution $u$ for $f$ small in $C^m$-norm with suitable $m$. As Definition 1.4 of \cite{KLU}, we define the source-to-solution map $L=L_{V; \square_g,H}$ as 
\beq
L(f) = u|_{V},
\eeq
where $u$ is the solution to \eqref{eqsem} with source $f$. Assume that we are given $V$ as a differentiable manifold and the map $L$. The inverse problem (of active measurements) we study in this work is whether one can determine the metric $g$ and the nonlinear term $H$ on $I(p_-, p_+)$ from these information. See Fig. \ref{figinv}.

\subsection{Determination of metrics and nonlinearities}
When $H(x, u) = a(x)u(x)^2$, the inverse problem was studied by Kurylev, Lassas and Uhlmann in \cite{KLU}. The same problem has been proposed and studied for the Einstein equations with matter sources in \cite{KLU1, KLU2, KLU3}. The main result Theorem 1.5 of \cite{KLU} states that if $H(x, u) = a(x)u(x)^2$ with $a(x)$ non-vanishing,  we can determine the conformal class of the metric $g$ up to diffeomorphisms. Under some additional assumptions e.g.\ the manifolds are Ricci flat, the authors in \cite{KLU} proved that the metric is uniquely determined up to diffeomorphisms. Similar results also hold for the Einstein equation with matter sources, see Theorem 1.1 of \cite{KLU1}.

In this work, we consider a general nonlinear term $H$. Roughly speaking, 
we prove that the Lorentzian metric can be determined up to diffeomorphisms if the nonlinearity $H(x, z)$ satisfy certain assumptions. Also, we show that on a given Lorentzian manifold, the source-to-solution map determines the nonlinear term $H$. We now state the precise theorems.

We start with the meaning of nonlinearity used in this work. 
\begin{definition}\label{defnl}
Let $H(x, z)\in C^\infty(U\times I)$ be real valued, where $U$ is open in $M$ and $I$ is a small neighborhood of $0$ in $\mbr$. We say $H$ is \textbf{genuinely nonlinear}\footnote{We thank Peter Hintz for this notion in private communication.} in $z$ on $U$ if $H(x, 0) = \p_z H(x, 0) = 0$  and for any $x\in U$, there is $k\in \mbn, k\geq 2$ such that 
\beq
\p^k_z H(x, 0) \neq 0.
\eeq
If there exists $k_0\geq 2$ such that $\p_z^k H(x, 0) = 0$ for all $k> k_0$ and $x\in U$, we say $k_0$ is the order of $H$. If there is no such $k_0$, we say $H$ is nonlinear of infinite order. 
\end{definition}

Our main result is the following theorem.  
\begin{theorem}\label{maingh}
Let $(M^{(j)}, g^{(j)}), j = 1, 2$ be two $4$-dimensional globally hyperbolic Lorentzian manifolds. Let $\hat \mu^{(j)}(t) \subset M^{(j)}$ be time-like geodesics where $t\in [-1, 1]$ and $V^{(j)}\subset M^{(j)}$ be open relatively compact neighborhood of $\hat \mu^{(j)}([s_-, s_+])$ where $-1<s_-<s_+<1$. Let $M_0^{(j)} = (-\infty, T_0)\times N^{(j)}, T_0 > 0$ such that $V^{(j)}\subset M_0^{(j)}$. Consider the semilinear wave equations with source terms
\beq
\begin{gathered}
\square_{g^{(j)}}u(x) + H^{(j)}(x, u(x))  = f(x), \text{ on } M^{(j)}_0,\\
u = 0 \text{ in } M_0^{(j)}\backslash J^+_{g^{(j)}}(\supp(f)),
\end{gathered}
\eeq
where $\supp(f)\subset V^{(j)}$. We assume that $H^{(j)}(x, z)$ are genuinely nonlinear on $I(p_-^{(j)}, p_+^{(j)})$ where $p^{(j)}_\pm = \hat \mu^{(j)}(s_\pm)$.  Suppose that there is a diffeomorphism $\Phi: V^{(1)}\rightarrow V^{(2)}$ such that $\Phi(p^{(1)}_\pm) = p^{(2)}_\pm$ and the source-to-solution maps $L^{(j)}$ satisfy
\beq
(\Phi^{-1})^*(L^{(1)}(\Phi^*f)) = L^{(2)}(f)
\eeq
for all $f$  in a small neighborhood of the zero function in $C_0^4(V^{(2)})$. Then there exists a diffeomorphism $\Psi: I(p^{(1)}_-, p^{(1)}_+)\rightarrow I(p^{(2)}_-, p^{(2)}_+)$ and $\gamma \in C^\infty(I(p^{(1)}_-, p^{(1)}_+))$ such that  for $x\in I(p^{(1)}_-, p^{(1)}_+)$ we have
\beqq\label{eqtrans}
\begin{gathered}
 g^{(1)} = e^{2\gamma}\Psi^*g^{(2)},\\
\p_z^k H^{(1)}(x, 0) =  e^{(k-3)\gamma(x)} \p_z^k H^{(2)}(\Psi(x), 0), \ \ \forall k\geq 4.
\end{gathered}
\eeqq
Also, for $k = 2, 3$, we have that for $x\in I(p^{(1)}_-, p^{(1)}_+)$ 
\begin{enumerate}
\item $\p_z^2 H^{(1)}(x, 0) \cdot \p_z^3 H^{(1)}(x, 0)=  e^{-\gamma(x)} \p_z^2 H^{(2)}(\Psi(x), 0) \cdot \p_z^3 H^{(2)}(\Psi(x), 0);$
\item $\p_z^2 H^{(1)}(x, 0) =  e^{-\gamma(x)} \p_z^2 H^{(2)}(\Psi(x), 0) \text{ if } \p_z^3 H^{(1)}(x, 0) = 0.$
\end{enumerate}
\end{theorem} 

We remark that in general linear terms in the wave equation do not affect the results and we give more precise statements including linear terms in Theorem \ref{mainconf} and Theorem \ref{maincoef}. 
However, the genuinely nonlinear condition is essential.  The theorem implies important consequences on unique determination of the Lorentzian metric and the nonlinear function $H$. We first consider the determination of the metric.

\begin{theorem}\label{main1}
Let $M$ be a $4$-dimensional manifold and $g^{(j)}, j = 1, 2$ be two globally hyperbolic Lorentzian metrics on $M$. Let $\hat \mu^{(j)}(t)$ be time-like geodesics on $(M, g^{(j)})$ where $t\in [-1, 1]$ and $V^{(j)}\subset M$ be open relatively compact neighborhood of $\hat \mu^{(j)}([s_-, s_+])$ where $-1<s_-<s_+<1$. Let $M_0^{(j)} = (-\infty, T_0)\times N^{(j)}, T_0 > 0$ such that $V^{(j)}\subset M_0^{(j)}$. Consider the semilinear wave equations with source terms
\beq
\begin{gathered}
\square_{g^{(j)}} u(x)  + H^{(j)}(x, u(x))  = f(x), \text{ on } M^{(j)}_0,\\
u = 0 \text{ in } M_0^{(j)}\backslash J^+_{g^{(j)}}(\supp(f)),
\end{gathered}
\eeq
where $\supp(f)\subset V^{(j)}$. Assume that $H^{(j)}(x, z)$ are genuinely nonlinear on $I(p_-^{(j)}, p_+^{(j)})$, where $p^{(j)}_\pm = \hat \mu^{(j)}(s_\pm)$. Suppose that there is a diffeomorphism $\Phi: V^{(1)}\rightarrow V^{(2)}$ such that $\Phi(p^{(1)}_\pm) = p^{(2)}_\pm$ and the source-to-solution maps $L^{(j)}$ satisfy
\beq
(\Phi^{-1})^*(L^{(1)}(\Phi^*f)) = L^{(2)}(f)
\eeq
for all $f$  in a small neighborhood of the zero function in $C_0^4(V^{(2)})$. Then there exists a diffeomorphism 
$\Psi: I(p^{(1)}_-, p^{(1)}_+)\rightarrow I(p^{(2)}_-, p^{(2)}_+)$  and $\gamma \in C^\infty(I(p^{(1)}_-, p^{(1)}_+))$ such that $ g^{(1)} = e^{2\gamma} \Psi^*g^{(2)}$ in $I(p^{(1)}_-, p^{(1)}_+)$.  Moreover, the diffeomorphism $\Psi$ is an isometry  i.e.\ $g^{(1)} = \Psi^*g^{(2)}$ on $I(p^{(1)}_-, p^{(1)}_+)$ under one of the following additional assumptions.
\begin{enumerate}
\item $H^{(i)}(x, z), i = 1, 2$ are independent of $x$ i.e.\ $H^{(i)}(x, z) = H^{(i)}(z)$ and $H^{(i)}(z)\neq b^{(i)}z^3$ for \indent some constants $b^{(i)}$;
\item The Ricci curvatures of $g^{(i)}$ are zero.
\end{enumerate}
\end{theorem} 

This theorem generalizes and improves the results obtained by Kurylev, Lassas and Uhlmann for $H(x, z) = a(x)z^2$ (Theorem 1.1 and Corollary 1.3 of \cite{KLU}). We refer to Theorem \ref{mainconf} and \ref{maincoef} for the statements with linear terms. We must point out that the unique determination of the metric is not true in general even if the nonlinear perturbations are known. We will demonstrate such examples in Section \ref{countex}, after introducing the gauge transformations. We emphasize that the case when $H$ is purely cubic needs special treatment. In particular, we can have equations with the same cubic function $H$ yet the source-to-solution maps are the same for any conformal metrics! This is related to the gauge invariance of the conformal wave equations in dimension $4$. 


Next, we state our result on the determination of the nonlinear term $H$. 
\begin{theorem}\label{main2}
Let $(M, g)$ be a $4$-dimensional globally hyperbolic Lorentzian manifold. Let $\hat \mu(t) \subset M$ be a time-like geodesic where $t\in [-1, 1]$ and  $V\subset M$ be an open relatively compact neighborhood of $\hat \mu$. Let $M_0 = (-\infty, T_0)\times N, T_0 > 0$ such that $V\subset M_0$. Consider semilinear wave equations with source terms
\beq 
\begin{gathered}
\square_{g} u(x) + H^{(j)}(x, u(x))  = f(x), \text{ on } M_0,\\
u = 0 \text{ in } M_0\backslash J^+_{g}(\supp(f)),
\end{gathered}
\eeq 
where $\supp(f)\subset V$. We assume that $H^{(j)}(x, z), j = 1, 2$ are genuinely nonlinear in $z$ on $I(p_-, p_+)$ with $p_\pm = \hat \mu(s_\pm), -1<s_-<s_+<1$. If the source-to-solution maps $L^{(j)}$ satisfy
\beq
L^{(1)} (f) = L^{(2)} (f)
\eeq
for all $f$  in a small neighborhood of the zero function in $C_0^4(V)$, then for $x \in I(p_-, p_+)$ we have 
\beq
\p^k_z H^{(1)}(x, 0) = \p^k_z H^{(2)}(x, 0) \text{ for } k\geq 4.
\eeq 
Moreover, if  $\p_z^3 H^{(1)}(x, 0) = \p_z^3 H^{(2)}(x, 0)$, we have
\beq
\begin{gathered}
\p_z^2 H^{(1)}(x, 0) = \p_z^2 H^{(2)}(x, 0).
\end{gathered}
\eeq
\end{theorem}

The theorem determines the coefficients of the Taylor expansion of $H$ in $z$ for $k\geq 4$ in general. 
We leave the determination of the cubic and quadratic terms as well as the linear term to future publications.

We remark that for all the inverse problems considered in this paper, the linear versions have not been solved yet. To solve these type of hyperbolic inverse problems, the boundary control (BC method) developed by Belishev has been used (see for example \cite{KKL}). The BC method depends on the unique continuation theorem of Tataru \cite{Ta1, Ta2} which assumes that the metric depends analytically on $t$. However, for globally hyperbolic Lorentzian manifolds, the coefficients of $\square_g$ are smooth in $t$ in general. Hence Tataru's theorem does not apply. See also Alinhac's counter-examples \cite{Al}. 

\subsection{Gauge invariance and non-linear Yamabe-type equations}
Now we discuss the gauge invariance and this is related to the nonlinear Yamabe equations. Let $R_g$ denote the scalar curvature of $g$. 
By \cite[Def.\ 3.5.9]{Bar}, the conformal wave operator (or the standard Yamabe operator multiplied by constant $\frac{n-1}{4n}=\frac 16$ with $n=3$)
\ba
\Y_gu(x)=\square_gu(x)+\frac16 R_g(x)u(x)
\ea
is conformally invariant in the following sense:
If $\varphi(x)$ is a  positive scalar function and  $\widetilde g_{jk}(x)=\varphi(x)^{p-2}g_{jk}(x)$ with
$p=\frac {2(n+1)}{n-1}$ ($p = 4$ when $n=3$), then
\beq
\Y_{\widetilde g}u=\varphi^{1-p} \Y_g(\varphi u).
\eeq
See also Theorem 5.1 of Appendix VI of \cite{Cb}. 
Let us consider the equation \eqref{eqsem}   of the form
\begin{eqnarray}\label{eq: wave-eq general 2}
&&\Y_{g}u(x) +H(x,u(x))=f(x), \quad\hbox{on }M_0,\hspace{-.5cm}\\ 
\nonumber &&\quad 
u(x)=0,\quad \hbox{for }x\in M_0\setminus J^+_{g}(\supp(f)).
\end{eqnarray}
We call the operator $\Y_{g}+H(x,\cdot)$ the non-linear Yamabe operator.
Note that $H$  does not contain linear terms.
 
Let us make now a gauge change:  We change metric $g_{jk}$ to the conformal metric 
\beq\label{eq gauge 1}
\widetilde  g_{jk}(x)=\varphi(x)^{p-2}g_{jk}(x).
\eeq
Then
\begin{eqnarray} \label{eq: wave-eq general 3}
&& \Y_{\widetilde  g}(\widetilde  u(x))+\widetilde  H(x,\widetilde  u(x))=\widetilde  f(x), 
\quad\hbox{on }M_0,\hspace{-.5cm}\\ \nonumber
&&\quad 
\widetilde  u(x)=0,\quad \hbox{for }x\in M_0\setminus J^+_{\widetilde  g}(\supp(\widetilde  f)),
\end{eqnarray}
where
\begin{eqnarray}\label{eq gauge 2}
\widetilde  u(x)&=& \varphi(x)^{-1} u(x),\\ \label{eq gauge 3}
\widetilde H(x,z) &=& \varphi(x)^{1-p} H(x, \varphi(x) z),
\\ \label{eq gauge 4}
\widetilde  f(x)&=& \varphi(x)^{1-p} f(x).
\end{eqnarray}

We consider now $\mathcal M_\varphi:u\mapsto \widetilde  u=\varphi(x)^{-1} u$ as a gauge transformation that
changes functions by the rules (\ref {eq gauge 1}) and  (\ref {eq gauge 2})-(\ref {eq gauge 4}).
Let $V\subset  M_0$ be sets where we do observations. Then the measurement map $L_{V;\Y_{g},H}:f\mapsto u|_V$ changes in the gauge  transformation as
$$
\mathcal M_{\varphi}L_{V;\Y_g,H} (\mathcal M_{\varphi^{p-1}})^{-1}=L_{V;\Y_{\tilde g},\tilde H}.
$$
Note that the measurement map depends on $V$, the metric $g$ and the nonlinear function $H$. We denote by $\Psi^*(\Y_g+H)$  the non-linear operator
$$
v\mapsto  \Y_{\Psi^*g}v+\Psi^*(H(\Psi(\,\cdotp),(\Psi^{-1})^*v)),
$$
that is obtained from $\Y_g+H(x,\,\cdotp)$ via a change of coordinates. Also, we denote by
$$
\mathcal M_{\varphi}(\Y_g+H)=\Y_{\tilde g}+\tilde H
$$
the gauge transformation of the non-linear operator $v\mapsto \Y_gv+H(\,\cdotp,v)$. Let 
$$
[\Y_g+H]=\{\mathcal M_{\varphi}(\Y_g+H) ;\ \varphi|_V=1\}
$$
denote   the class of operators, defined on the set $I_g(p_-, p_+)$, that are  gauge equivalence to  $\Y_g+H$. Then the measurement map does not change in the gauge transformation $M_{\varphi}$, that is,
\ba
L_{V,\Y_g,H}=L_{V;\Y_{\tilde g},\tilde H}.
\ea

Now we state a theorem of determining the metric and nonlinearity up to a gauge transformation.
\begin{theorem}\label{maingauge}
Let $(M^{(j)}, g^{(j)}), j = 1, 2$ be two $4$-dimensional globally hyperbolic Lorentzian manifolds. Let $\hat \mu^{(j)}(t) \subset M^{(j)}$ be time-like geodesics where $t\in [-1, 1]$ and $V^{(j)}\subset M^{(j)}$ be open relatively compact neighborhoods of $\hat \mu^{(j)}([s_-, s_+])$ where $-1<s_-<s_+<1$. Let $M_0^{(j)} = (-\infty, T_0)\times N^{(j)}, T_0 > 0$ such that $V^{(j)}\subset M_0^{(j)}$. Consider the nonlinear Yamabe equations with source terms
\beq
\begin{gathered}
\Y_{g^{(j)}} u(x) + H^{(j)}(x, u(x))  = f(x), \text{ on } M^{(j)}_0,\\
u = 0 \text{ in } M_0^{(j)}\backslash J^+_{g^{(j)}}(\supp(f)),
\end{gathered}
\eeq
where $\supp(f)\subset V^{(j)}$.  We assume that $H^{(j)}(x, z)$ are genuinely nonlinear on $I(p_-^{(j)}, p_+^{(j)})$  where $p^{(j)}_\pm = \hat \mu^{(j)}(s_\pm)$, the functions $z\mapsto H^{(j)}(x,z)$  are real-analytic and  $H^{(j)}(x, z) = O(z^4)$ as $z\rightarrow 0$. 

Suppose that there is a diffeomorphism $\Phi: V^{(1)}\rightarrow V^{(2)}$ such that $\Phi(p^{(1)}_\pm) = p^{(2)}_\pm$ and the source-to-solution maps $L^{(j)} = L_{V^{(j)};\Y_{g^{(j)}},H^{(j)}}$  satisfy
$$
L^{(2)}(f)=(\Phi^{-1})^*  (L^{(1)}(\Phi^*f))
$$
for all $f\in \mathcal W$, where $\mathcal W\subset C^4_0(V^{(2)})$ in a neighborhood of zero. Then there is a diffeomorphism $\Psi: I(p^{(1)}_-, p^{(1)}_+)\rightarrow I(p^{(2)}_-, p^{(2)}_+)$ such that
$$
\Psi^*(\Y_{g^{(2)}}+H^{(2)})\in [\Y_{g^{(1)}}+H^{(1)}],
$$
that is, the non-linear Yamabe operator $\Y_{g^{(2)}}+H^{(2)}$ is equal to
$\Y_{g^{(1)}}+H^{(1)}$ up to a combined gauge and coordinate transformation.
\end{theorem}

\subsection{Examples when the metric cannot be determined}\label{countex}
According to Theorem \ref{main1}, we know that for $H$ genuinely nonlinear, we can determine the conformal class of the metric. But in general this conformal factor cannot be determined even though the nonlinear functions are known, as demonstrated by the examples below. 

\textbf{Example 1}: Consider the conformal wave operator $\mcy_g = \square_g + \frac{1}{6}R_g$. In case when the scalar curvature vanishes, $\mcy_g = \square_g$. Consider the following nonlinear equation
\beqq\label{eqwave1}
\mcy_g u(x)+H(x, u)=f(x), \ \ H(x, z)= a(x) z^2.
\eeqq
Let $L_{g,a} f=u|_V$ be the source-to-solution map. We take $a(x) = (-\det g(x))^{-\frac{1}{4}} $ where $\det g$ denotes the determinant of the metric $g$. After the gauge transformation 
\beqq\label{gaugetrans}
\begin{gathered}
\tilde g_{jk}(x)=e^{2\gamma(x)} g_{jk}(x),\ \ \tilde u(x)=e^{-\gamma(x)} u(x), \ \ \tilde f(x)=e^{-3\gamma(x)} f(x),
\end{gathered}
\eeqq
we get $\det \tilde g(x) = e^{8\gamma(x)} \det g(x)$ and the equation \eqref{eqwave1} is transformed to
\beq
\begin{gathered}
\mcy_{\tilde g} \tilde u(x) +  e^{-3\gamma} e^{\gamma}(-\det \tilde g(x))^{\frac{1}{4}}  (e^{\gamma}\tilde u(x))^2= e^{-3\gamma}f(x),
\end{gathered}
\eeq
and we get
\beqq\label{eqwave2}
\mcy_{\tilde g} \tilde u(x) +  \tilde a(x) (\tilde u(x))^2= \tilde f(x),
\eeqq
where $\tilde a(x) = (-\det \tilde g(x))^{-\frac{1}{4}}$. Recall that on set $V$ where we perform the measurements, we have $\gamma = 0$, so that we get
\beq
\tilde u = u,\ \ \tilde f = f \text{ on } V.
\eeq
Thus we conclude that the two source-to-solution maps $L_{g, a} = L_{\tilde g, \tilde a}$ and this means we cannot determine two conformal metrics from the source-to-solution map in this case.

\textbf{Example 2}: Even in the case when two nonlinear functions are the same i.e.\ $H^{(1)} = H^{(2)}$ in Theorem \ref{main1}, we can still construct examples when the metric cannot be determined. Consider the following equation
\beqq\label{eqwave3}
\mcy_g u(x)+H(x, u)=f(x), \ \ H(x, z)=b z^3,
\eeqq
where $b$ can be any function of $x$. Let $L_{g, b}(f) = u|_V$ be the source-to-solution map. After the gauge transformation \eqref{gaugetrans}, the equation \eqref{eqwave3} is transformed to
\beq
\begin{gathered}
\mcy_{\tilde g} \tilde u(x) +  e^{-3\gamma}b (e^{\gamma}\tilde u(x))^3= e^{-3\gamma}f(x),
\end{gathered}
\eeq
and we get
\beqq\label{eqwave4}
\mcy_{\tilde g} \tilde u(x) +  b (\tilde u(x))^3= \tilde f(x).
\eeqq
Notice that the nonlinear terms in equations \eqref{eqwave3} and \eqref{eqwave4} are the same. Since $\gamma = 0$ on $V$, we get that
\beq
\tilde u = u,\ \ \tilde f = f \text{ on } V.
\eeq
Thus we conclude that the two source-to-solution maps $L_{g, b} = L_{\tilde g, b}.$ This means that when $H(x, z) = b z^3$, we cannot determine two conformal metrics from the source-to-solution map.

\subsection{Outline of the paper}
As in \cite{KLU}, we do not prove Theorem \ref{main1} by linearization but by producing artificial point sources thanks to the nonlinear interaction of linear waves. This is the reason we require $H$ to be genuinely nonlinear. However, compared with the analysis of singularities in \cite{KLU}, we carry out a more thorough microlocal analysis which enable us to characterize the type of the new singularities as well as to find their orders and principal symbols. The improvements are obtained from these new informations. We must mention that singularities due to nonlinear interactions in hyperbolic equations were actively studied in the 80's and 90's mainly for $1+2$ dimension by Bony \cite{Bo}, Melrose-Ritter \cite{MR, MR1}, Rauch-Reed \cite{RR}, etc. See Beals \cite{Bea} for an overview. However, in this work, we follow the idea in \cite{KLU, KLU1} to consider solutions of semilinear wave equations depending on some small parameters, and we analyze the singularities in the asymptotic expansion terms of the solution (instead of the solution itself). This simplifies our analysis and relates the singularities to the nonlinear term $H$. 
 
The paper is organized as following. In Section \ref{prel}, we collect some preliminaries from Lorentzian geometry and microlocal analysis. We derive the lower order asymptotic expansion of solutions to the semilinear wave equation. The main analysis lies in Section \ref{singu}, where we study the singularities produced by the nonlinear interaction of two, three and four conormal distributions. Here we carry out the analysis for a general set-up for the interaction of four linear waves. We start solving the inverse problem for lower order nonlinearities in Section \ref{secinv}, by making use of the analysis in Section \ref{singu} and following the approach of Kurylev-Lassas-Uhlmann \cite{KLU}. However, we encountered a problem that $H$ being cubic does not determine the metric by the analysis in Section \ref{singu}. To deal with this as well as higher order nonlinear terms, in Section \ref{singuhigh} we analyze in detail the singularities in higher order asymptotic expansions of the solutions. Finally, we prove the main theorems in Section \ref{pfmain}.

\section{Wave equations on Lorentzian manifolds}\label{prel}
\subsection{Lorentzian geometry}
We explain our assumptions on the Lorentzian manifold $(M, g)$ and introduce some notations. The general references are \cite{BEE, O, KLU}. 

Assume that $(M, g)$ is a $1+3$ dimensional Lorentzian manifold which is time oriented and globally hyperbolic. We take the signature of the metric as $(-, +, +, +)$. It is proved by Bernal and S\'anchez \cite{BS} that $(M, g)$ is globally hyperbolic if there is no closed causal paths in $M$ and for any $p, q\in M$ and $p<q$, the set $J(p, q)$ is compact. Also in \cite{BS0}, it is proved that $(M, g)$ is isometric to the product manifold $\mbr\times N$ with $g =-\beta(t, y)dt^2 + \kappa(t, y)$. Here $N$ is a $3$-dimensional manifold, $\beta: \mbr\times N\rightarrow \mbr_+$ is smooth and $\kappa$ is a Riemannian metric on $N$ and smooth in $t$. Without loss of generality, we identify $(M, g)$ with this isometric image. 
We shall use $x = (t, y) = (x^0, x^1, x^2, x^3)$ as the local coordinates on $M$. It is worth mentioning that for each $t\in \mbr$, the submanifold $\{t\}\times N$ is a Cauchy surface, i.e.\ every in-extendible causal curve intersects the submanifold only once. See for example \cite[Page 65]{BEE}. Besides the Lorentzian metric, we will take a complete Riemannian metric $g^+$ on $M$, whose existence is guaranteed by \cite{NO}. With this metric, we can introduce distances on $M$ and $TM$, and Sobolev spaces on $M$.  

For $p\in M$, we denote  the collection of light-like vectors at $p$ by $L_pM = \{\theta \in T_pM\backslash \{0\}: g(\theta, \theta) = 0\}$ and the bundle by $LM = \cup_{p\in M}L_pM$. The future (past) light-like vectors are denoted by $L^+_pM$ ($L^-_pM$), and the bundle $L^\pm M = \cup_{p\in M}L_p^\pm M$. Also, the set of light-like covectors at $p\in M$ is denoted by $L^*_pM$ and the bundles $L^{*}M, L^{*,\pm}(M)$ are defined similarly. Since the metric $g$ is non-degenerate, there is a natural isomorphism $i_p: T_pM\rightarrow T_p^*M$. For any $\theta \in L_pM$, $\xi = i_p(\theta)$ is in $L_p^*M$. With this isomorphism, we sometimes use vectors and co-vectors interchangeably. 
Let $\exp_p: T_pM  \rightarrow M$ be the exponential map. The geodesic from $p$ with initial direction $\theta$ is denoted by $\gamma_{p, \theta}(t) = \exp_p (t\theta), t\geq 0$. We denote the forward light-cone at $p\in M$ by
\beq
\mcl^+_pM = \{\gamma_{p, \theta}(t): \theta \in L^+_pM, t> 0\}
\eeq 
Notice that $p\notin \mcl^+_pM$ and $\mcl^+_p M$ is a subset of $M$.

\subsection{Lagrangian distributions}\label{seclagdis}
We will consider solutions to the wave equation with singularities on the conormal directions of a submanifold of $M$. Recall that the cotangent bundle $T^*M$ is a symplectic manifold with canonical two form given by $\omega = d\xi\wedge dx$ in local coordinates $(x, \xi)$. A submanifold $\La \subset T^*M$ is called Lagrangian if $\dim \La = 4$ and the canonical two form vanishes on $\La$. A simple example which we will consider later is the conormal bundle of a submanifold. Let $K\subset M$ be a submanifold. Then
\beq
N^*K = \{(x, \xi)\in T^*M\backslash 0: \langle \xi, \theta \rangle = 0, \theta \in T_xK\},
\eeq
the conormal bundle of $K$ is a Lagrangian submanifold of $T^*M\backslash 0$. Here $0$ represents the zero section of $T^*M$. We review the basic facts of conormal distributions and paired Lagrangian distributions. Our main references are \cite{DUV, Ho3, MU, GrU93, GrU90}.

\subsubsection{Conormal distributions} Let $X$ be a $n$-dimensional smooth manifold and $\La$ be a smooth conic Lagrangian submanifold of $T^*X\backslash 0$. Following the standard notation, we denote by $I^\mu(\La)$ the space of Lagrangian distributions of order $\mu$ associated with $\La$. In particular, for $U$ open in $X$, let $\phi(x, \xi): U\times \mbr^N \rightarrow \mbr$ be a smooth non-degenerate phase function (homogeneous of degree $1$ in $\xi$) that locally parametrizes $\La$ i.e.\ 
\beq
\{(x, d_x\phi)\in T^*_UX\backslash 0: x\in U, d_\xi \phi = 0\} \subset \La.
\eeq
Here $0$ denotes the zero-section. Then $u\in I^{\mu}(\La)$ can be locally written as a finite sum of oscillatory integrals
\beq
\int_{\mbr^N} e^{i\phi(x, \xi)} a(x, \xi) d\xi, \ \ a\in S^{\mu + \fnf - \frac{N}{2}}(U\times \mbr^N),
\eeq
where $S^\bullet(\bullet)$ denotes the standard symbol class, see \cite[Section 18.1]{Ho3}.  For $u\in I^\mu(\La)$, we know that the wave front set $\WF(u)\subset \La$ and $u\in H^s(X)$ for any $s< -\mu-\fnf$. The distribution $u$ has a principal symbol $\sigma(u)$ defined invariantly on $\La$, see \cite[Section 25.1]{Ho4}. 

For a submanifold $Y\subset X$ of codimension $k$,  the conormal bundle $N^*Y$ is a Lagrangian submanifold. We denote $I^{\mu}(Y) = I^\mu(N^*Y)$, which are called conormal distributions to $Y$. We remark that in our notation, $\mu$ is always the order of the distribution (instead of the order of the symbol which is used in \cite{GrU90} etc). If we take $u$ as a distributional half-density on $M$, the principal symbol of $u$  is well-defined in $S^{\mu + \fnf}(\La; \Omega^\ha)/S^{\mu + \fnf- 1}(\La; \Omega^\ha)$, where $\Omega^\ha$ denotes the half-density bundle on $\La$. 
Actually, we can find local coordinates $x = (x', x''), x'\in \mbr^k, x''\in \mbr^{n-k}$ such that $Y = \{x' = 0\}$. Let $\xi = (\xi', \xi'')$ be the dual variable, then $N^*Y = \{x' = 0, \xi'' = 0, \xi'\neq0 \}$. We can write $u\in I^\mu(Y)$ as
\beq
u = \int_{\mbr^k} e^{ix'\xi'} a(x'', \xi')d\xi', \ \ a\in S^{\mu + \fnf - \frac{k}{2}}(\mbr^{n-k}_{x''}; \mbr^k_{\xi'}).
\eeq
In this case, the principal symbol is 
\beq
\sigma(u) = (2\pi)^{\fnf - \frac{k}{2}} a_0(x'', \xi')|dx''|^\ha |d\xi'|^\ha,
\eeq
where  $a_0 \in S^{\mu + \fnf - \frac{k}{2}}(\mbr^{n-k}_{x''}; \mbr^k_{\xi'})$ is such that $a - a_0 \in  S^{\mu + \fnf - \frac{k}{2}- 1}(\mbr^{n-k}_{x''}; \mbr^k_{\xi'})$. See for example \cite[Section 18.2]{Ho3}. Later, we also use the notation $\sigma_{N^*Y}(u)$ to emphasize where the symbol is defined. We remark that for Lorentzian manifold $(M, g)$, there is a natural choice of the density bundle $d\vol_g$. Thus the half-density bundles can be trivialized and we shall consider the principal symbols of distribution $u$ as functions on $T^*M$.

\subsubsection{Paired Lagrangian distributions}\label{secplag} For two Lagrangians $\La_0, \La_1 \subset T^*X\backslash 0$ intersecting cleanly at a codimension $k$ submanifold i.e.\ 
\beq
T_p\La_0\cap T_p\La_1  = T_p(\La_0\cap \La_1),\ \ \forall p\in \La_0\cap \La_1,
\eeq
the paired Lagrangian distribution associated with $(\La_0, \La_1)$ is denoted by $I^{p, l}(\La_0, \La_1)$. For $u\in I^{p, l}(\La_0, \La_1)$, we know that $\WF(u)\subset \La_0\cup \La_1$. Also, microlocally away from the intersection $\La_0\cap \La_1$, $u \in I^{p+l}(\La_0\backslash \La_1)$ and $u \in I^p(\La_1\backslash \La_0)$ so $u$ has well-defined principal symbols $\sigma_{\La_0}(u)$ and $\sigma_{\La_1}(u)$ on the corresponding Lagrangians. They also satisfy the compatibility condition on $\La_0\cap \La_1$, see \cite{GU} and \cite{MU}.

Since all cleanly intersecting pairs of Lagrangians are locally equivalent (see \cite[Prop.\ 2.1]{GU}), we can write down such distributions explicitly as oscillatory integrals in certain model pairs of Lagrangians. For example (see (5.14) of \cite{DUV}), we let $x = (x', x'', x''')\in \mbr^k\times \mbr^{n-k-d}\times \mbr^d$ and take
\beqq\label{eqmodelag}
\tilde \La_0 = N^*\{x' = x''  = 0\}, \ \ \tilde \La_1 = N^*\{x''  = 0\}
\eeqq
as the model pair. For $u \in I^{p, l}(\tilde\La_0, \tilde\La_1)$, we can write 
\beqq\label{pairlag}
u = \int_{\mbr^n} e^{i(x'\cdot \xi' + x''\cdot \xi'')} b(x''', \xi', \xi'')d\xi'd\xi'' + u_0,
\eeqq
where $u_0 \in I^p(\tilde \La_1)$ and $b \in S^{M, M'}(\mbr_{x'''}^d; \mbr_{\xi''}^{n-k-d}; \mbr^k_{\xi'})$, $M = p-\fnf + \frac{k}{2} + \frac{d}{2}, M' = l - \frac{k}{2}$ is a symbol of product type i.e.\ for $x'''$ in a compact set $K$ and for multi-indices $\alpha, \beta, \gamma$, there is $C_{K, \alpha, \beta, \gamma}>0$ such that
\beq
|\p_{x'''}^\alpha \p_{\xi'}^\beta \p_{\xi''}^\gamma b(x''', \xi', \xi'')| \leq C_{K, \alpha, \beta, \gamma}\langle \xi', \xi''\rangle^{M-|\beta|}\langle \xi''\rangle^{M'-|\gamma|}.
\eeq
Away from $\tilde\La_0\cap \tilde\La_1$, the principal symbol of $u$ on $\tilde\La_1$ is given in $(5.16)$ of \cite{DUV}, which is
\beq
\sigma_{\tilde\La_1}(u) = (2\pi)^{\frac{3n+2k-2d}{4}} (\mcf')^{-1} b|dx'|^\ha |d\xi''|^\ha |dx'''|^\ha 
\eeq
in $S^{M}(\mbr_{x'''}^d; \mbr_{\xi''}^{n-k-d}; I^{l-\frac{k}{4}}(\mbr_{x'}^k; N^*\{0\}))$ modulo terms of order $M-1$. Here $\mcf'$ denotes the partial Fourier transform in $x'$, $(\mcf')^{-1}$ denotes the inverse transform in $\xi'$ variable  and the symbol space $S^M(\cdot)$ is the standard symbol space $S^{M}(\mbr_{x'''}^d; \mbr_{\xi''}^{n-k-d})$ but with distributional values in $I^{l - \frac{k}{4}}(\mbr^k_{x'}; N^*\{0\})$. We refer to (5.16) of \cite{DUV}  for other equivalent descriptions. On $\tilde\La_0\backslash \tilde \La_1$, $u$ is a Lagrangian distribution with amplitude $b\in S^{M + M'}(\mbr_{x'''}^d; \mbr_{\xi', \xi''}^{n-d})$ and we find that
\beq
\sigma_{\tilde\La_0}(u) = (2\pi)^{\fnf - \frac{n-d}{2}} b_0(x''', \xi', \xi'')|dx'''|^\ha|d\xi'd\xi''|^\ha,
\eeq
where $b_0\in S^{M + M'}(\mbr_{x'''}^d; \mbr_{\xi', \xi''}^{n-d})$ is such that $b - b_0 \in S^{M + M' - 1}(\mbr_{x'''}^d; \mbr_{\xi', \xi''}^{n-d})$.

\subsection{Linear wave equations and causal inverses}\label{seclinwa}
Let $\square_g$ be the Laplace-Beltrami operator on $(M, g)$. In local coordinates, we have 
\beq
\square_g = (-\det g(x))^{-\ha} \sum_{i, j = 0}^3\frac{\p}{\p x^i} ((-\det g(x))^\ha g^{ij}(x) \frac{\p}{\p x^j}).
\eeq
We consider the  linear wave equation
\beqq\label{eqlin}
\begin{split}
\square_g v  &= f \text{ on } M_0,\\
v &= 0 \text{ on } M_0\backslash J^+(\supp(f)).
\end{split}
\eeqq
where $f$ is a source term compactly supported in $t\geq 0$ and will be specified later.  The Schwartz kernel of the causal inverse of the wave operator $\square_g$ is a paired Lagrangian distribution we now review. We remark that later we do not distinguish the notations of operators and their Schwartz kernels unless it is necessary. 

We let $\mcp(x, \xi) = |\xi|^2_{g^*}$ be the principal symbol of $\square_g$, which is also the dual  metric function on $T^*M$. Here $g^* = g^{-1}$ denotes the dual Lorentzian metric on $T^*M$. Let $\Sigma_g$ be the characteristic set i.e.\ 
\beq
\Sigma_g = \{(x, \xi)\in T^*M: \mcp(x, \xi) = 0\}.
\eeq
Note $\Sigma_g$ consists of light-like co-vectors. The Hamilton vector field of $\mcp$ is denoted by $H_\mcp$ and in local coordinates
\beq
H_\mcp = \sum_{i = 0}^3( \frac{\p \mcp}{\p \xi_i}\frac{\p }{\p x^i} - \frac{\p \mcp}{\p x^i}\frac{\p }{\p \xi_i}).
\eeq
The integral curves of $H_\mcp$ in $\Sigma_g$ are called null bicharacteristics and their projections to $M$ are geodesics.

Consider the product manifold $M\times M$ and the cotangent bundle $T^*M\times T^*M$. Let $\pi$ be the projection to the left factor. We can regard $\mcp, \Sigma_g, H_\mcp$ as objects on product manifolds by pulling them back using $\pi$. Let $\diag = \{(z, z')\in M\times M: z = z'\}$ be the diagonal and 
\beq
N^*\diag = \{(z, \zeta, z', \zeta')\in T^*(M\times M)\backslash 0: z = z', \zeta' = -\zeta\}
\eeq 
be the conormal bundle of $\diag$ (minus the zero section). Then we let $\La_g$ be the Lagrangian submanifold obtained by flowing out $N^*\diag\cap \Sigma_g$ under $H_\mcp$. It is proved in \cite{MU} (see also \cite{DUV}) that $\square_g$ has a parametrix $Q_g \in I^{-\frac{3}{2}, -\ha}(N^*\diag, \La_g)$. In particular, this means that for any $f\in \mcd'(M)$ the distribution space on $M$, we have $\square_g Q_g f = f + f_0$ where $f_0 \in C^\infty(M)$. 

On globally hyperbolic manifolds, the wave operator $\square_g$ (and more generally normally hyperbolic operators) has a unique causal inverse which we denote by $\square_g^{-1}$, see for example \cite[Theorem 3.3.1]{Bar} and \cite{Fr}. Then   $\square_g^{-1} - Q_g$ is a smoothing operator, and $\square_g^{-1}\in I^{-\frac{3}{2}, -\ha}(N^*\diag, \La_g)$. For convenience, we take $Q_g = \square_g^{-1}$ as the causal inverse. 

If the source $f$ is a Lagrangian distribution, we can describe the solution $v$ easily by the following proposition, which is essentially Prop. 2.1 of \cite{GrU93}.
\begin{prop}\label{winv}
Suppose $\La_0\subset T^*M\backslash 0$ is a conic Lagrangian intersecting $\Sigma_g$ transversally and such that each bicharacteristics of $\mcp$ intersect $\La_0$ a finite number of times. Then 
\beq
Q_g : I^\mu(\La_0) \rightarrow I^{-\frac{3}{2} + \mu, -\ha}(\La_0, \La_0^g),
\eeq
where $\La_0^g$ denotes the Lagrangian submanifold obtained from flowing out of $\La_0\cap \Sigma_g$ under the Hamiltonian flow. Furthermore, for $(x, \xi)\in \La_0^g\backslash \La_0$,
\beq
\sigma(Q_g u)(x, \xi) = \sum \sigma(Q_g)(x, \xi; y_j, \eta_j)\sigma(u)(y_j, \eta_j),
\eeq
where the summation is over the points $(y_j, \eta_j)\in \La_0$ which lie on the bicharacteristics from $(x, \xi)$.  
\end{prop}

Recall that we can define Sobolev spaces on $M$ using the Riemannian metric $g^+$. From Prop.\ 5.6 of \cite{DUV} or Theorem 3.3 of \cite{GrU90}, we also have
\begin{prop}\label{propqg}
For $m\in \mbr$, $Q_g: H^{m}_{\comp}(M)\rightarrow H^{m+1}_{\loc}(M)$ is continuous.  
\end{prop}

\subsection{Semilinear wave equations and the asymptotic analysis}\label{semiwave}
Consider the semi-linear wave equation
\beqq\label{eqmain}
\begin{gathered}
\square_g u + H(x, u)  = f(x) \text{ on } M_0 = (-\infty, T_0)\times N, T_0 > 0\\
u = 0 \text{ in } M_0\backslash J^+(\supp(f)),
\end{gathered}
\eeqq
where $H$ is smooth and $f$ is a source term supported in $t\geq 0$ to be specified later. Here we write $x = (x^0, x^1, x^2, x^3)$ where $x^0 = t$. The local well-posedness of \eqref{eqmain} has been analyzed in \cite{KLU}, see also Section 3.1.2 and Appendix B of \cite{KLU1}. In particular, let $B\subset N$ be compact, $T_0>0$ and $f\in C^r([0,T_0], H^{s}(B))\cap C^{r+1}([0, T_0], H^{s-1}(B)), r\geq 0$. Recall that the function spaces are defined using the Riemannian metric $g^+$. If $r+ s\geq 4$ is even and $f$ is small enough, there is a unique solution $u$ to \eqref{eqmain} such that $u\in C_0^r([0,T_0], H_0^{s}(N))\cap C_0^{r+1}([0, T_0], H_0^{s-1}(N))$ and 
\beq
\|u\|_{C_0^r([0,T_0], H_0^{s}(N))\cap C_0^{r+1}([0, T_0], H_0^{s-1}(N))} \leq C \|f\|_{C^r([0,T_0], H^{s}(B))\cap C^{r+1}([0, T_0], H^{s-1}(B))},
\eeq
for some constant $C>0$, see equation (27) of \cite{KLU1}. Hereafter, $C$ denotes a generic constant.  It is convenient to use Sobolev spaces on $M$. Let $m = s+ r$, we know that $u\in H^{m}(M_0)$ and $\|u\|_{H^{m}(M_0)}\leq C\|f\|_{H^{m}(M_0)}$ for $m\geq 4$ even. We remark that the regularity required in the local well-posedness results may not be optimal but this is not our main concern.   

Next, we carry out the asymptotic analysis of $u$. We are able to compute the first few asymptotic terms explicitly, which are sufficient for many purposes. So we first consider the case when
\beqq\label{eqhlow}
H(x, z) = a(x)z^2 + b(x)z^3 + c(x)z^4,
\eeqq
for $z\in \mbr$ sufficiently small and $a, b, c$ are smooth functions in $x$. Indeed, higher order nonlinear terms will not affect the first four terms in the asymptotic expansion. Later in Section 5, we will return to the general case and use different techniques to analyze the higher order asymptotic terms.

We assume that $f_i \in C_0^4(M_0)\subset H_0^{4}(M_0), i = 1, 2, 3, 4$. Notice that $H^{4}(M_0)$ is an algebra since $\text{dim} M = 4$, see e.g.\ \cite[Theorem 8.3.1]{Ho}. Let $\eps_i > 0, i = 1, 2, 3, 4$ be small parameters and 
\beq
f = \sum_{i = 1}^4 \eps_i f_i \in C_0^{4}(M_0)\subset H_0^4(M_0)
\eeq
be the source term in \eqref{eqmain}. Here we used Sobolev embedding. Then the linear equation $\square_g v = f$ has a solution $v = \sum_{i = 1}^4 \eps_i v_i \in H^{5}(M_0)$ where $\square_g v_i = f_i$. The semilinear equation \eqref{eqmain} has a solution $u$ satisfying
\beq
\|u\|_{H^{4}(M_0)} \leq C(\eps_1+ \eps_2+\eps_3+\eps_4) \sum_{i = 1}^4\|f_i\|_{H^{4}(M_0)}.
\eeq
Now we derive the asymptotic expansion of $u$ as $\eps_i\rightarrow0$. Using that
\beq
\square_g (u - v) + H(x, u) = 0,
\eeq
we have
\beqq\label{eqite}
u = v - Q_g(H(x, u)) = v - Q_g(au^2 + bu^3 + cu^4).
\eeqq
We substitute $u$ back in the right hand side of \eqref{eqite}. We first compute
\beqq\label{eqite1}
\begin{split}
&u^2 = v^2 -2vQ_g(au^2) - 2vQ_g(bu^3) + Q_g(au^2)Q_g(au^2) + \mcr, \\
&u^3 = v^3 -3v^2Q_g(au^2) + \mcr,\\
&u^4 = v^4 + \mcr. 
\end{split}
\eeqq
Here $\mcr$ denotes the collection of terms which are $o(\eps_1\eps_2\eps_3\eps_4)$ in $H^{4}(M_0)$. Actually, since $H^{4}(M_0)$ is an algebra, $u^5\in H^{4}(M_0)$. Applying $Q_g$, we know that $Q_g(u^5)\in H^{5}(M_0)$. In general, $Q_g$ increases the regularity by $1$ thus the terms in $\mcr$ are all in $H^{4}(M_0)$. Now we have
\beqq\label{exp3}
\begin{gathered}
u = v - Q_g(av^2) + 2Q_g(avQ_g(au^2)) + 2Q_g(avQ_g(bu^3)) - Q_g(aQ_g(au^2)Q_g(au^2)) \\
- Q_g(bv^3) + 3Q_g(bv^2Q_g(au^2)) - Q_g(cv^4) + \mcr.
\end{gathered}
\eeqq
Finally, by substituting \eqref{eqite1} into \eqref{exp3}, we obtain
\beqq\label{exp4}
\begin{gathered}
u = v - Q_g(av^2) + 2Q_g(avQ_g(av^2)) - 4Q_g(avQ_g(avQ_g(av^2))) - Q_g(aQ_g(av^2)Q_g(av^2)) \\
+ 2Q_g(avQ_g(bv^3)) - Q_g(bv^3) + 3Q_g(bv^2Q_g(av^2)) - Q_g(cv^4) + \mcr.
\end{gathered}
\eeqq
This is the asymptotic expansion of the solution $u$. Indeed, we only need to consider the following terms where $i, j, k, l$ are distinct.
\beqq\label{eqinter}
\begin{split}
\mcu^{(2)}_{ij} &= -Q_g(a v_i v_j);\\
\mcu^{(3)}_{ijk} &=  2Q_g(av_iQ_g(av_j v_k)) - Q_g(bv_iv_jv_k);\\
\mcu^{(4)}_{ijkl} &= - 4Q_g(av_iQ_g(av_jQ_g(av_kv_l))) - Q_g(aQ_g(av_iv_j)Q_g(av_kv_l))\\
&+ 2Q_g(av_iQ_g(bv_jv_kv_l)) + 3Q_g(bv_iv_jQ_g(av_kv_l)) - Q_g(cv_iv_jv_kv_l).
\end{split}
\eeqq
Then we can write $u$ as
\beqq\label{uasym}
u = v + \sum_{i, j = 1, i\neq j}^4 \eps_i\eps_j \mcu^{(2)}_{ij} + \sum_{i, j, k = 1, i\neq j\neq k}^4  \eps_i\eps_j\eps_k \mcu^{(3)}_{ijk} + \eps_1\eps_2\eps_3\eps_4 \sum_{i, j, k, l = 1, i\neq j\neq k\neq l}^4 \mcu^{(4)}_{ijkl} + \mcr.
\eeqq
For convenience, we shall denote
\beqq\label{interms}
\begin{split}
\mcu^{(2)} &= \sum_{i, j = 1, i\neq j}^4 \mcu^{(2)}_{ij} = \sum_{i, j = 1, i\neq j}^4 \p_{\eps_i}\p_{\eps_j} u|_{\{\eps_i = \eps_j = 0\}}, \\
 \mcu^{(3)} &= \sum_{i, j, k = 1, i\neq j\neq k}^4 \mcu^{(3)}_{ijk} = \sum_{i, j, k = 1, i\neq j\neq k}^4 \p_{\eps_i}\p_{\eps_j}\p_{\eps_k}u|_{\{\eps_i = \eps_j = \eps_k = 0\}}, \\
  \mcu^{(4)} &= \sum_{i, j, k, l = 1, i\neq j\neq k\neq l}^4 \mcu^{(4)}_{ijkl} = \p_{\eps_1}\p_{\eps_2}\p_{\eps_3}\p_{\eps_4}u|_{\{\eps_1 = \eps_2 = \eps_3 = \eps_4 = 0\}}.
\end{split}
\eeqq

We remark that $\mcu^{(L)}, L = 2, 3, 4$ are the terms involving the interaction of $L$ conormal waves, i.e.\ they involve multiplications of $L$ conormal distributions. From their expressions, it is clear where the coefficients of the Taylor expansion of $H(x, u)$ contributed to the interactions. This will be important for the inverse problem.

\section{Analysis of the singularities in the nonlinear interactions}\label{singu}
In this section, our goal is to understand the singularities coming from the terms in \eqref{interms}. Because these terms also appear in the analysis of other nonlinear equations, for example Einstein equations studied in \cite{KLU1}, and their analyses are similar, we will make some general assumptions on the conormal distributions $v_i$ in Section \ref{assump}. Later in Section \ref{secinv}, we will construct concrete $v_i$ which satisfy the general assumptions. 

The key part in the analysis is to understand the multiplications of several conormal distributions and paired Lagrangian distributions. The multiplication of two such distributions are analyzed in Greenleaf and Uhlmann \cite{GrU93}. Our main interest is the new singularities of $\mcu^{(L)}, L\geq 3$ which involves the multiplication of more than three such distributions. In particular, we will characterize the type of these distributions and find their principal symbols. These results are used to solve the inverse problem in Section \ref{secinv}.

\subsection{Assumptions and notations}\label{assump}
Recall that two submanifolds $X, Y$ of $M$ intersect transversally if 
\beq
T_q X + T_q Y = T_q M, \ \ \forall q\in X\cap Y.
\eeq

We make the following definition on the intersection of four submanifolds.
\begin{definition}\label{trans4}
Assume that $K_i, i = 1, 2, 3, 4$ are codimension $1$ submanifolds of $M$ such that $N^*K_i$ i.e.\ the co-vectors normal to $K_i$ are light-like. We say that $K_i$ intersect transversally if the following are satisfied.
\begin{enumerate}
\item $K_i, K_j, i< j$ intersect transversally at $K_{ij}$, which is a codimension $2$ submanifold of M;
\item $K_i, K_j, K_k, i<j<k$ intersect transversally at  $K_{ijk}$, which is a codimension $3$ submanifold of M;
\item $K_i, i = 1, 2, 3, 4$ intersect transversally at a point $q_0$. 
\end{enumerate}
\end{definition}
In particular, the last condition means that the four submanifolds intersect at a point $q_0$ and the normal co-vectors $\zeta_i$ to $K_i$ at $q_0$ are linearly independent. We remark that for any $q\in M$, we can find $K_i$ intersecting transversally at $q$.
For $i = 1, 2, 3, 4$, we shall denote 
\beqq\label{lakj}
\begin{gathered}
\La_i = N^*K_i; \ \ \La_{ij} = N^*K_{ij}, i< j; \ \ \La_{ijk} = N^*K_{ijk}, i<j<k.
\end{gathered}
\eeqq
All of these are Lagrangian submanifolds in $L^{*}M$.  At $q_0$, we let $\La_{q_0} = T_{q_0}^*M\backslash 0$ which is a conic Lagrangian submanifold. 
We will use the following notations 
\beqq\label{eqlak}
\begin{gathered}
\La^{(1)} = \cup_{i = 1}^4 \La_i; \ \ \La^{(2)} = \cup_{i, j = 1, i<j}^4 \La_{ij}; \ \ \La^{(3)} = \cup_{i, j, k = 1, i<j<k}^4\La_{ijk};\\
K^{(1)} = \cup_{i = 1}^4 K_i; \ \ K^{(2)} = \cup_{i, j = 1, i<j }^4 K_{ij}; \ \ K^{(3)} = \cup_{i, j, k = 1, i<j<k }^4 K_{ijk}.
\end{gathered}
\eeqq

Now we consider the normal form of four transversally intersecting Lagrangians near $q_0$. This is convenient for local computations. In $\mbr^4$, we take $\tilde K_i = \{x^i = 0\}, i = 1, 2, 3, 4$. Let $\xi = (\xi_1, \xi_2, \xi_3, \xi_4)$  be the dual variables to $x$ in $T^*\mbr^4$. Then $\tilde \La_i = N^*\tilde K_i \subset T^*M\backslash 0$ can be expressed as 
\beqq\label{eqtildela}
\begin{gathered}
\tilde \La_1 =  \{x^1 = 0, \xi_2 = \xi_3 = \xi_4 = 0, \xi_1\neq 0\},\ \ \tilde \La_2 = \{x^2 = 0, \xi_1 = \xi_3 = \xi_4 = 0, \xi_2\neq 0\},\\ 
\tilde \La_3 =  \{x^3 = 0, \xi_1 = \xi_2 = \xi_4 = 0, \xi_3\neq 0\},\ \ \tilde \La_4 = \{x^4 = 0, \xi_1 = \xi_2 = \xi_3 = 0, \xi_4\neq 0\}.
\end{gathered}
\eeqq

\begin{lemma}\label{symp}
Let $K_i, i = 1, 2, 3, 4$ intersect transversally at $q_0$. Then there exists a neighborhood $\mco$ of $q_0$ and diffeomorphism $\phi : \mco\rightarrow \mbr^4$ such that $\phi(q_0) = 0$ and $\phi(K_i) \subset \tilde K_i, i = 1, 2, 3, 4$. Therefore, $\phi$ induces a symplectomorphism $\chi: T^*_\mco M\rightarrow T^*\mbr^4$ such that $\chi(\La_i)\subset \tilde \La_i, i = 1, 2, 3, 4$.
\end{lemma}

\bpf
Since $K_i$ intersect transversally at $q_0$, we can find (in local coordinate patches) smooth functions $f_1, f_2, f_3, f_4$ such that $K_i = \{f_i = 0\}$. See e.g.\ \cite[Appendix\ C.3]{Ho3}. Moreover, the differentials $df_i$ are linearly independent at $T_q^*M$.  By the inverse mapping theorem, we can find local coordinates $x = (x^1, x^2, x^3, x^4)$ in a coordinate patch $(\mco, \phi)$ of $q_0$ such that $\phi(K_i)\subset \{x^i = 0\}, i = 1, 2, 3, 4.$ This finishes the proof. 
\epf

We introduce some notations used throughout the rest of the paper. With the Riemannian metric $g^+,$ we can define the unit cotangent bundle $S^*M$. Let $\eps>0$ be a small parameter. For any set $\Gamma$ in $T^*M\backslash 0$, we denote by $\Gamma(\eps)$  a conic neighborhood of $\Gamma$ such that $\Gamma(\eps)\cap S^*M$ is an $\eps$ neighborhood of $\Gamma \cap S^*M$. In particular, $\Gamma(\eps)$ tends to the closure of $\Gamma$ as $\eps\rightarrow 0$. Also, for any conic set $\Gamma \subset T^*M$, we use the standard notation $\mcd'(M; \Gamma)$ to denote distributions $u$ with $\WF(u)\subset \Gamma.$  

Recall the Lagrangian submanifold $\La_g$ in Section \ref{seclinwa}. For any $\Gamma\subset T^*M$,  we denote the flow out of $\Gamma$ under $\La_g$ by
\beqq\label{flowset}
\Gamma^g = \La_g' \circ (\Gamma\cap \Sigma_g).
\eeqq
Here as usual in microlocal analysis, for $\La \subset T^*M\times T^*M$, 
\beq
\La' = \{(x, \xi, y, \eta) \in (T^*M\backslash 0)\times (T^*M\backslash 0): (x, \xi, y, -\eta)\in (T^*M\backslash 0)\times (T^*M\backslash 0)\}.
\eeq
 
Assume that $K_i, i = 1, 2, 3, 4$ are codimension $1$ submanifolds of $M$ such that $N^*K_i$
Finally, for $K_i, i = 1, 2, 3, 4$ intersecting transversally, we let $v_i \in I^\mu(K_i)$ be supported in $t > 0$. Later in Section 4, we will see how to construct $v_i$ from compactly supported conormal distributions $f_i$ so that 
$v_i = Q_g(f_i)$ and the asymptotic expansion \eqref{uasym} in Section \ref{semiwave} holds.  


\subsection{Singularities in two waves interacting}\label{double}
We study the following term in $\mcu^{(2)}$
\beq
\mcu^{(2)}_{12} = -Q_g(av_1v_2).
\eeq
The analysis works for the other terms $\mcu^{(2)}_{ij}$ in $\mcu^{(2)}$. This term involves the multiplication of two conormal distributions and the application of a paired Lagrangian distribution $Q_g$. First consider the multiplication. The following is essentially Lemma 1.1 of \cite{GrU93} (see also \cite{DUV}). We briefly repeat the proof to find the symbols. 
\begin{lemma}\label{prod2}
For $\La_j$ defined in \eqref{lakj}, let  $u\in I^\mu(\La_1), v\in I^{\mu'}(\La_2)$. Then we can write $w = uv$ as
\beq
\begin{gathered}
w = w_1 + w_2, \ \ w_1 \in I^{\mu, \mu' + \frac{n}{4}}(\La_{12}, \La_1), \ \ \WF(w_1)\cap \La_2 = \emptyset, \\
 w_2\in I^{\mu', \mu+ \frac{n}{4}}(\La_{12}, \La_2), \ \ \WF(w_2) \cap \La_1 = \emptyset.
 \end{gathered} 
\eeq
Moreover, for any $(q, \zeta)\in \La_{12}\backslash(\La_1\cup \La_2)$, we can write $\zeta = \zeta_1 + \zeta_2$ in a unique way such that $\zeta_i\in N^*_qK_i, i = 1, 2$. Microlocally away from $\La_1\cup \La_2$, $uv \in I^{\mu+\mu'+\frac{n}{4}}(\La_{12})$ and the principal symbol of $uv$  satisfies
\beq
\sigma_{\La_{12}}(uv)(q, \zeta) = (2\pi)^{-1}\sigma_{\La_1}(u)(q, \zeta_1)\cdotp\sigma_{\La_2}(v)(q, \zeta_2). 
\eeq 
\end{lemma}

We remark that in this lemma (as well as the rest of the paper), we fix a choice of the density bundle on $(M, g)$ to trivialize the half-density factors in distributions and principal symbols. 

\bpf[Proof of Lemma \ref{prod2}]
For any $q\in K_{12}$, we can choose local coordinates $x$ such that $K_1 = \{x^1 = 0\}$ and $K_2 = \{x^2 = 0\}$. Let $\xi$ be the dual variable to $x$. We can write $u, v$ as oscillatory integrals.
\beq
u = \int_{\mbr} e^{i x^1\xi_1}A(x^2, x^3, x^4; \xi_1) d\xi_1, \ \ v =  \int_{\mbr} e^{i x^2\xi_2}B(x^1, x^3, x^4; \xi_2) d\xi_2,
\eeq
where $A\in S^{\mu+ \fnf-\ha}(\mbr^3; \mbr), B\in S^{\mu'+\fnf-\ha}(\mbr^3; \mbr)$. Thus 
\beqq\label{equv}
uv = \int_{\mbr^2} e^{i (x^1\xi_1+x^2\xi_2)}A(x^2, x^3, x^4; \xi_1)B(x^1, x^3, x^4; \xi_2) d\xi_1d\xi_2.
\eeqq
By introducing cut-off functions as in the proof of Lemma 1.1 in \cite{GrU93}, we obtain the first statement. 

Next, we notice that $N_q^*K_{12}$ is spanned by $N_q^*K_1$ and $N_q^*K_2$ and the vectors are linearly independent at $q$ by the transversality assumption. Thus any $\zeta\in \La_{12}$ can be written as a unique linear sum of $\zeta_1 \in N_q^*K_1$ and  $\zeta_2 \in N_q^*K_2$. The statement about the principal symbol follows from \eqref{equv} by  a  stationary phase argument. Indeed, away from $\La_1\cup \La_2$, we know that $\xi_1, \xi_2\neq0$. By  stationary phase, we find that  
\beq
\begin{gathered}
\sigma(uv)(q, \zeta) =  (2\pi)^{\fnf - 1}A(0, x^3, x^4; \xi_1)B(0, x^3, x^4; \xi_2),\\
\text{where } \sigma(u) = (2\pi)^{\fnf - \frac{1}{2}}A(0, x^3, x^4; \xi_1),\ \ \sigma(v) = (2\pi)^{\fnf - \frac{1}{2}}B(0, x^3, x^4; \xi_1),
\end{gathered}
\eeq
modulo lower order terms. This finishes the proof.
\epf

Next, we consider the action of $Q_g$ on paired Lagrangian distributions, see \cite[Section 2 ]{GrU93}. 
\begin{lemma}\label{inter2}
Let $u\in  I^{p, l}(\La_{12}, \La_1)$ and $\WF(u)\cap \La_2  = \emptyset$. For $Q_g\in I^{-\frac{3}{2}, -\ha}(N^*\diag, \La_g)$, we have
\beq
Q_gu \in I^{p-1, l-1}(\La_{12}, \La_1). 
\eeq
Moreover, let $(q, \zeta)\in \La_{12}\backslash(\La_1\cup \La_2)$, the principal symbol $\sigma_{\La_{12}}(Q_gu)(q, \zeta) = |\zeta|_{g^*(q)}^{-2} \sigma_{\La_{12}}(u)(q, \zeta)$, where $g^*$ is the dual Lorentzian metric to $g$. 
\end{lemma}
\bpf
Let $\zeta_1\in N^*_qK_1$ and $\zeta_2\in N^*_qK_2$ where $q\in K_{12}$. Then $\zeta_1, \zeta_2$ are linearly independent light like co-vectors and they span $N^*_q(K_{12})$. If $\zeta$ is a linear combination of $\zeta_1, \zeta_2$ and $\zeta\in \Sigma_g$ i.e.\ light-like, then $\zeta$ is proportional to $\zeta_1$ or $\zeta_2$. Therefore, $\La'_g\circ ((\La_{12}\backslash\La_2)\cap \Sigma_g) \subset \La_{1}$. Now we apply Proposition 2.2 of \cite{GrU93}  to get the first statement. 

The principal symbol on $\La_{12}$ away from $\La_1$  can actually be found in the proof of Prop.\ 2.2 and Prop.\ 2.1 of \cite{GrU93}. We give the proof below for completeness. By microlocalizing and conjugating by an elliptic Fourier integral operator, we can assume that $M = \mbr^n, n = 4$ with local coordinates $x = (x^0, x')$, $\La_{12} = T_0^*\mbr^n\backslash 0$ and $\Sigma_g = \{(x, \xi): \xi_0 = 0\}$. Therefore, 
\beq
\La_1 = \{(x^0, 0; 0, \xi'): x^0\in \mbr, \xi'\in \mbr^{n-1}, \xi'\neq 0\}. 
\eeq
In this model pair, we can write $u$ as 
\beq
u(x) = \int_{\mbr^n} e^{ix\xi} A(x; \xi', \xi_0)d\xi, \ \ A\in S^{p-\fnf + \ha, l-\ha}(\mbr^n_x; \mbr_{\xi'}^{n-1}; \mbr_{\xi_0}).
\eeq
Also, we can write $Q_g$ as
\beq
Q_g w(x) = \int_{\mbr^n\times \mbr^n} e^{i(x-y)\eta} B(x, y; \eta', \eta_0) w(y) dy d\eta, \ \ B \in S^{-\frac{3}{2}  + \ha, -\ha -\ha}(\mbr^{2n}; \mbr_{\eta'}^{n-1}; \mbr_{\eta_0}).
\eeq
Note that $Q_g$ on $N^*\diag\backslash \La_g$ is a pseudo-differential operator and the principal symbol is $|\zeta|_{g^*(q)}^2$ times the half-density factor. Then we have  
\beq
\begin{gathered}
Q_g u(x) = \int_{\mbr^n\times \mbr^n\times \mbr^n} e^{i(x-y)\eta} B(x, y; \eta', \eta_0) e^{iy\xi} A(y; \xi', \xi_0)d\xi dy d\eta = \int_{\mbr^n} e^{ix\eta} C(x, \eta)d\eta, \\
\text{where } C(x, \eta) = \int_{\mbr^n\times \mbr^n} e^{iy(\xi-\eta)} B(x, y; \eta', \eta_0) A(y; \xi', \xi_0)d\xi dy.
\end{gathered}
\eeq
On $\La_{12}\backslash\La_1$ where $x = 0, \eta_0 \neq 0$, the principal symbol of $Q_g(u)$ is $C(0; \eta)$ and this can be found by the stationary phase lemma as
\beq
C(0; \eta) = B(0, 0; \eta)A(0; \eta),
\eeq
modulo lower order terms. To finish the proof, we just need to observe that $B(0, 0; \eta) = |\eta|_{g^*(0)}^{-2}$ is the principal symbol of $Q_g$ on $\La_{12}$ and $\sigma_{\La_{12}}(u) = A(0; \eta)$. 
\epf

Applying the above results, we obtain that  
\beq
\mcu_{ij}^{(2)} \in  I^{\mu-1, \mu}(\La_{ij}, \La_i) +  I^{\mu-1, \mu}(\La_{ij}, \La_j).
\eeq
Hence $\WF(\mcu_{ij}^{(2)})\subset \La_{ij}\cup \La_i\cup\La_j$, and the singular support of $\mcu^{(2)}$ is contained in $K^{(1)}$. Therefore, the interaction of two conormal waves does not produce new propagating singularities.

\subsection{Singularities in three waves interacting}\label{trip}
Next, we analyze the term $\mcu^{(3)}$. This term is not analyzed carefully in \cite{KLU}. It suffices to study the following term 
\beqq\label{ite3}
\begin{gathered}
\mcu^{(3)}_{321} = 2\mcv_1 - \mcv_2, \text{ where }\\
\mcv_{1} =  Q_g(av_3Q_g(av_2v_1)), \ \ \mcv_2 = Q_g(bv_1v_2v_3),
\end{gathered}
\eeqq
because the other terms in $\mcu^{(3)}$ are similar. 
The key point is to understand the multiplication of conormal distributions and paired Lagrangian distributions. The result is a new type of distribution associated with three intersecting Lagrangians. We do not have a convenient theory for such distributions at hand. However, since our major concern is the new singularities produced in the interaction, we can avoid the difficulty by cutting off the product distribution away from the old singularities. The reason we separate the two terms in $\mcu_{321}^{(3)}$ is that we will show  the new singularities in $\mcv_2$ are stronger than those in $\mcv_1$. Heuristically, we expect such result because the term $\mcv_2$ is due to the stronger nonlinearities in $H(x, u)$ while $\mcv_1$ is obtained by iterating the lower order nonlinearities. In addition, $Q_g$ increases the regularity by $1$.

For the three wave interactions, it suffices to assume there are three submanifolds $K_i, i = 1, 2, 3$ intersecting transversally, meaning 
\begin{itemize}
\item $K_i, K_j, i<j$ intersect transversally at $K_{ij}$ which is a codimension $2$ submanifold of $M$;
\item $K_1, K_2, K_3$ intersect transversally at $K_{123}$ which is a codimension $3$ submanifold of $M$. 
\end{itemize}
We can find the normal form near $K_{123}$ as in Lemma \ref{symp}.
\begin{lemma}\label{symp1}
Let $\tilde K_i = \{x^i = 0\} \subset \mbr^4, i = 1, 2, 3$ and $\tilde \La_i$ be as defined in \eqref{eqtildela}. For any $q\in K_{123}$, there exists a neighborhood $\mco$ of $q$ and diffeomorphism $\phi : \mco\rightarrow \mbr^4$ such that $\phi(q) = 0$ and $\phi(K_i) \subset \tilde K_i, i = 1, 2, 3$. Therefore, $\phi$ induces a symplectomorphism $\chi: T^*\mco  \rightarrow T^*\mbr^4$ such that $\chi(\La_i)\subset \tilde \La_i, i = 1, 2, 3$.
\end{lemma}


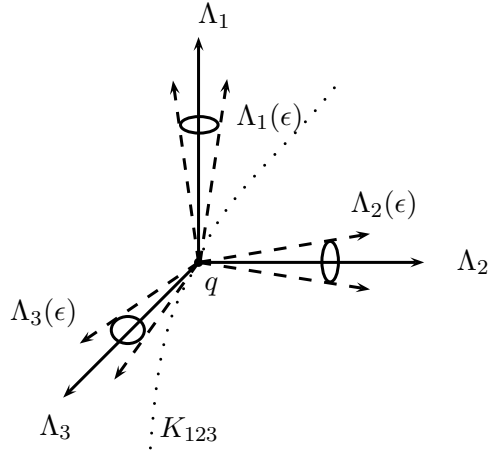
\begin{figure}[htbp]
\centering
\scalebox{1} 
{
\begin{pspicture}(-2,-3.1)(6.70291,3.1)
\psline[linewidth=0.04cm,arrowsize=0.05291667cm 2.0,arrowlength=1.4,arrowinset=0.4]{<-}(1.9210156,2.4196484)(1.9210156,-0.5803516)
\psline[linewidth=0.04cm,arrowsize=0.05291667cm 2.0,arrowlength=1.4,arrowinset=0.4]{->}(1.9210156,-0.5803516)(4.9210157,-0.5803516)
\psline[linewidth=0.04cm,arrowsize=0.05291667cm 2.0,arrowlength=1.4,arrowinset=0.4]{->}(1.9210156,-0.5803516)(0.12101562,-2.3803515)
\psellipse[linewidth=0.04,dimen=outer](1.9310156,1.2496485)(0.27,0.13)
\psellipse[linewidth=0.04,dimen=outer](3.6710157,-0.57035154)(0.13,0.29)
\psellipse[linewidth=0.04,dimen=outer](0.9810156,-1.4803516)(0.24,0.2)
\psline[linewidth=0.04cm,linestyle=dashed,dash=0.16cm 0.16cm,arrowsize=0.05291667cm 2.0,arrowlength=1.4,arrowinset=0.4]{->}(1.9410156,-0.5803516)(2.3010156,1.8596485)
\psline[linewidth=0.04cm,linestyle=dashed,dash=0.16cm 0.16cm,arrowsize=0.05291667cm 2.0,arrowlength=1.4,arrowinset=0.4]{->}(1.9210156,-0.5803516)(1.5810156,1.8396485)
\psline[linewidth=0.04cm,linestyle=dashed,dash=0.16cm 0.16cm,arrowsize=0.05291667cm 2.0,arrowlength=1.4,arrowinset=0.4]{->}(1.9210156,-0.5803516)(4.2010155,-0.20035157)
\psline[linewidth=0.04cm,linestyle=dashed,dash=0.16cm 0.16cm,arrowsize=0.05291667cm 2.0,arrowlength=1.4,arrowinset=0.4]{->}(1.9210156,-0.5803516)(4.2010155,-0.94035155)
\psline[linewidth=0.04cm,linestyle=dashed,dash=0.16cm 0.16cm,arrowsize=0.05291667cm 2.0,arrowlength=1.4,arrowinset=0.4]{->}(1.9210156,-0.5803516)(0.34101564,-1.6603515)
\psline[linewidth=0.04cm,linestyle=dashed,dash=0.16cm 0.16cm,arrowsize=0.05291667cm 2.0,arrowlength=1.4,arrowinset=0.4]{->}(1.9210156,-0.5803516)(0.8010156,-2.1403515)
\usefont{T1}{ptm}{m}{n}
\rput(2.1324707,2.7046484){$\La_1$}
\usefont{T1}{ptm}{m}{n}
\rput(5.572471,-0.57535156){$\La_2$}
\usefont{T1}{ptm}{m}{n}
\rput(0.0224707,-2.7753515){$\La_3$}
\usefont{T1}{ptm}{m}{n}
\rput(4.3924706,0.18464845){$\La_2(\eps)$}
\usefont{T1}{ptm}{m}{n}
\rput(2.8724708,1.3446485){$\La_1(\eps)$}
\usefont{T1}{ptm}{m}{n}
\rput(-0.1324707,-1.2553515){$\La_3(\eps)$}
\psdots[dotsize=0.12](1.9210156,-0.5803516)
\psbezier[linewidth=0.04,linestyle=dotted,dotsep=0.16cm](1.2810156,-3.049199)(1.3810157,-1.7691993)(1.5010157,-1.2291992)(1.9210156,-0.5891992)(2.1810157,0.29080078)(3.0810156,0.8708008)(3.8410156,1.8908008)
\usefont{T1}{ptm}{m}{n}
\rput(2.0924706,-0.9){$q$}
\usefont{T1}{ptm}{m}{n}
\rput(1.8,-2.8){$K_{123}$}
\end{pspicture} 
}
\caption{Explanation of the cut-off on $\La_{123}$. The dotted curve is $K_{123}$ which is a codimension $3$ submanifold. The picture is the cotangent space over $q\in K_{123}$. The cut-off functions $\psi_i$ are supported in $\La_i(\eps)$. }
\label{figcut}
\end{figure}

We consider the multiplication of a conormal distribution and paired Lagrangian distribution.
\begin{lemma}\label{prod3}
Assume that $u\in I^\mu(\La_3), v\in I^{p, l}(\La_{12}, \La_1)$ are compactly supported near $K_{123}$. For  $\eps>0$ sufficiently small, we can write $w = uv$ as
\beqq\label{eq3}
\begin{gathered}
   w = w_0 + w_1 + w_2, \ \ w_0\in I^{\mu+p+l - \ha}(\La_{123}), \ \ w_1\in \mcd'(M; \La^{(1)}(\eps)),\\
   w_2 \in I^{p, \mu+\fnf}(\La_{13}, \La_1) + I^{\mu, p + \fnf}(\La_{13}, \La_3) + I^{p, l}(\La_{12}, \La_1).
   \end{gathered}
\eeqq
Moreover, for $q\in K_{123}$ and $\zeta \in N_q^*K_{123}\backslash(\cup_{i = 1}^3\La_i)$,  we can write $\zeta = \zeta_1+\zeta_2+\zeta_3$ uniquely for $\zeta_i\in N_q^*K_i$. The principal symbol of $w_0$ satisfies
\beq
\sigma_{\La_{123}}(w_0)(q, \zeta) = (2\pi)^{-1}\sigma_{\La_3}(u)(q, \zeta_3)\cdotp\sigma_{\La_{12}}(v)(q, \zeta_1+\zeta_2).
\eeq
\end{lemma}

\bpf
By choosing local coordinates as in Lemma \ref{symp}, it suffices to consider the distributions on model Lagrangians $\tilde \La_i$. Also, we notice that away from $K_{12}$, the distribution in $I^{p, l}(\La_{12}, \La_1)$ can be written as a sum of  distributions in $I^{p}(\La_1)$. Its multiplication with $v_3$ gives terms in $w_2$ by Lemma \ref{prod2}. So we will focus on the part near $\La_{12}\cap \La_3$.

For $u \in I^\mu(\tilde\La_3)$, we can write
\beq
u = \int_{\mbr} e^{ix^3\xi_3}A(x^1, x^2, x^4; \xi_3)d\xi_3,
\eeq
where $A\in S^{m}(\mbr_x^3; \mbr_{\xi_3})$ with $m = \mu + \fnf - \ha$. Also, for $v\in I^{p, l}(\La_{12}, \La_1)$, we can write (modulo a distribution in $I^p(\tilde \La_1)$) that
\beq
v = \int_{\mbr^2} e^{i(x^1 \xi_1 + x^2\xi_2)} B(x^3, x^4; \xi_1, \xi_2)d\xi_1d\xi_2,
\eeq
with $B\in S^{M, M'}(\mbr^2_{x^3, x^4}; \mbr_{\xi_2}, \mbr_{\xi_1})$ being a symbol of product type where $M = p-\fnf + \ha + 1, M' = l-1$, see Section \ref{secplag} and the references there. Therefore, the product can be written as
\beq
w =  \int_{\mbr^3} e^{i(x^1 \xi_1 + x^2\xi_2 + x^3\xi_3)}A(x^1, x^2, x^4; \xi_3)B(x^3, x^4; \xi_1, \xi_2)d\xi_1d\xi_2d\xi_3.
\eeq

Let $\hat\psi: \mbr\rightarrow [0, \infty)$ be a smooth cut off function such that $\hat\psi(t) = 1$ for $|t|<\ha$ and $\hat\psi(t) = 0$ for $|t|>1$. For $\eps>0$, we let 
\beq
\psi_i(\xi) = \hat\psi(\frac{\sum_{i = j}^4|\xi_j| - |\xi_i|}{\eps \sum_{j = 1}^4|\xi_j|}),\ \ i = 1, 2, 3.
\eeq
Then we see that each $\psi_i$ is supported in an $\eps$-neighborhood of $\La_i$. See Figure \ref{figcut}. Without loss of generality, we can assume that this neighborhood is $\La_i(\eps)$. Now we let $\psi(x, \xi) =  \prod_{i = 1}^3\psi_i(\xi)$ be supported in the union of $\La_i(\eps), i = 1, 2, 3,$ which is contained in $\La^{(1)}(\eps)$ with $\La^{(1)}$ defined in \eqref{eqlak}. We consider
\beq
w_0 = \int_{\mbr^3} e^{i(x^1 \xi_1 + x^2\xi_2 + x^3\xi_3)} (1 - \psi) A(x^1, x^2, x^3; \xi_3) B(x^3, x^4; \xi_1, \xi_2)d\xi_1d\xi_2d\xi_3.
\eeq
By the symbol estimates of $A$ and $B$, we know that $(1-\psi)AB\in S^{M+M'+m}(\mbr_x^4; \mbr^3_{\xi_1, \xi_2, \xi_3})$. Thus $w_0\in I^{\mu+p+l -\fnf + \frac{3}{2}}(\La_{123}\backslash \La^{(1)}(\eps/2)) \subset I^{\mu+p+l -\fnf + \frac{3}{2}}(\La_{123})$ for any $\eps>0$. 
For the other part, 
\beq
w_1 = \int_{\mbr^3} e^{i(x^1 \xi_1 + x^2\xi_2 + x^3\xi_3)} \psi(\xi) A(x^1, x^2, x^4; \xi_3)B(x^3, x^4; \xi_1, \xi_2)d\xi_1d\xi_2d\xi_3.
\eeq
we can use the standard stationary phase method to conclude that $\WF(w_1)\subset \La^{(1)}(\eps)$. 

Finally, we consider the symbols. This is similar to the proof of Lemma \ref{prod2}. We obtain using the stationary phase lemma that 
\beq
\begin{gathered}
\sigma(w_0) =  (2\pi)^{\fnf - \frac{3}{2}}A(0, 0, x^4; \xi_3)B(0, x^4; \xi_1, \xi_2),\\
\text{ where } \sigma(u) = (2\pi)^{\fnf - \ha}A(0, 0, x^4; \xi_3), \ \ \sigma(v) =  (2\pi)^{\fnf - 1}B(0, x^4; \xi_1, \xi_2),
\end{gathered}
\eeq
modulo lower order terms. This proves the relation of the symbols. 
\epf

We remark that from the proof we actually have that $\WF(w_1)$ is contained in the union of $\La^{(1)}$ and $\La^{(1)}(\eps)\cap \La_{123}$. This means away from the intersection $K_{123}$, $\WF(w_1)\subset \La^{(1)}$, and the wave front of $w_1$ at $K_{123}$ is in an $\eps$-neighborhood of $\La^{(1)}$. For our analysis, we do not need such a precise statement.

Finally, consider the triple interaction terms. 
\begin{prop}\label{inter3}
Let $\La_{123}^g$ be the flow out of $\La_{123}$ under $\La_g$, see \eqref{flowset}. Away from the union $\cup_{i = 1}^3 \La_i$, we have  
\beq
\mcv_1\in I^{3\mu - 3, -\ha}(\La_{123}, \La_{123}^g), \ \ \mcv_2\in I^{3\mu -  \frac{3}{2}, -\ha}(\La_{123}, \La_{123}^g).
\eeq
In particular, away from $\La_{123}$, we have $\mcv_1\in I^{3\mu - 3}(\La_{123}^g)$ and $\mcv_2\in I^{3\mu -  \frac{3}{2}}(\La_{123}^g).$
Moreover, away from $\cup_{i = 1}^3 \La_i$ and $\La_{123}$, we have
\begin{enumerate}
\item If $b$ is non-vanishing on $K_{123}$, then $\mcu^{(3)}_{321} \in I^{3\mu -  \frac{3}{2}}(\La_{123}^g)$.
\item If $b$ vanishes in a neighborhood of $K_{123}$ where $a$ is non-vanishing, then $\mcu^{(3)}_{321} \in I^{3\mu -  3}(\La_{123}^g)$.
\end{enumerate}
Here $a, b$ are the coefficients in the nonlinear function $H$, see \eqref{eqhlow}.
\end{prop}

\bpf
We start with $\mcv_1$. From the previous subsection, we know that 
\beq
Q_g(av_2v_1) \in I^{\mu-1, \mu+\fnf-1}(\La_{12}, \La_1) + I^{\mu-1, \mu+\fnf-1}(\La_{12}, \La_2).
\eeq
Here $n = 4$. By Lemma \ref{prod3}, for any $\eps>0$, we write $w = av_3Q_g(av_2v_1)  = w_0 + w_1 + w_2$ such that $w_0\in I^{3\mu -\frac{3}{2}}(\La_{123})$, $\WF(w_1)\subset \La^{(1)}(\eps)$ and $w_2$ is a sum of paired Lagrangian distribution with $\WF(w_2)\subset \La^{(1)}\cup\La^{(2)}$. Finally, we apply Prop.\ 2.1 of \cite{GrU93} to get that $Q_g(w_0) \in I^{3\mu-\frac{3}{2} -\frac{3}{2}, -\ha}(\La_{123}, \La^g_{123})$. For $w_2$, we can apply Lemma \ref{inter2} to see that $\WF(Q_gw_2) \subset \La^{(1)}\cup \La^{(2)}$.  For $w_1$, we apply the standard calculus of wave front sets, e.g.\ \cite[Corollary\ 1.3.8 ]{Du} to get 
\beq
\WF(Q_g(w_1)) \subset \La^{(1)}(\eps)\cup  (\La'_g\circ(\La^{(1)}(\eps)\cap \Sigma_g)).
\eeq 
The right hand side is a small neighborhood of $\La^{(1)}$ and tends to $\La^{(1)}$ as $\eps\rightarrow 0$. Here we used the fact that the Hamiltonian flow is a smooth map. Therefore, away from $\La^{(1)}$, $\mcv_1\in I^{3\mu - 3}(\La_{123}^g\backslash\La_{123})$. Next, for $\mcv_2$, the analysis is the same and the only difference is the order. We know that $v_2v_1\in I^{\mu, \mu+\fnf}(\La_{12}, \La_1) + I^{\mu, \mu+\fnf}(\La_{12}, \La_2)$. By Lemma \ref{prod3}, we obtain that $bv_3v_2v_1 \in I^{3\mu + 1}(\La_{123})$ modulo a distribution whose wave front set is contained in $\La^{(1)}(\eps)$. So after applying $Q_g$, we know that $\mcv_2\in I^{3\mu - \frac{3}{2}}(\La_{123}^g\backslash \La_{123})$ modulo a distribution whose wave front set is in $\La^{(1)}(\eps)$ for any $\eps$ small. Finally, if $b$ vanishes in a neighborhood of $K_{123}$, the term $\mcv_2$ is smooth. Then $\mcv_1 \in I^{3\mu - 3, -\ha}(\La_{123}, \La_{123}^g)$ carries the leading singularities. This completes the proof of the Proposition. 
\epf

We remark that since $K_{123}$ is a one-dimensional submanifold, the terms $\mcv_1, \mcv_2$ have conic singularities along $K_{123}$ and the singular support is contained in the projection of $\La_{123}^g$ to $M$.  To  see that these singularities are non-trivial and $\mcv_2$ actually has a stronger singularity than $\mcv_1$, we will compute their principal symbols. 

First we use Lemma \ref{prod2}, \ref{inter2} and \ref{prod3}. For any $q \in K_{123}$ and $\zeta\in N_q^*K_{123}$, we can write $\zeta = \sum_{i = 1}^3\zeta_i$ where $\zeta_i\in N^*_qK_i$. Also, we let $A_i$ be the principal symbols of $v_i$ (recall that the half-densities on $\La_i$ are trivialized). Then 
\beq
\sigma_{\La_{123}}(\mcv_2)(q, \zeta) = (2\pi)^{-2} b(q)A_1(q, \zeta_1)A_2(q, \zeta_2)A_3(q, \zeta_3).
\eeq
Here and after, we shall ignore the $2\pi$ factors in the symbol computations. Now let $\sigma_{\La_g}(Q_g)$ be the principal symbol of $Q_g$ on $\La_g$ away from $N^*\diag$. Since on globally hyperbolic manifolds there is no closed causal curve,  we can apply Prop.\ \ref{winv} to get 
 \beq
 \sigma_{\La_{123}^g}(\mcv_2)(y, \eta) =  (2\pi)^{-2}\sigma_{\La_g}(Q_g)(y, \eta, q, \zeta)b(q)A_1(q, \zeta_1)A_2(q, \zeta_2)A_3(q, \zeta_3),
 \eeq
 where $(y, \eta)$ is joined with $(q, \zeta)$ by bicharacteristics of $\square_g$. Similarly, one can show that
 \beq
  \sigma_{\La_{123}^g}(\mcv_1)(y, \eta) =  (2\pi)^{-2}\sigma_{\La_g}(Q_g)(y, \eta, q, \zeta)a^2(q) A_3(q, \zeta_3) |\zeta_1+\zeta_2|_{g^*(q)}^{-2} A_1(q, \zeta_1)A_2(q, \zeta_2).
 \eeq
 We remark that since $\sigma_{\La_g}(Q_g)$ is an invertible matrix, if the symbols $A_i$ are non-vanishing and $a$ or $b$ is non-vanishing, the principal symbols of $\mcv_1$ or $\mcv_2$ hence $\mcu^{(3)}_{321}$ will be non-vanishing on $\La_{123}^g$. In other words, by properly choosing $v_i \in I^\mu(K_i)$, the term $\mcu^{(3)}_{321}$ has a non-vanishing conic type singularity at $\La^{g}_{123}$.  
For the inverse problem, it would be interesting to know the set $\La_{123}^g$ because that is where we can measure the singularities. In this case, we can use only three waves and the analysis is simpler.

\subsection{Singularities in four waves interacting}\label{quad}
Among the terms in $\mcu^{(4)}$, we analyze $\mcu^{(4)}_{1234}$ as a model term which is given by
\beqq\label{u1234}
\begin{split}
\mcu^{(4)}_{1234} = &- 4Q_g(av_1Q_g(av_2Q_g(av_3v_4))) - Q_g(aQ_g(av_1v_2)Q_g(av_3v_4))\\
&+ 2Q_g(av_1Q_g(bv_2v_3v_4)) + 3Q_g(bv_1v_2Q_g(av_3v_4)) - Q_g(cv_1v_2v_3v_4),
\end{split}
\eeqq
where $a, b, c$ are the coefficients in the nonlinear function $H$ in \eqref{eqhlow}. We observe that all these terms involve two kind of basic operations. One is the product of two paired Lagrangian distributions such as $I^{*, *}(\La_{12}, \La_1)$ and $I^{*, *}(\La_{34}, \La_3)$, and the other one is the multiplication of a conormal distribution and the three wave interaction term we analyzed in Prop.\ \ref{inter3}. These multiplications result in a distribution associated with four intersecting Lagrangians. We again use cut-off techniques to prove that the new singularities are conormal to $\mcl^+_{q_0}$. Also, we show that the term $Q_g(cv_1v_2v_3v_4)$ produces the strongest singularity if $c(q_0)\neq 0$. This can be used to simplify the analysis for more complicated nonlinear equations. 

We first consider the multiplication of two paired Lagrangian distributions. This is similar to Lemma \ref{prod3}. 
\begin{lemma}\label{prod41}
Let $u\in  I^{p, l}(\La_{12}, \La_1), v\in  I^{p', l'}(\La_{34}, \La_3)$. For $\eps>0$ sufficiently small, we can write $w = uv$ as
\beq
\begin{gathered}
w = w_0 + w_1 + w_2, \ \ w_0\in I^{p+l+p'+l'+1}(\La_{q_0}), \ \ w_1 \in \mcd'(M; \La^{(1)}(\eps)),\\
w_2\in I^{p+ l + p'+\ha}(\La_{123})+I^{p+ l' + p'+\ha}(\La_{134})+I^{p, p'+1}(\La_{13}, \La_1) + I^{p', p+1}(\La_{13}, \La_3)\\
+ I^{p, l}(\La_{12}, \La_1) + I^{p', l'}(\La_{34}, \La_3).
\end{gathered}
\eeq
Moreover, for $\zeta\in \La_{q_0}\backslash \La^{(1)}$,  we can write $\zeta = \sum_{i = 1}^4\zeta_i$ uniquely for $\zeta_i\in N_{q_0}^*K_i$. The principal symbol of $w_0$ satisfies
\beq
\sigma_{\La_{q_0}}(w_0)(q_0, \zeta) = (2\pi)^{-1}\sigma_{\La_{12}}(u)(q_0, \zeta_1+\zeta_2)\cdotp\sigma_{\La_{34}}(v)(q_0, \zeta_3+\zeta_4).
\eeq
\end{lemma}
\bpf
It suffices to consider the distributions on the model Lagrangians $\tilde \La_i$. For $u \in I^{p, l}(\tilde \La_{12}, \tilde\La_1)$, we can write 
\beq
u = u_0+u_1, \ \ u_0 = \int_{\mbr^2} e^{i(x^1 \xi_1 + x^2\xi_2)} A(x^3, x^4; \xi_1, \xi_2)d\xi_1d\xi_2,\ \ u_1\in I^p(\tilde \La_1),
\eeq
with $A\in S^{M, M'}(\mbr^2_{x^3, x^4}; \mbr_{\xi_2}, \mbr_{\xi_1})$ being a symbol of product type where $M = p-\fnf + \ha + 1, M' = l-1$. Also, for $v\in I^{p', l'}(\La_{34}, \La_3)$, we can write 
\beq
v = v_0+v_1, \ \ v_0 = \int_{\mbr^2} e^{i(x^3 \xi_3 + x^4\xi_4)} B(x^1, x^2; \xi_3, \xi_4)d\xi_3d\xi_4, \ \ v_1 \in I^{p'}(\tilde \La_3),
\eeq
with $B\in S^{m, m'}(\mbr^2_{x^1, x^2}; \mbr_{\xi_3}, \mbr_{\xi_4})$ being a symbol of product type where $m = p'-\fnf + \ha + 1, m' = l'-1$.
The product of $u_0v_1$ and $u_1v_0$ and $v_1u_1$ are already studied. Actually, 
\beq
\begin{gathered}
u_0v_1 \in I^{p+ l + p'+\ha}(\La_{123}) + \mcd'(M; \La^{(1)}(\eps)),\\
 u_1v_0 \in I^{p+ l' + p'+\ha}(\La_{134}) + \mcd'(M; \La^{(1)}(\eps)),\\
 v_1u_1 \in I^{p, p'+1}(\La_{13}, \La_1) + I^{p', p+1}(\La_{13}, \La_3).
 \end{gathered}
\eeq
Thus their wave front sets are known and we only need to find the product
\beq
 u_0v_0 = \int_{\mbr^4} e^{i(x^1 \xi_1 + x^2\xi_2 + x^3\xi_3 + x^4\xi_4)} A(x^1, x^2; \xi_3, \xi_4)B(x^3, x^4; \xi_1, \xi_2)d\xi.
\eeq
We take the cut-off function $\hat\psi$ in Lemma \ref{prod3}. For $\eps>0$, we let 
\beq
\tilde\psi_i(\xi) = \hat\psi(\frac{\sum_{j = 1}^4|\xi_j| - |\xi_i|}{\eps\sum_{j = 1}^4|\xi_j|}),\ \ i = 1, 2, 3, 4.
\eeq
Then we see that each $\tilde\psi_i$ is supported in $\La_i(\eps)$, an $\eps$-neighborhood of $\La_i$. Now we let $\tilde\psi(x, \xi) = \sum_{i = 1}^4\tilde\psi_i(\xi)$ be supported in $\La^{(1)}(\eps)$. We consider
\beq
w_0 = \int_{\mbr^4} e^{i(x^1 \xi_1 + x^2\xi_2 + x^3\xi_3 + x^4\xi_4)} (1 - \tilde\psi) A(x^1, x^2; \xi_3, \xi_4) B(x^3, x^4; \xi_1, \xi_2)d\xi.
\eeq
Then by the symbol estimates of $A$ and $B$, we know that $(1-\psi)AB\in S^{M+M'+m+m'}(\mbr_x^4; \mbr^4_{\xi})$. Thus $w_0\in I^{p'+l'+p+l +1}(\La_{q_0}\backslash \La^{(1)}(\eps/2))$. For the other part, 
\beq
w_1 = \int_{\mbr^4} e^{i(x^1 \xi_1 + x^2\xi_2 + x^3\xi_3 + x^4\xi_4)} \tilde \psi(\xi) A(x^1, x^2; \xi_3, \xi_4) B(x^3, x^4; \xi_1, \xi_2)d\xi.
\eeq
We conclude that $\WF(w_1)\subset \La^{(1)}(\eps)$. The symbols can be found as in the proof of Lemma \ref{prod3}.
\epf

Next, we consider the multiplication of $I^*(\La_4)$ and $I^{*, *}(\La_{123}, \La_{123}^g)$. We need a lemma to decompose paired Lagrangian distributions. 
\begin{lemma}\label{prod43}
Let $\La_0, \La_1$ be two transversally (or more generally, cleanly) intersecting Lagrangians on $T^*M$. For $\eps>0$ sufficiently small, we can write $u\in I^{p, l}(\La_0, \La_1)$ as
\beq
u = u_0 + u_1, \ \ u_0 \in I^{p+l}(\La_0), \ \ u_1\in \mcd'(M; \La_1(\eps)).
\eeq
\end{lemma}
\bpf
From Section \ref{seclagdis}, we know that $I^{p, l}(\La_0, \La_1)\subset I^{p+l}(\La_0)$ away from $\La_1$. 
\epf

We remark that in general we can not decompose the paired Lagrangians as the sum of Lagrangian distributions, as explained in  \cite[Section 5]{DUV}. Our decomposition involves a Lagrangian distribution and another distribution with known wave front set.

Using Lemma \ref{prod2}, \ref{prod43} or by repeating the proof of Lemma \ref{prod41}, we have the following result whose proof is omitted.  
\begin{lemma}\label{prod42}
Let $u\in  I^{\mu}(\La_4), v \in  I^{\mu'}(\La_{123})$. For $\eps>0$ sufficiently small, we can write  $w = uv$ as
\beq
w = w_0 + w_1, \ \ w_0\in I^{\mu+\mu'+1}(\La_{q_0}), \ \ w_1 \in \mcd'(M; \La_4(\eps)\cup\La_{123}(\eps)).
\eeq
Moreover, for $\zeta\in \La_{q_0} \backslash \La^{(1)}$,  we can write $\zeta = \sum_{i = 1}^4\zeta_i$ uniquely for $\zeta_i\in N_{q_0}^*K_i$. The principal symbol of $w_0$ satisfies
\beq
\sigma_{\La_{q_0}}(w_0)(q_0, \zeta) = (2\pi)^{-1}\sigma_{\La_{4}}(u)(q_0, \zeta_4)\cdotp\sigma_{\La_{123}}(v)(q_0, \zeta_1+\zeta_2+\zeta_3).
\eeq
\end{lemma}

Finally, we use the above lemmas to analyze $\mcu^{(4)}_{1234}$.
\begin{prop}\label{inter4}
Let $v_i\in I^\mu(\La_i), i = 1, 2, 3, 4$ and $\La^g_{q_0}$ be the flow out of $\La_{q_0}$ under the Hamiltonian flow, see \eqref{flowset}. Let $\Xi \doteq \La^{(1)}\cup \La^{(3), g} \cup \La_{q_0}$, we have the following conclusions for $\mcu^{(4)}_{1234}$.
\begin{enumerate}
\item If $c(q_0)\neq 0$, we have $\mcu^{(4)}_{1234}\in I^{4\mu+\frac{3}{2}}(\La^g_{q_0}\backslash \Xi)$;
\item If $c = 0$ in a neighborhood of $q_0$, where $b, a$ are non-vanishing, we have $\mcu^{(4)}_{1234}\in I^{4\mu-\ha}(\La^g_{q_0}\backslash \Xi)$;
\item If $b = c = 0$ in a neighborhood of $q_0$ where $a$ is non-vanishing, we have $\mcu^{(4)}_{1234}\in I^{4\mu-\frac{5}{2}}(\La^g_{q_0}\backslash \Xi)$.
\end{enumerate}
Here $a, b, c$ are the coefficients of the nonlinear function $H$ in \eqref{eqhlow}. Moreover, the same conclusions  hold for $\mcu_{ijkl}^{(4)}$ hence for $\mcu^{(4)}$.
\end{prop}
We remark that part (3) was obtained in \cite{KLU}. Also, we emphasize that we stay away from $\La^{(1)}$ where the wave front set of $v_i$ lie, and we are away from the union of $\La_{ijk}^g$ which appears due to the interaction of three waves. In other words, we only look at the new singularities produce by the four wave interactions.  
 
\bpf[Proof of Prop.\ \ref{inter4}]
(1) We consider the term $Q_g(cv_1v_2v_3v_4)$. From Lemma \ref{prod2}, we know that 
\beq
\begin{gathered}
v_1v_2\in I^{\mu, \mu+1}(\La_{12}, \La_1) + I^{\mu, \mu+1}(\La_{12}, \La_2);\\
v_3v_4\in I^{\mu, \mu+1}(\La_{34}, \La_3)  + I^{\mu, \mu+1}(\La_{34}, \La_4).
\end{gathered}
\eeq
Applying Lemma \ref{prod41}, we obtain
\beq
\begin{gathered}
v_1v_2v_3v_4 = w_0 + w_1 + w_2, \ \ w_0 \in I^{4\mu+ 3}(\La_{q_0}), \ \ w_1\in \mcd'(M;  \La^{(1)}(\eps)),\\
w_2 \in \sum_{i, j, k} I^{3\mu + \frac{3}{2}}(\La_{ijk}) + \sum_{i, j} I^{\mu, \mu +1}(\La_{ij}, \La_i) + \sum_{i, j} I^{\mu,\mu +1}(\La_{ij}, \La_j),
\end{gathered}
\eeq
in which the first summation is over $i, j, k = 1, 2, 3, 4$ with $i, j, k$ distinct and the rest two summations are over $i, j = 1, 2, 3, 4$ with $i\neq j$. Finally, applying $Q_g$ and letting $\eps\rightarrow 0$, we obtain that away from the union of $\La_i$ and $\La_{ijk}^g$, 
\beq
Q_g(cv_1v_2v_3v_4) \in I^{4\mu+\frac{3}{2}, -\ha}(\La_{q_0}, \La^g_{q_0}). 
\eeq

(2) Consider $Q_g(bv_1v_2Q_g(av_3v_4))$. First we have 
\beq
\begin{gathered}
Q_g(av_3v_4)\in I^{\mu-\frac{3}{2}, \mu+1-\ha}(\La_{34}, \La_3) + I^{\mu-\frac{3}{2}, \mu+1-\ha}(\La_{34}, \La_4).
\end{gathered}
\eeq
The multiplication with $v_1v_2$ is similar to part (1) and we get 
\beq
Q_g(bv_1v_2Q_g(av_3v_4)) \in I^{4\mu-\ha, -\ha}(\La_{q_0}, \La^g_{q_0}),
\eeq
microlocally away from the union of $\La^{(1)}, \La^{(3), g}$. Next consider $Q_g(av_1Q_g(bv_2v_3v_4))$. By the proof of Prop.\ \ref{inter3}, we know that
\beq
\begin{gathered}
Q_g(bv_2v_3v_4) = w_0 + w_1, \ \ w_0\in I^{3\mu-3, -\ha}(\La_{234}, \La_{234}^g), \\
\WF(w_1)\subset \La^{(1)}(\eps)\cup (\La'_g\circ (\La^{(1)}(\eps)\cap \Sigma_g)). 
\end{gathered}
\eeq
We consider the multiplication with $v_1$. Here we cannot apply Lemma \ref{prod3} because $\La_{234}^g$ is not a conormal bundle. So we use Lemma \ref{prod43} to write $w_0 = w_{0}'+ w_{0}''$ where $w_{0}' \in I^{3\mu-3-\ha}(\La_{234})$ and $\WF(w_{0}'')\subset \La_{234}^g(\eps)$ which denotes a small $\eps$-neighborhood of $\La_{234}^g$. Then by Lemma \ref{prod42}, $w_{0}'v_1\in I^{4\mu-3}(\La_{q_0}) + \mcd'(M; \La_1(\eps)\cup \La_{234}(\eps))$. Finally, by applying $Q_g$, we get 
\beq
\begin{gathered}
Q_g(av_1Q_g(bv_2v_3v_4)) \in I^{4\mu-3-\frac{3}{2}, -\ha}(\La_{q_0}, \La^g_{q_0}) + Q_g(w_{0}''v_1) + Q_g(w_1v_1) \\
+ Q_g(\mcd'(M; \La_1(\eps)\cup \La_{234}(\eps))). 
\end{gathered}
\eeq
The first term on the right hand side is as desired. We analyze the wave front set of the remaining terms. For the last term, we know $\La_1$ is a conic $\eps$-neighborhood of the union of $\La_1$. Under the flow of $Q_g$ i.e.\ $\La_g$, the wave front set is still a $\eps$-neighborhood of $\La_1$. By taking $\eps$ small enough, we see that the wave front set is close to $\La_1$. Similarly, $\La'_g\circ\La_{234}(\eps)$ is a small neighborhood of $\La_{234}^g$ which tends to $\La_{234}^g$ as $\eps\rightarrow 0$.

\begin{figure}[htbp]
\centering
\scalebox{1} 
{
\begin{pspicture} (0,-3.1)(6.8418946,3.1)
\psline[linewidth=0.04cm,arrowsize=0.05291667cm 2.0,arrowlength=1.4,arrowinset=0.4]{<-}(2.8210156,2.6196485)(2.8210156,-0.38035157)
\psline[linewidth=0.04cm,arrowsize=0.05291667cm 2.0,arrowlength=1.4,arrowinset=0.4]{->}(2.8210156,-0.38035157)(5.821016,-0.38035157)
\psline[linewidth=0.04cm,arrowsize=0.05291667cm 2.0,arrowlength=1.4,arrowinset=0.4]{->}(2.8210156,-0.38035157)(1.0210156,-2.1803515)
\psellipse[linewidth=0.04,dimen=outer](1.8810157,-1.2803515)(0.24,0.2)
\psline[linewidth=0.04cm,linestyle=dashed,dash=0.16cm 0.16cm,arrowsize=0.05291667cm 2.0,arrowlength=1.4,arrowinset=0.4]{->}(2.8410156,-0.38035157)(4.301016,-1.9203515)
\psline[linewidth=0.04cm,linestyle=dashed,dash=0.16cm 0.16cm,arrowsize=0.05291667cm 2.0,arrowlength=1.4,arrowinset=0.4]{->}(2.8210156,-0.38035157)(3.6010156,-2.2203515)
\psline[linewidth=0.04cm,linestyle=dashed,dash=0.16cm 0.16cm,arrowsize=0.05291667cm 2.0,arrowlength=1.4,arrowinset=0.4]{->}(2.8010156,-0.36035156)(5.321016,-0.100351565)
\psline[linewidth=0.04cm,linestyle=dashed,dash=0.16cm 0.16cm,arrowsize=0.05291667cm 2.0,arrowlength=1.4,arrowinset=0.4]{->}(2.8210156,-0.38035157)(5.2010155,-0.7803516)
\psline[linewidth=0.04cm,linestyle=dashed,dash=0.16cm 0.16cm,arrowsize=0.05291667cm 2.0,arrowlength=1.4,arrowinset=0.4]{->}(2.8210156,-0.38035157)(1.0210156,-1.6203516)
\psline[linewidth=0.04cm,linestyle=dashed,dash=0.16cm 0.16cm,arrowsize=0.05291667cm 2.0,arrowlength=1.4,arrowinset=0.4]{->}(2.8210156,-0.38035157)(1.5810156,-2.1003516)
\usefont{T1}{ptm}{m}{n}
\rput(6.2724705,-0.41535157){$\La_2$}
\usefont{T1}{ptm}{m}{n}
\rput(4.2524705,-2.7953515){$\La_3$}
\usefont{T1}{ptm}{m}{n}
\rput(0.77247073,-2.54753516){$\La_4$}
\usefont{T1}{ptm}{m}{n}
\rput(3.0824707,0.16464844){$\La_{13}(\eps)$}
\usefont{T1}{ptm}{m}{n}
\rput(4.682471,1.2246485){$\La_{12}(\eps)$}
\usefont{T1}{ptm}{m}{n}
\rput(1.3024707,0.42464843){$\La_{14}(\eps)$}
\psdots[dotsize=0.12](2.8210156,-0.38035157)
\psline[linewidth=0.04cm,arrowsize=0.05291667cm 2.0,arrowlength=1.4,arrowinset=0.4]{->}(2.8210156,-0.38035157)(4.1810155,-2.4603515)
\psellipse[linewidth=0.04,dimen=outer](3.6210155,-1.5603516)(0.26,0.16)
\psbezier[linewidth=0.04,linestyle=dashed,dash=0.16cm 0.16cm](2.8210156,1.6196485)(2.1810157,1.4196484)(1.6810156,-0.16035156)(1.7010156,-1.1403515)
\psbezier[linewidth=0.04,linestyle=dashed,dash=0.16cm 0.16cm](2.8210156,1.5796485)(2.4210157,1.0196484)(1.9410156,-0.56035155)(2.0610156,-1.4203515)
\psbezier[linewidth=0.04,linestyle=dashed,dash=0.16cm 0.16cm](2.8210156,1.6196485)(3.7810156,1.5596484)(4.5610156,0.49964845)(4.6410155,-0.18035156)
\psbezier[linewidth=0.04,linestyle=dashed,dash=0.16cm 0.16cm](2.8610156,1.6196485)(3.9610157,1.0250539)(4.4010158,-0.40035155)(4.341016,-0.62035155)
\psbezier[linewidth=0.04,linestyle=dashed,dash=0.16cm 0.16cm](2.8210156,1.6796484)(3.2410157,0.81964844)(3.2410157,-0.5203516)(3.3410156,-1.5403515)
\psbezier[linewidth=0.04,linestyle=dashed,dash=0.16cm 0.16cm](2.8210156,1.6396484)(3.8010156,1.0038043)(3.7610157,-1.1026893)(3.8210156,-1.4603516)
\psbezier[linewidth=0.04](4.8110156,-0.26432106)(4.9010158,-0.4482906)(4.361016,-0.82035154)(4.361016,-0.548688)(4.361016,-0.27702442)(4.7210155,-0.08035156)(4.8110156,-0.26432106)
\usefont{T1}{ptm}{m}{n}
\rput(3.224706,2.6246484){$\La_1$}
\end{pspicture} 
}
\caption{Explanation of the wavefront set on $\La_{q_0}$. Each solid circle represents the $\eps$-neighborhood of $\La_i\cap S^*M$. The regions bounded by dashed curves are the linear span of $\La_1$ and $\La_i(\eps)$. }
\label{figcut4}
\end{figure}
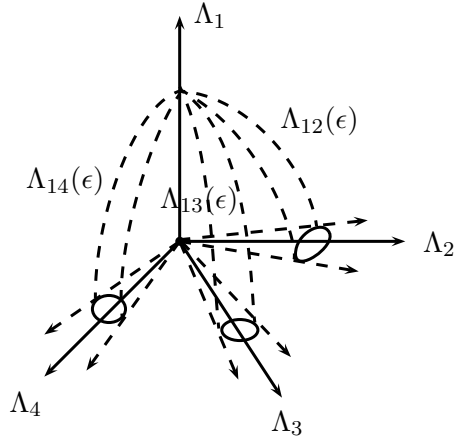

Next, consider the wave front set of $Q_g(w_1v_1)$. Away from $q_0$, $\WF(w_1v_1)$ is contained in a neighborhood $\La^{(1)}(\eps)$, the flow out of which is still close to $\La^{(1)}$. At $q_0$, $\WF(w_1v_1)$ is contained in the linear span of $\La^{(1)}(\eps)$ and $\La_1(\eps)$ over $q_0$, see Figure \ref{figcut4}. The result is a $\eps$-neighborhood of the union of $\La_{1i}, i = 2, 3, 4$. For $\eps$ small enough, this is also close to $\La^{(1)}\cup \La^{(2)}$. So the flow out under $\La^g$ is still close to $\La^{(1)}$.
 
 Finally, consider the wave front set of $Q_g(w_0''v_1)$. Away from $q_0$, $\WF(w_0''v_1)$  is contained in the span of $\La_{1}$ and $\La^g_{234}(\eps)$ which is empty for $\eps$ sufficiently small. At $q_0$, $\WF(w_{0}''v_1)$ is contained in the span of $\La_1$ and $\La_{234}(\eps)\cap \Sigma_g$ over $q_0$. By the similar argument in Lemma \ref{inter2}, the vector $\zeta\in \WF(w_{0}''v_1)$ is light-like if and only if $\zeta \in \La_1$ or $\zeta\in \La_{234}(\eps)\cap \Sigma_g$. Thus under the flow out of $\La_g,$ the wave front set is contained in $\La_1\cup \La_{234}^g(\eps)$. 

(3) We analyze the rest two terms in a similar fashion. First of all, we have
\beq
\begin{gathered}
Q_g(av_1v_2)\in I^{\mu-\frac{3}{2}, \mu+1-\ha}(\La_{12}, \La_1) + I^{\mu-\frac{3}{2}, \mu+1-\ha}(\La_{12}, \La_2);\\
Q_g(av_3v_4)\in I^{\mu-\frac{3}{2}, \mu+1-\ha}(\La_{34}, \La_3) + I^{\mu-\frac{3}{2}, \mu+1-\ha}(\La_{34}, \La_4).
\end{gathered}
\eeq
Thus the analysis of $Q_g(aQ_g(av_1v_2)Q_g(av_3v_4))$ is the same as in case (1). Next, from the proof of Prop.\ \ref{inter3}, we know that
$Q_g(av_2Q_g(av_3v_4))$ has a similar structure as $Q_g(bv_2v_3v_4)$ in case (2). Thus  $Q_g(av_1Q_g(av_2Q_g(av_3v_4)))$ can be analyzed as in case (2).
\epf

\subsection{Leading singularities and principal symbols}\label{nonsym}
In this subsection, we compute the principal symbols of terms in $\mcu^{(4)}$. Our purpose is twofold. First we show that the singularities are non-vanishing by a proper choice of $v_i$. More importantly, the symbols contains information of the metric which is further explored later.  

\begin{prop}\label{nonvan}
Let $K_i,$ $i = 1, 2, 3, 4$ be submanifolds of $M$  having codimension $1$ such that $\Lambda_i=N^*K_i$  are light-like.
Assume that the set $\bigcap_{i=1}^4 K_i$ contains only one point  $q_0$ 
and that the manifolds $K_i$  intersect transversally at $q_0$, see
Definition \ref{trans4}. Let $\xi_i\in N^*_{q_0}K_i$ be some fixed normal vectors and $v_i\in I^\mu( \Lambda_i)$. 
Let $(q, \eta) \in \La_{q_0}^g\backslash \Xi$ which is joined to $(q_0, \zeta)\in \La_{q_0}$ by bicharacteristics. If $\zeta  =\sum_{i = 1}^4\zeta_i$ with $\zeta_i = r_i \xi_i$ and $A_i(q_0, \zeta_i)$ be the principal symbols of $v_i$ at $(q_0, \zeta_i)$,  we have the following conclusions.
\begin{enumerate}
\item 
The principal symbol  of $\mcu^{(4)}$ can be written as
\beqq\label{symu4}
\sigma(\mcu^{(4)})(q, \eta) = P(\zeta_1, \zeta_2, \zeta_3, \zeta_4) Q_{\La_{q_0}^g}(q, \eta, q_0, \zeta) \prod_{i = 1}^4 A_i(q_0, \zeta_i),
\eeqq
where $P$ is specified below. 
\begin{enumerate}
\item If $c(q_0)\neq 0$, we have 
\beq
P(\zeta_1, \zeta_2, \zeta_3, \zeta_4) = -4! (2\pi)^{-3} c(q_0).
\eeq
\item If $c = 0$ in a neighborhood of $q_0$, where $b, a$ are non-vanishing,  we have 
\beq
\begin{gathered}
P(\zeta_1, \zeta_2, \zeta_3, \zeta_4) =   (2\pi)^{-3} a(q_0) b(q_0)  \sum_{i, j, k, l}[ \frac{3}{|\zeta_i+\zeta_j|_{g^*(q_0)}^{2}} + \frac{2}{|\zeta_j+\zeta_k+\zeta_l|_{g^*(q_0)}^{2}}],
\end{gathered}
\eeq
where the summation in $(i, j, k, l)$ runs over all permutations of $(1, 2, 3, 4)$. 
\item If $b = c = 0$ in a neighborhood of $q_0$ where $a$ is non-vanishing,  we have  
\beq
\begin{split}
P(\zeta_1, \zeta_2, \zeta_3, \zeta_4) &=  -  (2\pi)^{-3} a(q_0)^3  \cdot \\
& \sum_{i, j, k, l}[ \frac{4}{|\zeta_j+\zeta_k+\zeta_l|_{g^*(q_0)}^2 \cdot |\zeta_k+\zeta_l|_{g^*(q_0)}^{2}} + \frac{1}{|\zeta_i+\zeta_j|_{g^*(q_0)}^{2} \cdot |\zeta_k+\zeta_l|_{g^*(q_0)}^{2}}] ,
\end{split}
\eeq
where the summation in $(i, j, k, l)$ runs over all permutations of $(1, 2, 3, 4)$.
\end{enumerate}

\item The coefficients $P(\zeta_1, \zeta_2, \zeta_3, \zeta_4)$ considered in the above cases (1)-(3) can be regarded as real analytic functions of $(\xi_1, \xi_2, \xi_3, \xi_4, \zeta)$ defined on the set
\beq
\mcx = \{(\xi_1, \xi_2, \xi_3, \xi_4, \zeta) \in (L_{q_0}^*M)^5: \xi_1, \xi_2, \xi_3, \xi_4 \text{ are linearly independent. }\}
\eeq
Then $P$  do not vanish in any open subset of  $\mcx$. This in particular implies
that there is an open, conic, and dense set $\mathcal W\subset   ({L^{*}_{q_0}M})^4$ such  that for all
$(\xi_1,\xi_2,\xi_3,\xi_4)\in \mathcal W$ 
it holds that if the principal symbols $A_i(q_0, \zeta_i)$  of $v_i$,  at $\zeta_i=r_i\xi_i$, are all non-zero, then $\sigma(\mcu^{(4)})$ is non-vanishing on any open subset of $\La_{q_0}^g\backslash \Xi$.  
\end{enumerate}

\end{prop}

Roughly speaking, the proposition says that if $K_i, i = 1, 2, 3, 4$ intersect at $q_0$, we have in a generic  case that 
  the principal symbol of $\mcu^{(4)}$ does not vanish identically on any open subset of the future light cone at $q_0$. 
%
In particular,
we observe that for any $(q, \eta)\in \La_{q_0}^g\backslash \Xi$ which is joined to $(q_0, \zeta) \in \La_{q_0}$ by bicharacteristics, one can always choose $K_i$ intersecting at $q_0$ such that $\sigma(\mcu^{(4)})$ is non-vanishing at $(q, \eta)$. 
Below,  we say that the submanifolds $K_i$ intersect
in  a generic way when the co-normal vectors $\xi_j\in N^*K_j$ are such that $(\xi_1,\xi_2,\xi_3,\xi_4)\in \mathcal W$. 

\bpf[Proof of Prop.\ \ref{nonvan}]
\textbf{Part (1):} We know that
\beq
\sigma(\mcu^{(4)})(q, \eta) = \sum_{i, j, k, l} \mcu^{(4)}_{ijkl}(q, \eta),
\eeq
where the summation runs over all permutations of $(1, 2, 3, 4)$. The terms $\mcu^{(4)}_{ijkl}$ are similar to $\mcu^{(4)}_{1234}$ and we start with the principal symbols of $\mcu^{(4)}_{1234}$. For any $\zeta\in \La_{q_0}$, we can write $\zeta = \sum_{i = 1}^4\zeta_i$ where $\zeta_i\in N^*_{q_0}K_i$. Consider $\mcy_1 = Q_g(cv_1v_2v_3v_4)$ when $c\neq 0$. First of all, we have
\beq
\sigma_{\La_{q_0}}(cv_1v_2v_3v_4)(q_0, \zeta) = (2\pi)^{-3} c(q_0)\prod_{i = 1}^4A_i(q_0, \zeta_i).
\eeq
Again, we ignored the $2\pi$ factors in symbol computations. By \cite[Prop.\ 2.1]{GrU93}, we get 
 \beq
 \sigma_{\La_{q_0}^g}(\mcy_1)(q, \eta) = (2\pi)^{-3} \sigma_{\La_g}(Q_g)(q, \eta, q_0, \zeta)c(q_0)\prod_{i = 1}^4A_i(q_0, \zeta_i),
\eeq
 where $(q, \eta)$ is joined with $(q_0, \zeta)$ by bicharacteristics.   Next consider $\mcy_2 = Q_g(bv_1v_2Q_g(av_3v_4))$ and $\mcy_3 = Q_g(av_1Q_g(bv_2v_3v_4))$. These are the terms with $b$. By similar arguments, we find that
 \beq
 \begin{gathered}
  \sigma_{\La_{q_0}^g}(\mcy_2)(q, \eta)= (2\pi)^{-3}\sigma_{\La_g}(Q_g)(q, \eta, q_0, \zeta)b(q_0) A_1(q_0, \zeta_1) A_2(q_0, \zeta_2) |\zeta_3+\zeta_4|_{g^*(q_0)}^{-2}\\
 \cdot a(q_0)A_3(q_0, \zeta_3)A_{4}(q_0, \zeta_4),\\
    \sigma_{\La_{q_0}^g}(\mcy_3)(q, \eta) = (2\pi)^{-3}\sigma_{\La_g}(Q_g)(q, \eta, q_0, \zeta)a(q_0) A_1(q_0, \zeta_1)|\zeta_2+\zeta_3+\zeta_4|_{g^*(q_0)}^{-2} \\
    \cdot b(q_0)A_2(q_0, \zeta_2) A_3(q_0, \zeta_3)A_4(q_0, \zeta_4).
  \end{gathered}
 \eeq 
 Finally, we consider the terms in $\mcu^{(4)}$ which only have $a$. This is the case when $b = c = 0$ studied in \cite{KLU}. Let $\mcy_4 =Q_g(av_1Q_g(av_2Q_g(av_3v_4)))$, $\mcy_5 = Q_g(aQ_g(av_1v_2)Q_g(av_3v_4))$. Then we find
 \beq
 \begin{gathered}
  \sigma_{\La_{q_0}^g}(\mcy_4)(q, \eta)= (2\pi)^{-3}\sigma_{\La_g}(Q_g)(q, \eta, q_0, \zeta)a^3(q_0) A_1(q_0, \zeta_1)  |\zeta_2+\zeta_3+\zeta_4|_{g^*(q_0)}^{-2} A_2(q_0, \zeta_2)\\
\cdot|\zeta_3+\zeta_4|_{g^*(q_0)}^{-2} A_3(q_0, \zeta_3)A_4(q_0, \zeta_4).\\
    \sigma_{\La_{q_0}^g}(\mcy_5)(q, \eta)= (2\pi)^{-3}\sigma_{\La_g}(Q_g)(q, \eta, q_0, \zeta)a^3(q_0) |\zeta_3+\zeta_4|_{g^*(q_0)}^{-2} A_4(q_0, \zeta_4) A_3(q_0, \zeta_3) \\
 \cdot|\zeta_1+\zeta_2|_{g^*(q_0)}^{-2} A_1(q_0, \zeta_1)A_2(q_0, \zeta_2).
  \end{gathered}
 \eeq 

Now we find the principal symbols of $\mcu^{(4)}$ and show that it is non-vanishing on any open subset of the forward light-cone of $q_0.$ It suffices to show that for $\zeta\in L_{q_0}^{*, +}M$ away from $\La^{(1)}\cup \La^{(3)}$, the principal symbol $\mcu^{(4)}(q_0, \zeta)$ is non-vanishing. This is because
\beq
\sigma(\mcu^{(4)})(q, \eta) = Q_{\La_{q_0}^g}(q, \eta, q_0, \zeta)\sigma(\mcu^{(4)})(q_0, \zeta).
\eeq
First, when $c$ is non-vanishing, we have
\beq
\sigma(\mcu^{(4)})(q_0, \zeta) = -4! (2\pi)^{-3} c(q_0) \prod_{i = 1}^4 A_i(q_0, \zeta_i),
\eeq
where $\zeta = \sum_{i = 1}^4 \zeta_i$ with $\zeta_i, \zeta$  light-like vectors. Therefore the symbol is obviously non-vanishing if $A_i(q_0, \zeta_i)$ are non zero. This proves part (a).

Next consider part (b). We know that the principal part of $\mcu^{(4)}$ is a sum of terms like $\mcy_2, \mcy_3$ and we find explicitly that
\beqq\label{sym2}
\begin{gathered}
\sigma(\mcu^{(4)})(q_0, \zeta) =  \sum_{i, j, k, l}[ \frac{3}{|\zeta_i+\zeta_j|_{g^*(q_0)}^{2}} + \frac{2}{|\zeta_j+\zeta_k+\zeta_l|_{g^*(q_0)}^{2}}]   \cdot (2\pi)^{-3} a(q_0) b(q_0) \prod_{i = 1}^4 A_i(q_0, \zeta_i).
\end{gathered}
\eeqq

Part (c) is similar. We know that the principal part of $\mcu^{(4)}$ is a sum of terms like $\mcy_4, \mcy_5$. We can write the principal symbol as
\beqq\label{sym3}
\begin{gathered}
\sigma(\mcu^{(4)})(q_0, \zeta) =  -\sum_{i, j, k, l}[ \frac{4}{|\zeta_j+\zeta_k+\zeta_l|_{g^*(q_0)}^2 \cdot |\zeta_k+\zeta_l|_{g^*(q_0)}^{2}} + \frac{1}{|\zeta_i+\zeta_j|_{g^*(q_0)}^{2} \cdot |\zeta_k+\zeta_l|_{g^*(q_0)}^{2}}] \\
\cdot (2\pi)^{-3} a(q_0)^3  \prod_{i = 1}^4 A_i(q_0, \zeta_i),
\end{gathered}
\eeqq
where the summation is over permutations of $(1, 2, 3, 4)$.

\textbf{Part (2):} We start with the real analytic structure of $\mcx$. Let $\xi_1,\xi_2,\xi_3,\xi_4\in L^{*}_{q_0}M$ be a basis of the vector space $T^*_{q_0}M$ and consider 
vectors $\zeta_j=r_j \xi_j$, where $r_j\in \R$. We can write the future directed light cone
$ L^{*,+}_{q_0}M$  in the coordinates corresponding to the basis vectors  $\xi_1,\xi_2,\xi_3,\xi_4$ as
\ba
S=\{\vec r=(r_1,r_2,r_3,r_4)\in  \R^4\setminus \{0\}&;&\ \sum_{j,k=1}^4 r_jr_k g(\xi_j,\xi_k)=0,\  
\sum_{j, =1}^4 r_j  g(\xi_j,\eta)<0\},
\ea
where $\eta$ is a future directed time-like vector.
The set $ L^{*,+}_{q_0}M$ is connected, so its representation $S$  in the  coordinates corresponding to the basis vectors  $\xi_1,\xi_2,\xi_3,\xi_4$ is also connected. For all $\vec r\in S$ there is $j=j(\vec r)\in \{1,2,3,4\}$  such that
$$
\sum_{k\not =j} r_kg(\xi_j,\xi_k)=  g(\xi_j,\sum_{k=1}^4 r_k\xi_k)= g(\xi_j,\vec r)\not =0,
$$
as otherwise $\vec r$ would be zero. Thus near any $\vec r\in S$ we can use the three variables $r_k,$ $k\in \{1,2,3,4\}\setminus \{ j(\vec r)\}$,
as local coordinates and in these coordinates $S$ is given in a neighborhood of a point $\vec r$ by
\ba
r_j=-\frac {\sum_{k,l\in \{1,2,3,4\}\setminus \{ j\}} r_kr_l g(\xi_k,\xi_l)}{\sum_{i\in \{1,2,3,4\}\setminus \{ j\}} r_i g(\xi_i,\xi_j)}  ,\quad \hbox{where }j={j(\vec r)}.
\ea
Thus $S$ is a connected real analytic manifold. Also, functions $ P$ considered above can be written as 
$$
P_j(\zeta_1,\zeta_2,\zeta_3,\zeta_4)=\frac {Q(r_1,r_2,r_3,r_4)} {R(r_1,r_2,r_3,r_4)},
$$
where $ {Q(r_1,r_2,r_3,r_4)}$ and $ {P(r_1,r_2,r_3,r_4)}$ are polynomials and thus real-analytic functions on $S$.
As real analytic function on a real analytic manifold vanishes in an open set only if it vanishes in a topological component of $S$, it suffices to show that functions $P(\zeta_1,\zeta_2,\zeta_3,\zeta_4)$ are non-zero at some points of $S$. Then they are non-vanishing in an open and dense set.


To do the computation, without loss of generality, we can assume that the metric at $q_0$ is Minkowski, so $g(q_0) = \text{diag}(-1, 1, 1, 1)$. In this case, the dual metric $g^*(q_0) = g(q_0)$, so we can identify vectors and covectors. We choose $\xi_i$ as following
\ba
\xi_1=(1,0,1,0),\ \ \xi_2=(1,0,0,1), \\ 
\xi_3=(-1, -1, 0,0),\ \ \xi_4=(1,-1,0,0).
\ea
These are light-like vectors and linearly independent.  To make the leading term simpler, now we choose $\alpha_1 = 1, \alpha_3 = \rho, \alpha_4 = \rho^{10}$ with $\rho$ a small parameter, and solve for
\beq
\alpha_2 = \frac{\alpha_3 + \alpha_4}{1-\alpha_3-\alpha_4} = \rho + O(\rho^{10}).
\eeq
such that $\zeta = \sum_{i = 1}^4 \alpha_i\xi_i$ is light-like. We take $\zeta_i = \alpha_i\xi_i$ i.e.\
\beqq\label{zetas}
\begin{gathered}
\zeta_1=(1,0,1,0),\ \ \zeta_2=(\rho + O(\rho^{10}))(1,0,0,1), \\ 
\zeta_3= \rho (-1, -1, 0,0),\ \ \zeta_4= \rho^{10}(1,-1,0,0).
\end{gathered}
\eeqq
Notice that $|\zeta_i + \zeta_j|_{g}^2 = 2g(\zeta_i, \zeta_j)$, and we compute that
\beq
\left.\begin{array}{ll}
g(\zeta_1,\zeta_2)=-\rho + O(\rho^{10}), \ \ & g(\zeta_2,\zeta_3)=\rho^2 + O(\rho^{11}), \\ 
g(\zeta_1,\zeta_3)=\rho, & g(\zeta_2,\zeta_4)=-\rho^{11} + O(\rho^{20}), \\
g(\zeta_1,\zeta_4)=-\rho^{10}, & g(\zeta_3,\zeta_4)=2\rho^{11}.
\end{array}\right.
\eeq
Also, we have that 
\ba
&&|\zeta_1 + \zeta_2+ \zeta_3|_g^2 = 2g(\zeta_1, \zeta_2) + 2g(\zeta_1, \zeta_3) + 2g(\zeta_2, \zeta_3) = 2\rho^2 + O(\rho^{10}),\\
&&|\zeta_1 + \zeta_3+ \zeta_4|_g^2 = 2\rho+ O(\rho^{10}),\\
&&|\zeta_2 + \zeta_3+ \zeta_4|_g^2 = 2\rho^2 + O(\rho^{10}).
\ea

In the symbol \eqref{sym2}, we consider 
\beq
\mcp_2 = \sum_{i, j, k, l}[ \frac{3}{|\zeta_i+\zeta_j|_{g^*(q_0)}^{2}} + \frac{2}{|\zeta_j+\zeta_k+\zeta_l|_{g^*(q_0)}^{2}}].
\eeq
The leading terms can be determined as $\rho\rightarrow 0$ because the smallest terms are $g(\zeta_2, \zeta_4)$ and $g(\zeta_3, \zeta_4)$. So we find that 
\beq
\begin{split}
\mcp_2   & = 2(\frac{3}{|\zeta_2+\zeta_4|_{g^*(q_0)}^{2}} + \frac{3}{|\zeta_3+\zeta_4|_{g^*(q_0)}^{2}})  + O(\rho^{-10}) \\
 & = -\frac{3}{2\rho^{11}} + O(\rho^{-10}).
\end{split}
\eeq
For $\rho$ small, this is non-vanishing. 

In the symbol \eqref{sym3}, we consider the leading term in 
\beq
\begin{gathered}
\mcp_3 =  -\sum_{i, j, k, l}[ \frac{4}{|\zeta_j+\zeta_k+\zeta_l|_{g^*(q_0)}^2 \cdot |\zeta_k+\zeta_l|_{g^*(q_0)}^{2}} + \frac{1}{|\zeta_i+\zeta_j|_{g^*(q_0)}^{2} \cdot |\zeta_k+\zeta_l|_{g^*(q_0)}^{2}}].
\end{gathered}
\eeq
We observe that 
\beq
\begin{gathered}
\sum_{i, j, k, l}\frac{1}{|\zeta_i+\zeta_j|_{g^*(q_0)}^{2} \cdot |\zeta_k+\zeta_l|_{g^*(q_0)}^{2}} = O(\rho^{-12}),
\end{gathered}
\eeq
but
\beq
\frac{4}{|\zeta_2+\zeta_3+\zeta_4|_{g^*(q_0)}^{2} \cdot |\zeta_3+\zeta_4|_{g^*(q_0)}^{2}} = \rho^{-13}(\frac{1}{1 + O(\rho^8)}) = \rho^{-13} + O(\rho^{-12})
\eeq
and the other terms are all of the order $O(\rho^{-12})$. Therefore, 
\beq
\begin{gathered}
\mcp_3 =  -2\rho^{-13} + O(\rho^{-12}). 
\end{gathered}
\eeq
As $\rho\rightarrow 0$, we see this term is non-zero. This completes the proof of the proposition.

\epf

\section{Solution to inverse problems for lower order nonlinearities}\label{secinv}
In this section, we use the singularities analyzed in Section \ref{singu} and the method in Kurylev-Lassas-Uhlmann \cite{KLU} to prove our main theorems for a special case when $H$ has lower order nonlinearities.
\begin{theorem}\label{main11}
Let $(M^{(j)}, g^{(j)}), j = 1, 2$ be two $4$-dimensional globally hyperbolic Lorentzian manifolds. Let $\hat \mu^{(j)}(t) \subset M^{(j)}$ be time-like geodesics where $t\in [-1, 1]$ and $p_\pm^{(j)}=\hat\mu^{(j)}(s_\pm), -1<s_-<s_+<1$. Let $V^{(j)}\subset M$ be open relatively compact neighborhoods of $\hat \mu^{(j)}([s_-, s_+])$  and $M_0^{(j)} = (-\infty, T_0)\times N^{(j)}, T_0 > 0$ be such that $V^{(j)}\subset M_0^{(j)}$. Consider the semilinear wave equations with source terms
\beqq\label{eqsem2}
\begin{gathered}
\square_{g^{(j)}} u(x) + H^{(j)}(x, u(x))  = f(x), \text{ on } M^{(j)}_0,\\
u = 0 \text{ in } M_0^{(j)}\backslash J^+_{g^{(j)}}(\supp(f)),
\end{gathered}
\eeqq
where $\supp(f)\subset V^{(j)}$ and
\beq
H^{(j)}(x, z) = a^{(j)}(x)z^2 + b^{(j)}(x)z^3 + c^{(j)}(x)z^4,\ \ i = 1, 2,
\eeq
where $[a^{(j)}]^2 + [c^{(j)}]^2$ are non-vanishing. Assume that there is a diffeomorphism $\Phi: V^{(1)}\rightarrow V^{(2)}$ such that $\Phi(p_\pm^{(1)}) = p_\pm^{(2)}$ and the source-to-solution maps satisfy
\beq
(\Phi^{-1})^*(L^{(1)}(\Phi^*f)) = L^{(2)}(f)
\eeq
for all $f$  in a small neighborhood of the zero function in $C_0^4(V^{(2)})$. Then we have the following conclusions. 
\begin{enumerate}
\item There exists a diffeomorphism $\Psi: I(p^{(1)}_-, p^{(1)}_+)\rightarrow I(p^{(2)}_-, p^{(2)}_+)$ such that 
$\Psi^*g^{(2)} = e^{2\gamma} g^{(1)}$ in $I(p^{(1)}_-, p^{(1)}_+)$ for some $\gamma\in C^\infty(I(p^{(1)}_-, p^{(1)}_+))$. Moreover, $\Psi = \Phi$ on $I(p^{(1)}_-, p^{(1)}_+) \cap V^{(1)}$.
\item In addition, if $H^{(j)}(x, z), i = 1, 2$ are independent of $x$, i.e.\ $a^{(j)}, b^{(j)}, c^{(j)}$ are constants, then the conformal diffeomorphism is an isometry, meaning  $\Psi^*g^{(2)} = g^{(1)}$ in $I(p_-^{(1)}, p_+^{(1)})$. 
\end{enumerate}
\end{theorem} 

We also prove the analogue of Theorem \ref{main2} and \ref{maingh} in this section, and we shall return to the general case in Section 6. We remark that the strategy for solving the inverse problems is the same for the general case. However, to deal with higher order nonlinear terms, we will need higher order asymptotic expansions of the solution. The analysis of those singularities involves other techniques. Also, notice that in Theorem \ref{main11} we leave out the case when $H(x, z) = b(x) z^3$ i.e.\ $H$ is cubic. We treat this case after we determine the higher order nonlinear terms. 

To begin with, we prove the analogue of Theorem 1.5 in \cite{KLU} for our semilinear equation, which says that the source-to-solution map determines the conformal class of the metric. Except for the analysis of singularities, the proof of Theorem \ref{main11} heavily relies on the work of Kurylev-Lassas-Uhlmann \cite{KLU}. We will not repeat their proofs here. Instead, we will point out which arguments are used below and refer the interested reader to \cite{KLU} for more details.

\subsection{Distorted plane waves}\label{disto}
We construct conormal distributions propagating along geodesics as in \cite{KLU}.  
Let $(x_0, \theta_0)\in L^+M_0$ and $\gamma_{x_0, \theta_0}(t), t\geq 0$ be the geodesic from $x_0$ with direction $\theta_0$. For a small parameter $s_0>0$, we let 
\beq
K(x_0, \theta_0; t_0, s_0) = \{\gamma_{x', \theta'}(t)\in M_0; \theta' \in \mco(s_0), t\in (0, \infty)\},
\eeq
where $(x', \theta') = (\gamma_{x_0, \theta_0}(t_0), \gamma_{x_0, \theta_0}'(t_0))$ and $\mco(s_0) \subset L^+_{x'}M$ is a open neighborhood of $\theta'$ consisting of $\zeta\in L_{x'}^+M$ such that $\|\zeta - \theta'\|_{g^+}< s_0$. Notice that as $s_0\rightarrow 0$, $K(x_0, \theta_0; t_0, s_0)$ tends to the geodesic $\gamma_{x_0, \theta_0}$.  Next, let 
\beq
Y(x_0, \theta_0; t_0, s_0) = K(x_0, \theta_0; t_0, s_0)  \cap \{t = 2t_0\},
\eeq
be a $2$-dimensional surface intersecting the geodesic at $\gamma_{x_0, \theta_0}(2t_0)$. We let $\La(x_0, \theta_0; t_0, s_0)$ be the Lagrangian submanifold obtained by flowing out $N^*K(x_0, \theta_0; t_0, s_0)\cap N^*Y(x_0, \theta_0; t_0, s_0)$ under the Hamilton vector field $H_\mcp$ in $\Sigma_g$ (see Section \ref{seclinwa}).

Assume that $\mu\leq -11$ and let $f_0\in I^{\mu +1}(Y) \subset H^{8}(M_0)$ be supported in a neighborhood $U$ of $\gamma_{x_0, \theta_0}\cap Y$. By Sobolev embedding, we know that $f_0 \in C_0^4(M_0)$. Assume the principal symbol of $f_0$ vanishes outside of $\mco(s_0)\subset T^*M$. According to Lemma 3.1 of \cite{KLU}, $v_0 = Q_g f_0$ belongs to $I^{\mu -\ha}(M_0\backslash Y; \La)\subset H^{9}(M_0)$. Moreover, the principal symbol of $v_0$ satisfies
\beqq\label{wavesym}
\sigma(v_0)(x, \xi) = \sigma(Q_g)(x, \xi, y, \eta)\sigma(f_0)(y, \eta),
\eeqq
where $(x, \xi)$ and $(y, \eta)$ lie on the same bicharacteristics. Actually this is a consequence of Prop.\ \ref{winv}. We emphasis that since we will take $s_0$ small, the conormal distribution $v_0$  indeed should be regarded as associated with the geodesic $\gamma_{x_0, \theta_0}$.

Now we consider the conjugate points along $\gamma_{x_0, \theta_0}$. Let $\bold{t_0}=\bold{t_0}(x_0, \theta_0)>0$ be such that $\gamma_{x_0, \theta_0}(\bold{t_0})$ is the first conjugate point of $x_0$ along $\gamma_{x_0, \theta_0}$. Then the exponential map $\exp_{x_0}$ is a local diffeomorphism from a neighborhood of $t\theta_0\in T_{x_0}M$ to a neighborhood of $\gamma_{x_0, \theta_0}(t)$ for $t< \bold{t_0}$. Therefore, $K(x_0, \theta_0; t_0, s_0)$ is a codimension $1$ submanifold near $\gamma_{x_0, \theta_0}(t)$ and 
\beq
N^*K(x_0, \theta_0; t_0, s_0) = \La \text{ near } \gamma_{x_0, \theta_0}(t) \text{ for }  t< \bold{t_0}.
\eeq
Therefore, before the first conjugate point of $x_0$ along $\gamma_{x_0, \theta_0}$, $v_0$ is a distribution conormal to $K(x_0, \theta_0; t_0, s_0)$. It is worth mentioning that in general $v_0$ is a Lagrangian distribution and only before the first conjugate point, it is a conormal distribution. In \cite{KLU}, a stronger notion of (null) cut points was used, see \cite[Section 2.1]{KLU}. Recall that on a globally hyperbolic Lorentzian manifold, the first null cut point $\hat x$ of $x_0$ along the geodesic $\gamma_{x_0, \theta_0}$ is either the first conjugate point or there are at least two light-like geodesics joining $x_0$ and $\hat x$. In particular, the first cut point appears on or before the first conjugate point.

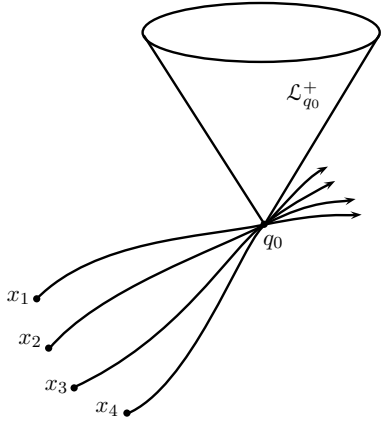
\begin{figure}[htbp]
\scalebox{.8} 
{
\begin{pspicture} (0,-3.5188477)(6.0210156,3.4788477) 
\psellipse[linewidth=0.04,dimen=outer](4.0310154,2.9688478)(1.99,0.51)
\psline[linewidth=0.04cm](2.0810156,2.9188476)(4.0610156,-0.18115234)
\psline[linewidth=0.04cm](6.0010157,2.9388475)(4.1010156,-0.20115234)
\psdots[dotsize=0.12](4.0810156,-0.22115235)
\psbezier[linewidth=0.04,arrowsize=0.05291667cm 2.0,arrowlength=1.4,arrowinset=0.4]{<-}(5.6010156,0.19884765)(5.321016,0.13884765)(4.9010158,0.19884765)(4.1010156,-0.22115235)(3.3010156,-0.6411523)(1.3010156,-1.3211523)(0.5010156,-2.3011522)
\psbezier[linewidth=0.04,arrowsize=0.05291667cm 2.0,arrowlength=1.4,arrowinset=0.4]{<-}(5.2610154,0.49884766)(5.0410156,0.37884766)(4.619152,0.18169029)(4.0579753,-0.23587316)(3.4967985,-0.6534366)(2.6610155,-2.1411524)(0.90101564,-2.9411523)
\psbezier[linewidth=0.04,arrowsize=0.05291667cm 2.0,arrowlength=1.4,arrowinset=0.4]{<-}(5.7010155,-0.061152343)(5.341016,-0.08115234)(5.095066,-0.021152344)(4.066549,-0.23026125)(3.0380318,-0.43937016)(1.2678777,-0.47876698)(0.30101562,-1.4611523)
\psbezier[linewidth=0.04,arrowsize=0.05291667cm 2.0,arrowlength=1.4,arrowinset=0.4]{<-}(5.1410155,0.7388477)(4.7610154,0.5388477)(4.4010158,0.078847654)(4.0410156,-0.28115234)(3.6810157,-0.6411523)(2.8410156,-2.9611523)(1.7810156,-3.3811524)
\psdots[dotsize=0.12](0.30101562,-1.4611523)
\psdots[dotsize=0.12](0.5010156,-2.2811522)
\psdots[dotsize=0.12](0.9210156,-2.9411523)
\psdots[dotsize=0.12](1.8010156,-3.3611524)
\usefont{T1}{ptm}{m}{n}
\rput(4.2324707,-0.51615233){$q_0$}
\usefont{T1}{ptm}{m}{n}
\rput(0,-1.4361523){$x_1$}
\usefont{T1}{ptm}{m}{n}
\rput(0.1924707,-2.1561522){$x_2$}
\usefont{T1}{ptm}{m}{n}
\rput(0.6324707,-2.9161522){$x_3$}
\usefont{T1}{ptm}{m}{n}
\rput(1.4724707,-3.2961524){$x_4$}
\usefont{T1}{ptm}{m}{n}
\rput(4.7424707,1.9238477){$\mcl^{+}_{q_0}$}
\end{pspicture} 
}
\caption{The interaction of four distorted plane waves at $q_0$. The waves propagate along the geodesics $\gamma_{x_i, \theta_i}$. The interaction produces new singularities at $q_0$ which are propagated to the forward light-cone $\mcl^+_{q_0}$ by $Q_g$.}
\label{fourinter}
\end{figure}
 
Now let $V\subset M_0$ be a neighborhood of a time-like geodesic $\hat \mu([-1, 1])$. Assume $x_j \in V$ and $(x_j, \theta_j)\in L^+M, j = 1, 2, 3, 4$ such that 
\beq
\gamma_{x_j, \theta_j}([0, t_0])\subset V, \ \ x_j(t_0)\notin J^+(x_k(t_0)),
\eeq
which means that the points are causally independent. We define $K_j = K(x_j, \theta_j; t_0, s_0)$ and $\La_j, j = 1, 2, 3, 4$ similar to $K$ and $\La$.  Also, we let $f_j\in I^{\mu+1}(Y_j), \mu\leq -11$ be constructed as $f_0$ above and $v_j = Q_g(f_j)\in I^{\mu-\ha}(\La_j)$. Let $\bold{t_j}, j = 1, 2, 3, 4$ be such that $\gamma_{x_j, \theta_j}(\bold{t_j})$ is the first conjugate point of $x_j$ along the geodesics and $\bold{t} = \min_{j = 1, 2, 3, 4}(\bold{t_j})$. We see that outside the future of the point $\gamma_{x_j, \theta_j}(\bold{t})$, we have
\beq
v_j \in I^{\mu-\ha}(N^*K_j).
\eeq
The interaction of such conormal waves are analyzed in Section \ref{singu}. However, beyond the first conjugate points, the situation is much more complicated. For example, the distributions $v_j$ may interact at conjugate points or interact many times. These interactions may also produce new singularities which we haven't analyzed yet, and these singularities may affect the ones we analyzed. To avoid such complexities, we follow the approach of \cite{KLU} to consider the interactions only in the following set  
\beq
\begin{gathered}
\mathcal{N}((\vec x, \vec \theta), t_0) = M_0\backslash \cup_{j = 1}^4 J^+(\gamma_{x_j, \theta_j}(\bold{t_j})), \\
\vec x = (x_1, x_2, x_3, x_4), \ \ \vec \theta = (\theta_1, \theta_2, \theta_3, \theta_4),
\end{gathered}
\eeq
i.e.\ away from the causal future of points after the conjugate points.

\subsection{Determination of the conformal class}\label{sec42}
Let $\eps_i, i = 1, 2, 3, 4$  be four small parameters and $\vec \eps = (\eps_1, \eps_2, \eps_3, \eps_4)$. We take $f_{\vec \eps} = \sum_{i = 1}^4 \eps_i f_i \in C_0^4(M_0)$ as constructed in Section \ref{disto} to be the source term in \eqref{eqsem2}. The solution $u_{\vec \eps}$ of the equation \eqref{eqsem2} on $V$ is indeed $u_{\vec \eps} = L(f_{\vec \eps}).$ Here we use subscript $\vec \eps$ to emphasis the dependence on $\vec \eps$. However, in the analysis below, we also use $u = u_{\vec \eps}, f = f_{\vec \eps}$ to simplify the notations.  The asymptotic analysis in Section \ref{semiwave} applies and we denote the fourth order interaction term by
\beq
\mcu^{(4)} = \p_{\eps_1}\p_{\eps_2}\p_{\eps_3}\p_{\eps_4} u|_{\{\eps_1=\eps_2=\eps_3=\eps_4 = 0\}} = \p_{\eps_1}\p_{\eps_2}\p_{\eps_3}\p_{\eps_4}L(f)|_{\{\eps_1=\eps_2=\eps_3=\eps_4 = 0\}}, 
\eeq
see \eqref{eqinter} and \eqref{interms}. Note that this term is determined by $L$. Now we can prove the analogue of Theorem 3.3 and Prop.\ 3.4 in \cite{KLU}, which says that $\mcu^{(4)}$ has new singularities. From our analysis in Section \ref{quad} especially Prop.\ \ref{inter4}, we expect this term to contain singularities due to three wave interactions. To make things clear, we introduce some notations. Let $\pi: T^*M\rightarrow M$ be the standard projection and
\beq
\La^{(3),g} = \cup_{i<j<k} \La^g_{ijk}, \ \ i, j, k = 1, 2, 3, 4.
\eeq
Then we let $\mck^{(3)} = \pi(\La^{(3), g})\subset M$. In particular, $\mck^{(3)}$ is the set in $M$ carrying the singularities produced by three wave interactions. We remark that here we do not know if $K_i$ intersect transversally.

\begin{prop}\label{inter5}
Under the above assumptions and for $s_0>0$ sufficiently small, we have 
\begin{enumerate}
\item If $\bigcap_{j = 1}^4\gamma_{x_j, \theta_j}([0, \bold{t_0}(x_j, \theta_j))) = \emptyset$, i.e. the four geodesics do not intersect before first conjugate points,  then $\mcu^{(4)}$ is smooth in $\mathcal{N}((\vec x, \vec \theta), t_0)$ away from $\mck^{(3)}$ and $K^{(1)}$;
\item If $\bigcap_{j = 1}^4\gamma_{x_j, \theta_j}([0, \bold{t_0}(x_j, \theta_j)))=\{ q_0\}$ and the tangent vectors 
of geodesics $\gamma_{x_j, \theta_j}$ at $q_0$ are linearly independent, then  in $\mathcal{N}((\vec x, \vec \theta), t_0)$ away from $\mck^{(3)}$ and $K^{(1)}$, we have 
 $\mcu^{(4)} \in I^{\mu_0}(\La^g_{q_0}\backslash\La_{q_0})$. Moreover, $\mu_0$ is determined as below
\begin{enumerate}
\item $\mu_0 = 4\mu - \ha$ if $c(q_0)\neq 0$;
\item $\mu_0 = 4\mu - \frac{5}{2}$ if $c = 0$ in a neighborhood of $q_0$ where $b, a$ are non-vanishing;
\item $\mu_0 = 4\mu - \frac{9}{2}$ if $b = c = 0$ in a neighborhood of $q_0$ where $a$ is non-vanishing.
\end{enumerate}
\end{enumerate}
\end{prop}
\bpf

(1) Under the assumptions, for $s_0$ sufficiently small, $\cap_{i = 1}^4K_i$ is empty. We will consider the rest of the cases when $K_i$ may intersect. First, we consider the case where there are three sets intersecting with each other. Without loss of generality, we assume that $K_1, K_2$ and $K_3$ intersect transversally at $K_{123}$. Then $v_4$ is smooth near $K_{123}$ and from the expression of  $\mcu_{1234}^{(4)}$ (see \eqref{u1234}), we know the term is reduced to the triple interactions studied in Section \ref{trip}. From the analysis there, we know that $\WF(\mcu_{1234}^{(4)})\subset \La^{(3),g}\cup \La^{(1)}$. Next, assume that $K_1, K_2, K_3$ intersect at $K_{123}$ but not transversally. In this case, there are no new singularities produced. Actually, we can assume $K_1\cap K_2$ transversally and that
\beq
N^*_q K_3\subset N^*_q K_1 + N^*_q K_2, \ \ q\in K_{123}.
\eeq 
Then we can find the wave front set $\WF(\mcu_{1234}^{(4)})\subset  \La^{(1)}$ by the calculus of wave front sets, see for example \cite[Section 1.3]{Du}. This finishes the proof when there are three sets $K_i$ intersect. 

Finally, if $K_{ijk} = \empty$ for $i<j<k$ i.e.\ there is no three $K_i$ intersect, then the analysis of $\mcu^{(4)}$ is reduced to two wave interactions analyzed in Section \ref{double}. In this case, it is easy to see that $\WF(\mcu^{(4)})\subset \La^{(1)}$. This finishes the proof of (1).

(2) By taking $s_0$ small, we can assume that $q_0$ is the only intersection point. Since the tangent vectors of $\gamma_{x_j, \theta_j}$ at $q_0$ are linearly independent, we see that $K_j$ intersect transversally. All the analysis in Section 3 apply here. The conclusion follows from Prop.\ \ref{inter4}.
\epf

The last ingredient we need for the proof of Theorem \ref{main11} is the important concept of earliest light observation set, see Def. 1.1 of \cite{KLU}. This is used to deal with the problems caused by the conjugate points. Recall that $V$ is a neighborhood of a time-like geodesic $\hat\mu[-1, 1]$. The light observation set of $q\in M$ in $V$ is defined as $\mcp_V(q) = \mcl_q^+\cap V$. The earliest light observation set is 
\beq
\begin{gathered}
\mce_V(q) = \{x\in \mcp_V(q): \text{ there is no $y\in \mcp_V(q)$ and future-pointing time-like path $\alpha:[0, 1] \rightarrow V$} \\
\text{such that $\alpha(0) = y$ and $\alpha(1) = x$}\}\subset V.
\end{gathered}
\eeq
For $W\subset M$ open, the collection of the earliest light observation sets with source points in $W$ is 
\beq
\mce_V(W) = \{\mce_V(q): q\in W\}.
\eeq
In particular, we observe that if $q_0$ is the interaction point as in Prop.\ \ref{inter5}, then 
\beq
\mce_V(q_0) \subset \mathcal{N}((\vec x, \vec \theta), t_0). 
\eeq
This is a consequence of the definition of $\mce_V(q_0)$ and the short cut argument in Section 2.1 of \cite{KLU}. Also, from Section 2.2.1 of \cite{KLU}, we know that $\mce_V(q_0)$ contains a $3$-dimensional submanifold hence is not empty.  

Now we prove the first part of Theorem \ref{main11} that the conformal class can be determined. This is the analogue of Theorem 1.5 of \cite{KLU}.
\bpf[Proof of Theorem \ref{main11}, Part (1)]
The proof now follows the argument  in \cite{KLU}. Our Prop.\ \ref{inter5} and Prop. \ref{nonvan} is almost equivalent to Theorem 3.3 and Proposition 3.4 of \cite{KLU}. However, we need to pay attention to the set in $I(p_-, p_+)$ where the assumptions in Prop.\ \ref{inter5} do not hold. 
It suffices to consider the problem on $(M, g)$ (instead of two manifolds $(M^{(i)}, g^{(i)}), i = 1, 2$ as in the statement of the theorem). Suppose $H(x, z) = a(x)z^2 + b(x)z^3 + c(x)z^4$, where $a, b, c$ are smooth functions in $I(p_-, p_+)$. Consider the set 
\beq
\mcs(c)=\{q\in I(p_-, p_+): c(q)\neq0\} \cup \text{int}(\{q\in I(p_-, p_+): c(q) = 0\}),
\eeq
where $\text{int}(\bullet)$ denotes the interior of the set. We define $\mcs(b)$ similarly. It is easy to see that $\mcs(c), \mcs(b)$ are both dense and open subsets of  $I(p_-, p_+)$. Prop.\ \ref{inter5} (a)(b) holds on $\mcs(c)$ and Prop.\ \ref{inter5} (c) holds on $\mcs(b)\cap \text{int}(\{q\in I(p_-, p_+): c(q) = 0\})$ which is a dense and open subset of $\text{int}(\{q\in I(p_-, p_+): c(q) = 0\})$. Therefore, for $q$ in a dense subset of $I(p_-, p_+)$, there exists $f_i$ so that the term $\mcu^{(4)}$ produces non-trivial singularities on $\mcl_q^+$.  


Consider $\vec x=(x_j)_{j=1}^4,$  $x_j\in V$  and $\vec \theta=(\theta_j)_{j=1}^4,$  where
$\theta_j\in T_{x_j}M$ are future directed light-like vectors. Also, let $t_0>0$ be small enough.
Similarly to  
\cite[Section 3.5]{KLU},
we say that a point $y\in V$,
   satisfies the singularity {\it detection condition} ($D_{}$)  
with light-like directions $(\vec x,\vec \theta)$, 
   and $t_0,\hat s>0$
  if
  \medskip  
  
  \noindent
(${\bf D}_{}$) {\it For any $s,s_0\in (0,\hat s)$ and   $j=1,2,3,4$
there exists $(x_j^{\prime},\theta_j^{\prime})$
 in the
$s$-neighborhood of $(x_j,\theta_j)$, open sets $B_j\subset B_{g^+}(\gamma_{x_j^{\prime},\theta_j^{\prime}}(t_0),s)$, satisfying $B_j\cap J ^+(B_k)=\emptyset$ for 
$j\not =k$, such that the following is valid: There are  ${f}_{j}\in { I}^{\mu+1}(Y((x_j^{\prime},\theta_j^{\prime});t_0,s_0))$ 
such that $\supp(f_j)\subset B_j$, the wavefront set
of $f_j$  is in $s_0$-neighborhood of  $(x_j^{\prime},\theta_j^{\prime})$, and  for the solution 
 $u=u_{\vec \epsilon}$ 
 of (\ref{eqsem2})
with the source ${f}_{\vec \epsilon}=\sum_{j=1}^4
\epsilon_j{f}_{j}$
we have that
$\mcu^{(4)} = \p_{\eps_1\eps_2\eps_3\eps_4}^4 u_{\vec \epsilon}|_{\{\eps_1=\eps_2=\eps_3=\eps_4 = 0\}} $
is not
$C^\infty$-smooth at $y$.}
  \medskip

Our above considerations show that if geodesics $\gamma_{x_j,\theta_j}(\R_+)$ intersect in a point 
$q\in \mathcal{N}((\vec x, \vec \theta), t_0)$, then the set
$$
S(\vec x,\vec \theta,t_0):=\{y\in V:\ \hbox{there is  $\hat s>0$ such that $y$ satisfies (D) with $(\vec x,\vec \theta)$ and 
$t_0,\hat s$}\}
$$
has the property that 
$$
S(\vec x,\vec \theta,t_0)\cap \bigg(\mathcal{N}((\vec x, \vec \theta), t_0) \setminus (\mck^{(3)}\cup \bigcup_{k=1}^4 K_j)\bigg)=
 \mcl_q^+\cap V\cap \bigg(\mathcal{N}((\vec x, \vec \theta), t_0) \setminus (\mck^{(3)}\cup \bigcup_{k=1}^4 K_j) \bigg).
 $$
Roughly speaking, this means that the linearized waves $v_j=Q_g f_j$ interact at the point $q$ and produce a wave $\mcu^{(4)}$  that in the set $\mathcal{N}((\vec x, \vec \theta), t_0)$ may be singular only on the future 
light cone $ \mcl_q^+$  emanating from $q$. Moreover, at any point $y\in  \mcl_q^+\cap  \mathcal{N}((\vec x, \vec \theta), t_0)$
the wave $\mcu^{(4)}$  is surely  non-smooth near $y$  if one makes a suitable perturbation to sources $f_j$.

Next, without loss of generality we can assume that the neighborhood $V$ of the time-like geodesic 
$\hat \mu$ is a union of some time-like geodesics $\mu_a$, $a\in A$.

Define 
 $S_{reg} (\vec x,\vec \theta,t_0)$ be the set  
of the points  $y\in S(\vec x,\vec \theta,t_0)$ having a neighborhood $W\subset U_{g}$  
such that the intersection $W\cap S(\vec x,\vec \theta,t_0)$  
is a non-empty $C^\infty$-smooth 3-dimensional submanifold.  
Moreover, let  $S_{cl} (\vec x,\vec \theta,t_0)$ be the closure of the set $S_{reg} (\vec x,\vec \theta,t_0)$  in $V$
and define  $S_{e} (\vec x,\vec \theta,t_0)$ to be the set of those $y\in S_{cl} (\vec x,\vec \theta,t_0)$
for which any  geodesics $\mu_a$, $a\in A$,   containing $y$  does not intersect $S_{cl} (\vec x,\vec \theta,t_0)$ in the chronological past of $y$.

The proof of Lemma 4.4 of \cite{KLU}, combined with the fact that $\mcs(b)\cap \mcs(c)\subset I(p^-,p^+)$ is an open and dense subset shows
the following: First, in the case when all four geodesics $\gamma_{x_j,\theta_j}(\R_+)$, $j=1,2,3,4$ intersect in some point 
$q\in \mathcal{N}((\vec x, \vec \theta), t_0)$, the above constructed set  $S_{e} (\vec x,\vec \theta,t_0)$ coincides
with $\mathcal E_V(q)$. Second, in the case when all four geodesics $\gamma_{x_j,\theta_j}(\R_+)$ do not intersect 
at any point of
$\mathcal{N}((\vec x, \vec \theta), t_0)$, the above constructed set  $S_{e} (\vec x,\vec \theta,t_0)$ does not intersect 
$\mathcal{N}((\vec x, \vec \theta), t_0)$. Roughly speaking, this means that using the operator $L$ we can construct the earliest light observation sets
corresponding to the intersection point $q$ of any four geodesics $\gamma_{x_j,\theta_j}(\R_+)$, assuming that the four geodesics intersect at a same point $q$ in the set $\mathcal{N}((\vec x, \vec \theta), t_0)$, that is, 
when the intersection point $q$ exists and is before the conjugate points of the geodesics. 
This is the very same conclusion that was made in end of Section 3 of of \cite{KLU}.

By the arguments in Section 4 of \cite{KLU}, we see that the source-to-solution map $L$ determines the earliest light observation sets $\mce_V(q)$ where $q$  runs over the set $I(p_-, p_+)$,
that is, $L$  determines uniquely the collection $\{\mce_V(q): q \in I(p_-, p_+)\}$.
 The problem is thus reduced to the inverse problem with passive measurements.  The conclusion now follows from Theorem 1.2 and Remark 2.2 of \cite{KLU} that the differential structure of $I(p_-, p_+)$ and the conformal class of the metric can be uniquely determined up to diffeomorphisms. This finishes the proof of part (1).
\epf

\subsection{Determination of the nonlinearity and conformal factor}
Next we prove an (simplified) analogue of Theorem \ref{maingh} when $H(x, z)$ only has lower order nonlinearities. Then we complete the proof of Theorem \ref{main11}. 
\begin{theorem}\label{main13}
Let $g^{(1)}, g^{(2)}$ be two globally hyperbolic Lorentzian metrics on a $4$-dimensional manifold $M$. We assume that $g^{(1)} = e^{2\gamma} g^{(2)}$ where $\gamma \in C^\infty(M)$ and $\gamma = 0$ on an open relatively compact set $V \subset M_0 = (-\infty, T_0)\times N, T_0 > 0$.  In particular, $g^{(1)}$ and $g^{(2)}$ are isometric on $V$. Let $\hat \mu$ be a time-like geodesics and $\hat \mu([-1, 1]) \subset V$. Let $p_\pm = \hat \mu(s_\pm)$ with $-1<s_-<s_+<1$.  Consider the semilinear wave equations with source terms
\beq 
\begin{gathered}
\square_{g^{(i)}} u(x) + H^{(i)}(x, u(x))  = f(x), \text{ on } M_0,\\
u = 0 \text{ in } M_0\backslash J^+_{g^{(i)}}(\supp(f)),
\end{gathered}
\eeq
where $\supp(f)\subset V$ and
\beq
H^{(i)}(x, z) = a^{(i)}(x) z^2 + b^{(i)}(x) z^3 + c^{(i)}(x) z^4, \ \ i = 1, 2,
\eeq
where $[a^{(i)}]^2 + [c^{(i)}]^2$ are nowhere vanishing. Let $L^{(1)}, L^{(2)}$ be the source-to-solution map with respect to the metrics $g^{(1)}, g^{(2)}$ and assume that they satisfy
\beq
L^{(1)}(f) = L^{(2)}(f)
\eeq
for all $f$ in a small neighborhood of the zero function in $C_0^4(V)$. We have the following conclusions. 
\begin{enumerate}
\item For $q_0 \in I(p_-, p_+)$, we have $c^{(1)} = e^{\gamma} c^{(2)}$ at $q_0$.
\item If $c^{(1)} = 0$ in a neighborhood of $q_0 \in I(p_-, p_+)$, then $a^{(1)}b^{(1)} = e^{-\gamma} a^{(2)}b^{(2)}$ at $q_0$.
\item If $c^{(1)} = b^{(1)} = 0$ in a neighborhood of $q_0\in I(p_-, p_+)$, then $a^{(1)}  = e^{-\gamma} a^{(2)}$ at $q_0$.
\end{enumerate} 
\end{theorem} 

First of all, we derive some information on the nonlinear term using the order and the principal symbols of $\mcu^{(4)}$. Since the manifold may have caustics, the key point below is to consider $\mcu^{(4)}$ only on the set $\mce_V(q_0), q_0\in I(p_-, p_+)$.  
\begin{prop}\label{main12}
Under the same assumptions as in Theorem \ref{main11}, there is a (conformal) diffeomorphism $\Psi: I(p^{(1)}_-, p^{(1)}_+)\rightarrow I(p^{(2)}_-, p^{(2)}_+)$ such that on $I(p^{(1)}_-, p^{(1)}_+)$ we have  
\begin{enumerate}
\item $c^{(1)}\neq0$ if and only if $\Psi^*c^{(2)}\neq 0$.
\item If $c^{(1)} = \Psi^*c^{(2)}  =0$ in a neighborhood $\mco$ of $q\in I(p^{(1)}_-, p^{(1)}_+)$, then $b^{(1)}\neq 0$ in $\mco$ if and only if $\Psi^*b^{(2)} \neq 0$ in $\mco$. 
\end{enumerate}
\end{prop}
\bpf
It suffices to consider one manifold $(M, g)$. For each $q_0  \in I(p_-, p_+)$, we can find $(x_j, \theta_j) \in T^*V, j = 1, 2, 3, 4$ such that the geodesics $\gamma_{x_j, \theta_j}$ intersect at $q_0$ before their first conjugate points,
that is, $\gamma_{x_j, \theta_j}(t_j)=q_0$ with $0<t_0<{\bf t_0}(x_j, \theta_j)$,
and co-vectors $\xi_j=(\dot \gamma_{x_j, \theta_j}(t_j))^\flat$ satisfy 
 $(\xi_1,\xi_2,\xi_3,\xi_4)\in \mathcal W$, where  $\mathcal W\subset   ({L^{*}_{q_0}M})^4$ is the open and dense 
 set given in Proposition \ref{nonvan}, see also Section 2.2.3 and Section 4 of \cite{KLU}.

By Prop.\ \ref{inter5} and \ref{nonvan}, we can choose the principal symbols of $v_i$ so that the principal symbol of $\mcu^{(4)}$ is non-vanishing on $\mce_V(q_0)\cap (\La^g_{q_0}\backslash (\mck^{(3)}\cup K^{(1)}))$. Recall that $\mck^{(3)}, K^{(1)}$ depend on the parameter $s_0$ and for $s_0\rightarrow 0$, the sets $\mck^{(3)}\cup K^{(1)}$ tend to a set of Hausdorff dimension $2$, while $\mce_V(q_0)\cap \La^{g}_{q_0}$ is of Hausdorff dimension $3$. From \eqref{wavesym}, we know that the principal symbols of $v_i$ and $f_i$ can be determined from each other on $\mce_V(q_0)\cap \La^{g}_{q_0}$. Therefore, we can find $f_i, i = 1, 2, 3, 4$ such that the order of $\mcu^{(4)}$ on $\mce_V(q_0)$ is given by $\mu_0$ in Prop.\ \ref{inter5}. The proposition is finished by comparing the order of the singularities of $\mcu^{(4)}$ corresponding to two semilinear wave equations.
\epf

By examining further the symbols of $\mcu^{(4)}$, we can determine the conformal factor. Here the only missing component is the symbol of $Q_g$ under conformal transformations. The transformation of the wave operators under conformal transformations can be found in Theorem 5.1 of Appendix VI of \cite{Cb} and Appendix A.3 of \cite{Ho}. See also Section 4.6 of \cite{Fr}. Recall that in this paper, we already trivialized the half-density bundles using the volume form on $(M, g)$. 

\begin{prop}\label{symconf}
Let $g, \tilde g$ be two Lorentzian metrics on $M$ such that $g = e^{2\gamma} \tilde g$ where $\gamma \in C^\infty(M).$ Let $Q_g, Q_{\tilde g}$ be the causal inverse of $\square_g, \square_{\tilde g}$ respectively. Then the Lagrangians $\La_g = \La_{\tilde g}$ and the principal symbols of $Q_g, Q_{\tilde g} \in I^{-2}(N^*\diag\backslash \La_g)$ satisfy
\beq
\begin{gathered}
\sigma(Q_g) = e^{2\gamma} \sigma(Q_{\tilde g}).
\end{gathered}
\eeq
For their principal symbols in $I^{-\frac{3}{2}}(\La_g \backslash N^*\diag)$, we have
\beq
\sigma(Q_g)(x, \xi, y, \eta) = e^{-\gamma(x)} \sigma(Q_{\tilde g})(x, \xi, y, \eta) e^{3\gamma(y)},
\eeq
for $(x, \xi), (y, \eta)$ on the same bicharacteristics on $\La_g$.
\end{prop}
\bpf
We first show $\La_g = \La_{\tilde g}$. Let $\mcp(x, \xi) = |\xi|^2_{g^*}$ and $\tilde \mcp(x, \xi) = |\xi|_{\tilde g^*}^2$ be the dual metric function on $T^*M$. Then $\mcp = e^{-2\gamma} \tilde \mcp$. Therefore, the characteristic sets $\Sigma_{g} = \Sigma_{\tilde g}$ and the Hamilton vector fields satisfy
\beq
H_{\mcp} = e^{-2\gamma}H_{\tilde \mcp} + \tilde \mcp H_{e^{-2\gamma}}.
\eeq
Notice that on $\Sigma_g$, we have $H_\mcp = e^{-2\gamma}H_{\tilde \mcp}$. Thus on $\Sigma_g$, the integral curves of $H_\mcp$ and $H_{\tilde \mcp}$ are the same but with different parameterizations. Hence we proved $\La_g = \La_{\tilde g}$. 

Now we know that $Q_g, Q_{\tilde g} \in I^{-\frac{3}{2}, -\ha}(N^*\diag, \La_g)$. We will show that the principal symbols of $
Q_g$ and $e^{-\gamma}Q_{\tilde g}e^{3\gamma}$ are the same on $N^*\diag$ and $\La_g$. First of all, because $\square_g Q_g = \id$, we know that $\sigma(\square_g)\sigma(Q_g) = 1$ on $N^*\diag$. Thus $\sigma(Q_g) = \mcp^{-1}$ on $N^*\diag$. Similarly, $\sigma(Q_{\tilde g}) = \tilde \mcp^{-1} = e^{2\gamma} \mcp^{-1}$. Thus we proved $\sigma(Q_g) = e^{2\gamma}\sigma(Q_{\tilde g})$ on $N^*\diag$. 

Next, consider the principal symbols on $\La_g$. According to the proof of Prop.\ 6.6 in \cite{MU}, we know that $\sigma(Q_g)$ on $\La_g$ satisfies 
\beqq\label{eqlie}
H_\mcp \sigma(Q_g) + i\sigma_{sub, g}(x, \xi) \sigma(Q_g) = 0,
\eeqq
where 
\beq
\sigma_{sub, g}(x, \xi) = - \frac{1}{2i} \sum_{i, j = 1}^4 \frac{\p^2 \mcp(x, \xi)}{\p x_i \p \xi_j}  
\eeq
is the sub-principal symbol of $\square_g$ on $\La_g$. See also \cite[Prop.\ 4.3.1]{Du}. Note here we trivialized the half-density factors and write the Lie derivative as $H_\mcp$. This is an ordinary differential equation along the integral curves of $H_\mcp$. The initial condition is determined by $\sigma(Q_g)$ at $N^*\diag$, see \cite[Prop.\ 6.6]{MU}. Therefore, to see the principal symbols are the same, we only need to show that $\sigma(e^\gamma Q_{\tilde g} e^\gamma)$ satisfies \eqref{eqlie}. This follows from the following computations.
\beq
\begin{split}
H_{\tilde \mcp} \sigma(Q_{\tilde g}) + i\sigma_{sub, \tilde g}(x, \xi) \sigma(Q_{\tilde g}) &= e^{2\gamma}H_\mcp  \sigma(Q_{\tilde g}) -  \frac{1}{2} \sum_{i, j = 1}^4 [e^{2\gamma} \frac{\p^2 \mcp(x, \xi)}{\p x_i \p \xi_j} + \frac{\p e^{2\gamma}}{\p x_i} \frac{\p \mcp(x, \xi)}{\p \xi_j}]\sigma(Q_{\tilde g})\\
& = e^{3\gamma} [ H_\mcp (e^{-\gamma} \sigma(Q_{\tilde g})) -  \frac{1}{2} \sum_{i, j = 1}^4 \frac{\p^2 \mcp(x, \xi)}{\p x_i \p \xi_j} (e^{-\gamma}\sigma(Q_{\tilde g}))]\\
& = e^{3\gamma}[H_\mcp (e^{-\gamma} \sigma(Q_{\tilde g})) + i\sigma_{sub, g} (e^{-\gamma} \sigma(Q_{\tilde g}))].
\end{split}
\eeq
This proves that $\sigma(Q_g)$ and  $\sigma(e^{-\gamma} Q_{\tilde g} e^{3\gamma})$ satisfy the same equation on $\La_g$. Since we proved the symbols are the same on $N^*\diag$, by solving the transport equations, we see that the symbols are the same on $\La_g$. This finishes the proof.
\epf

Now we prove the relation of the conformal factor and the nonlinear terms.
\bpf[Proof of Theorem \ref{main13}]
Although we assumed $g^{(1)}$ and $g^{(2)}$ are isometric on $V$, this actually follows by linearizing the source-to-solution map, see Remark 3.1 of \cite{KLU}.  Next, we assume $f = \sum_{i = 1}^4\eps_i f_i$ constructed as in Section \ref{sec42} and denote
\beq
\mcu^{(4), j} = \p_{\eps_1}\p_{\eps_2}\p_{\eps_3}\p_{\eps_4}L^{(j)}(f)|_{\{\eps_1=\eps_2=\eps_3=\eps_4 =0\}}, \ \ j = 1, 2. 
\eeq
For any $q_0 \in I(p_-, p_+)$, we will compare the principal symbols of $\mcu^{(4), 1}$ and $\mcu^{(4), 2}$ on $\mce_{V}(q_0)\backslash (\mck^{(3)}\cup K^{(1)})$ for $s_0\rightarrow 0$ (see the proof of Prop.\ \ref{main12}) using the computations  in Section \ref{nonsym}. We remark that since conformal transformations of Lorentzian metrics preserves light-like (pre)geodesics, the sets $\mce_V(q_0)$ are the same for $g^{(1)}, g^{(2)}$. 

(1): If $c^{(1)}(q_0)\neq0$, from Prop.\ \ref{main12} we know that $c^{(2)}(q_0)\neq0$. For $(q, \eta)\in \mce_V(q_0)$ which is joined with $(q_0, \xi) \in L^{*,+}M$ and $q_0\in I(p^{(1)}_-, p^{(1)}_+)$, we know the principal symbol of $\mcu^{(4), 1}(q, \eta)$ from Prop.\ \ref{nonvan}. In particular, we can write 
\beq
\sigma(\mcu^{(4), 1})(q, \eta) = -4! (2\pi)^{-3}\sigma(Q_{g^{(1)}})(q, \eta, q_0, \xi)c^{(1)}(q_0)\prod_{i = 1}^4A_i(q_0, \xi_i),
\eeq
where $A_i$ are the principal symbols of $v_i$, which satisfies
\beq
A_i(q_0, \xi_i) = \sigma(Q_{g^{(1)}})(q_0, \xi_i, x_i, \zeta_i)B_i(x_i, \zeta_i), \ \ i = 1, 2, 3, 4,
\eeq
where $x_i\in V$, $(q_0, \xi_i)$  and $(x_i, \zeta_i)$ are joined by bicharacteristics and $B_i$ are the principal symbols of $f_i$. Also, $\sigma(\mcu^{(4), 2})$ has a similar expression by changing $c^{(1)}$ to $c^{(2)}$ and $Q_{g^{(1)}}$ to $Q_{g^{(2)}}$. By Prop.\ \ref{symconf}, we have
\beqq\label{symrel}
\begin{gathered}
\sigma(Q_{g^{(1)}})(q, \eta, q_0, \xi) = \sigma(Q_{g^{(2)}})(q, \eta, q_0, \xi)e^{3\gamma(q_0)},\\
\sigma(Q_{g^{(1)}})(q_0, \xi_i, x_i, \zeta_i) = e^{-\gamma(q_0)}\sigma(Q_{g^{(2)}})(q_0, \xi_i, x_i, \zeta_i).
\end{gathered}
\eeqq
Therefore, we obtain the following relation
\beq
\sigma(\mcu^{(4), 1})(q, \eta)c^{(2)}(q_0) e^{\gamma(q_0)} = c^{(1)}(q_0) \sigma(\mcu^{(4), 2})(q, \eta).
\eeq
By the discussion in Section \ref{nonsym}, we know that $\sigma(\mcu^{(4), i}), i = 1, 2$ satisfy the same relation. But from $L^{(1)}(f) = L^{(2)}(f)$, we know the principal symbols of $\mcu^{(4), i}$ should be the same hence we proved (1).

(2): From Prop.\ \ref{nonvan}, we can find the principal symbol of $\mcu^{(4), i}, i = 1, 2$ and by \eqref{symrel}, they should satisfy
\beq
\sigma(\mcu^{(4), 1})(q, \eta) b^{(2)}(q_0)a^{(2)}(q_0) = b^{(1)}(q_0)a^{(1)}(q_0) e^{(3-4)\gamma(q_0)} e^{2\gamma(q_0)} \sigma(\mcu^{(4), 2})(q, \eta).
\eeq
As the source-to-solution maps are the same, the symbols are the same and we proved part (2).

(3): By Prop.\ \ref{nonvan} and the relations \eqref{symrel}, we obtain the following relation
\beq
\sigma(\mcu^{(4), 1})(q, \eta) [a^{(2)}(q_0)]^{3}= [a^{(1)}(q_0)]^3 e^{(3-4)\gamma(q_0)} e^{4\gamma(q_0)} \sigma(\mcu^{(4), 2})(q, \eta).
\eeq
As the source-to-solution maps are the same, the symbols are the same. So we proved (3) and completed the proof of the theorem.
\epf

Finally, we finish the proof of Theorem \ref{main11}.
\bpf[Proof of Theorem \ref{main11}, Part (2)]
Using the result in part (1) and Remark 3.1 of \cite{KLU}, we know that $g^{(1)} = \Psi^*g^{(2)}$ on $V^{(1)}\cap I(p^{(1)}_-, p^{(1)}_+)$ i.e.\ $\gamma = 0$ there. This is because the derivative of the source-to-solution map $f\rightarrow L(f)$ of semilinear wave equations determines the source-to-solution map of the linearized wave equation which determines the metric. From Theorem \ref{main13}, we know that the coefficients $a^{(1)} = a^{(2)}, b^{(1)} = b^{(2)}$ and $c^{(1)} = c^{(2)}$ because $e^\gamma = 1$ on $V^{(1)}\cap I(p^{(1)}_-, p^{(1)}_+)$. It follows from Theorem \ref{main13} again that $e^{\gamma} = 1$ on a dense subset of $I(p^{(1)}_-, p^{(1)}_+)$ (see the proof of  part (1) of Theorem \ref{main11}). By the continuity of $\gamma$, we see that $e^{\gamma}=1$ in $I(p_-, p_+)$. This completes the proof.
\epf

\section{Singularities in higher order asymptotic expansions}\label{singuhigh}
Now we return to the full generality of the nonlinear term $H$, i.e.\ we assume $H(x, z)$ is genuinely nonlinear. As in the lower order nonlinearity case, we will make use of the singularities generated by nonlinear interactions, but in higher order terms of the asymptotic expansion. The analysis bears some similarities with the classical treatment of interaction of conormal singularities, see for example \cite[Theorem 4.1]{Bea}. However, we emphasis the difference is that we analyze singularities of every term of the asymptotic expansion of $u$ rather than $u$ itself.  

\subsection{The asymptotic expansion}
Consider the semilinear wave equation  
\beq
\begin{gathered}
\square_g u(x) + H(x, u(x))  = f(x), \text{ in } M_0\\
u = 0 \text{ in } M_0\backslash J^+(\supp(f)),
\end{gathered}
\eeq
where $M_0 = (-\infty, T_0)\times N, T_0 > 0$ and $H$ is genuinely nonlinear. We continue with the same notations and assumptions as in Section \ref{semiwave}. For $k\geq 4$, we need the following terms in the asymptotic expansion of the solution $u$ as $\eps_i \rightarrow 0, i = 1, 2, 3 ,4$,
\beqq\label{eqmcuk}
\mcu^{(k)} = \p^{k-3}_{\eps_1}\p_{\eps_2}\p_{\eps_3}\p_{\eps_4} u|_{\{\eps_1 = \eps_2 = \eps_3 = \eps_4 = 0\}}.
\eeqq
For $k = 4$, this term was computed explicitly in Section \ref{semiwave}. Here we focus on the cases when $k\geq 5$. We start from the equation
\beq
u = v - Q_g(H(x, u)).
\eeq
We expand $H$ into Taylor series in $z$
\beq
H(x, z) = \sum_{j = 2}^k h_j(x) z^j + O(z^{k+1}) \text{ for $z$ small}. 
\eeq
Then we have
\beqq\label{eqitea}
u = v - Q_g(h_k(x) u^k) - Q_g(\sum_{j = 2}^{k-1} h_j(x) u^j) + \mcr,
\eeqq
where $\mcr$ denotes terms which are in $H^{4}(M)$ but not of the order $\eps_1^{k-3}\eps_2\eps_3\eps_4$. We shall ignore such terms later. By iterating \eqref{eqitea}, we obtain
\beq
\begin{gathered}
u = v - Q_g(h_k(x) v^k) - Q_g(\sum_{\alpha_1 = 2}^{k-1} h_{\alpha_1}(x) (v - Q_g(\sum_{{\alpha_2} = 2}^{k-1} h_{\alpha_2}(x) u^{\alpha_2}))^{\alpha_1}) + \mcr\\
 = v - Q_g(h_k(x) v^k) - Q_g(\sum_{\alpha_1 = 2}^{k-1} h_{\alpha_1}(x) (v - Q_g(\sum_{{\alpha_2} = 2}^{k-1} h_{\alpha_2}(x) (v - Q_g(\sum_{{\alpha_3} = 2}^{k-1} h_{\alpha_3}(x) u^{\alpha_3}))^{\alpha_2}))^{\alpha_1}) + \mcr.
\end{gathered}
\eeq
Continuing the iteration, we get that
\beqq\label{eqmcug}
\begin{gathered}
u  = v - Q_g(h_k(x) v^k) - \mcu + \mcr, \text{ where } \\
\mcu = Q_g(\sum_{\alpha_1 = 2}^{k-1} h_{\alpha_1}(x) (v - Q_g(\sum_{{\alpha_2} = 2}^{k-1} h_{\alpha_2}(x) (v - \cdots (v - Q_g(\sum_{{\alpha_{k-1}} = 2}^{k-1} h_{\alpha_{k-1}}(x) v^{\alpha_{k-1}}))^{\alpha_{k-2}} \cdots )^{\alpha_2}))^{\alpha_1}).
\end{gathered}
\eeqq
We observe that in general, $\mcu^{(k)}$ defined in \eqref{eqmcuk} consists of many terms which are hard to write down. Instead of analyzing each term as we did in Section \ref{singu}, we write  
\beqq\label{highnon}
\begin{gathered}
\mcu^{(k)} = \mcu^{(k)}_0 + \mcu^{(k)}_1, \\
\text{ where } \mcu^{(k)}_0 = -\binom{k}{2} Q_g(h_k(x)v_1^{k-3}v_2v_3v_4),
\end{gathered}
\eeqq
and we will show that the terms in $\mcu_1^{(k)}$ which come from $\mcu$ are smoother than $\mcu^{(k)}_0$ near $\La_{q_0}^g$.
Thus, we can identify the leading singularities in $\mcu^{(k)}, k\geq 4$. This idea is the same as in Theorem 4.1 of \cite{Bea}. For this purpose, we will need some microlocal function spaces to describe the regularities conveniently. 

\subsection{Leading singularities and principal symbols}\label{sechs}
The leading term $\mcu^{(k)}_0$ in \eqref{highnon} contains powers of conormal distributions. We deal with this first.  Recall that the bundle of light-like co-vectors $L^*M$ is the (disjoint) union of $L^{*,+}M$ and $L^{*,-}M$. Thus, if $K\subset M$ and $N^*K\subset L^*M$, the Lagrangian submanifold $N^*K$ is the union of two disjoint Lagrangian submanifolds $N^*K\cap L^{*, +}M$ and $N^*K\cap L^{*, -}M$.  In this subsection, we shall assume that the distorted plane wave $v_1\in I^{\mu-\ha}(K_1)$ and $\WF(v_1)\subset N^*K_1\cap L^{*, +}M$.

\begin{lemma}\label{lm50}
Let $K\subset M$ be a codimension one submanifold. Let $u\in I^\mu(K), v\in I^{\mu}(K)$ and $\WF(u), \WF(v)\subset N^*K \cap L^{*,+}M$. Then for $\mu<-\frac{3}{2}$, $uv$ is a well-defined distribution in $I^{\mu}(K)$. Moreover, the principal symbols satisfy
\beq
\sigma(uv) = (2\pi)^{-\ha}\sigma(u)\ast\sigma(v).
\eeq
Here the convolution is over the fiber variable of $N^*K$ i.e.\
\beq
\sigma(u)\ast\sigma(v)(x, \xi) = \int_{\mbr} \sigma(u)(x, \xi-\eta) \sigma(v)(x, \eta)d\eta, \ \ (x, \xi), (x, \eta) \in N^*K.
\eeq
\end{lemma}
\bpf
Since the wave front set of $u$ and $v$ are future-pointing light like vectors, the multiplication $uv$ is well-defined. Next, choose local coordinates $x = (x^1, x^2, x^3, x^4)$ such that $K = \{x^1 = 0\}$. Let $\xi$ be the dual variable to $x$. We can write
\begin{gather*}
u(x) = \int_{\mbr} e^{ix^1\xi_1} A(x^2, x^3, x^4; \xi_1) d\xi_1, \ \ v(x) = \int_{\mbr} e^{ix^1\eta_1} B(x^2, x^3, x^4; \eta_1) d\eta_1,
\end{gather*}
where $A, B\in S^{\mu + \fnf - \ha}(\mbr^3; \mbr)$. Therefore,
\begin{gather*}
u(x)v(x) = \int_{\mbr^2} e^{ix^1(\xi_1 + \eta_1)} A(x^2, x^3, x^4; \xi_1)B(x^2, x^3, x^4; \eta_1) d\xi_1 d\eta_1 = \int_{\mbr} e^{ix^1\zeta_1} C(x^2, x^3, x^4; \zeta_1) d\zeta_1,
\end{gather*}
where $\zeta_1 = \xi_1 + \eta_1$ and 
\beq
C(x^2, x^3, x^4; \zeta_1) = \int_{\mbr} A(x^2, x^3, x^4; \xi_1)B(x^2, x^3, x^4; \zeta_1 - \xi_1) d\xi_1 = A\ast B,
\eeq
where the convolution is in the bundle variables. We will show that $C\in S^{\mu + \fnf - \ha}(\mbr^3; \mbr)$ if $\mu<-\frac{3}{2}$.  Actually, for $|\zeta_1 - \xi_1| < |\xi_1|$, we get $|\zeta_1| < 2|\xi_1|$. Then we get 
\beq
|\langle \zeta_1\rangle^{-\mu-\fnf+ \ha} C(x^2, x^3, x^4, \zeta_1)| \lesssim \int_{\mbr} \langle\zeta_1 - \xi_1\rangle^{\mu+\fnf- \ha} d\xi_1< C_0,
\eeq
where $C_0$ denotes a generic constant. On the other hand, if $|\zeta_1 - \xi_1| > |\xi_1|$, we get $|\zeta_1| < 2|\zeta_1- \xi_1|$ and 
\beq
|\langle \zeta_1\rangle^{-\mu-\fnf+ \ha} C(x^2, x^3, x^4, \zeta_1)| \lesssim \int_{\mbr}  \langle\xi_1\rangle^{\mu+\fnf- \ha} d\xi_1 < C_0.
\eeq
Therefore, $uv\in I^{\mu}(K)$. The principal symboles is $C$ modulo lower order terms in $S^{\mu + \fnf - \ha-1}(\mbr^3; \mbr)$
\epf

Now we introduce a function space specially designed to describe the singularities in the interaction of four conormal waves $v_i\in I^\mu(K_i), i = 1, 2, 3, 4$ where $K_i$ intersect transversally at $q_0$, see Def. \ref{trans4}. For convenience, we define 
\beq
\Theta = \cup_{i = 1}^3 \La^{(i)}\cup \La_{q_0}.
\eeq
We have seen in Section \ref{singu} that the wave front set of the multiplication of $v_i$ is contained in $\Theta(\eps)$ for any $\eps$ small and after the action of $Q_g$, the wave front set is contained in $\Theta(\eps)\cup \Theta^g(\eps)$. Also, we stayed away from the singularities on 
\beq
\mck \doteq K^{(1)}\cup \mck^{(3)},
\eeq
which are possibly stronger than the singularities on $\La_{q_0}^g$. We shall consider the space of distributions with these  properties.  

\begin{definition}
For $s >2 $, we define a microlocal function space $\hml^s(M)$ consisting of distributions $u$ such that $\WF(u)\subset \Theta(\eps)\cup \Theta^g(\eps)$ for $\eps > 0$ sufficiently small and $u\in H^s_{\loc}(M\backslash \mck)$.
\end{definition} 

In particular, for $v_i\in I^{\mu}(K_i)$, we have $v_i \in \hml^{-\infty}(M)$. We prove some properties for $\hml^s(M)$.
\begin{lemma}\label{lm51}
For $s > \frac{n+1}{2} = 2$, $\hml^s(M)$ is an algebra. 
\end{lemma}
\bpf
Let $u, v\in \hml^s(M)$. Away from $\mck$, $u, v\in H^s_{\loc}(M)$ which is an algebra for $s>2$. By the calculus of wave front set (see e.g.\ \cite[Section 1.3]{Du}), we know that $\WF(uv)\subset \Theta(\eps)\cup \Theta^g(\eps)$ for some $\eps>0$ sufficiently small. Thus $uv\in \hml^s(M)$.
\epf

To understand the action of $Q_g$ on $\hml^s(M)$, we recall the classical propagation of singularities for operators of real principal type due to H\"ormander, see e.g.\ Theorem 26.1.4 of \cite{Ho4}. We restate the theorem for $\square_g$ which is all we need and we use the $H^s$ wave front set $\WF^s(u)$ of $u$. By definition, $\WF^s(u)$ is  the complement of the set consisting of $(x_0, \zeta_0)\in T^*M\backslash 0$ such that there exists a conic neighborhood $\Gamma$ of $(x_0, \zeta_0)$ such that $\chi_\Gamma u\in H^s(M)$ for any cut-off function $\chi$ which is supported in $\Gamma$.
\begin{theorem}\label{hopro}
Let $\square_g$ be the Laplace-Beltrami operator on $(M, g)$ and $\mcp$ be the principal symbol of $\square_g$, see Section \ref{seclinwa}. Let $(x_0, \zeta_0)\in T^*M\backslash 0$ such that $\mcp(x_0, \zeta_0) = 0$ and $\gamma$ be the null bicharacteristic through $(x_0, \zeta_0)$. If $\square_g u = f$, $\gamma \cap \WF^{s-1}(f) = \emptyset$ and $(x_0, \zeta_0)\notin \WF^s(u)$, then $\gamma\cap \WF^s(u) = \emptyset.$
\end{theorem}

\begin{lemma}\label{lm52}
Let $Q_g$ be the causal inverse of $\square_g$. We have 
\beq
Q_g: \hml^s(M)\rightarrow \hml^{s+1}(M).
\eeq
\end{lemma}
\bpf
Suppose $f\in \hml^s(M)$. The wave front set of $Q_g(f)$ is in $\Theta(\eps)\cup \Theta^g(\eps)$. Let $\chi$ be a microlocal smooth cut-off function supported away from $\La^{1}(\eps)\cup \La^{(3)}(\eps)\cup \La^{(3), g}(\eps)$ for $\eps>0$ small. Let $\text{Op}(\chi)$ be a pseudo-differential operator with principal symbol $\chi$. Then $\text{Op}(\chi) f\in H^s_{\loc}(M)$ and $Q_g\circ\text{Op}(\chi) f \in H^{s+1}_{\loc}(M)$ by Prop.\ \ref{propqg}. However, by the propagation of singularities (Theorem \ref{hopro}), $Q_g\circ(\id - \text{Op}(\chi))(f) \in H^{s+1}(M\backslash \mck)$. This finishes the proof.
\epf

Now we  analyze the singularities in $\mcu^{(k)}$.
\begin{prop}\label{highasym}
For any $k\geq 4$, consider the term 
\beq
\mcu^{(k)} = \p^{k-3}_{\eps_1}\p_{\eps_2}\p_{\eps_3}\p_{\eps_4} u|_{\{\eps_1 = \eps_2 = \eps_3 = \eps_4 = 0\}}.
\eeq 
We can write 
\beq
\mcu^{(k)} = \mcu^{(k)}_0 + \mcu^{(k)}_1, \text{ where }\mcu^{(k)}_0 = -\binom{k}{2}Q_g(h_k(x)v_1^{k-3}v_2v_3v_4)
\eeq
such that away from $\mck$ 
\beq
\mcu^{(k)}_0 \in I^{4\mu + \frac{3}{2}}(\La_{q_0}^g\backslash \La_{q_0}), \ \ \mcu^{(k)}_1 \in H_{\loc}^{-4\mu - 2}(M). 
\eeq
\end{prop}

\begin{proof}
We start with $\mcu^{(k)}_0$. By Lemma \ref{lm50}, we know that $v_1^{k-3} \in I^\mu(K_1)$. Therefore, by Prop.\ \ref{inter4}, we know that $\mcu^{(k)}_0 \in I^{4\mu + \frac{3}{2}}(\La_{q_0}^g\backslash \La_{q_0})$ away from $\mck$. It remains to show the regularity of $\mcu_1^{(k)}$. Since $\mcu_{1}^{(k)}$ is contained in $\mcu$ in \eqref{eqmcug}, it suffices to understand the regularity of $\mcu$ given by
\beq
\begin{gathered}
\mcu = Q_g(\sum_{\alpha_1 = 2}^{k-1} h_{\alpha_1}(x) (v - Q_g(\sum_{{\alpha_2} = 2}^{k-1} h_{\alpha_2}(x) (v - \cdots (v - Q_g(\sum_{{\alpha_{k-1}} = 2}^{k-1} h_{\alpha_{k-1}}(x) v^{\alpha_{k-1}}))^{\alpha_{k-2}} \cdots )^{\alpha_2}))^{\alpha_1}).
\end{gathered}
\eeq

Recall that $v = \sum_{i = 1}^4\eps_i v_i, v_i\in I^\mu(K_i).$ Therefore,  with $s= -\mu- 1$, we know from Lemma \ref{lm50}, \ref{lm51}, \ref{lm52} that 
\beq
v^l \in \hml^{4s}(M), \ \ Q_g(v^l) \in \hml^{4s+1}(M), \ \ \forall l\geq 0.
\eeq
By the algebraic properties of the spaces, we know that 
\beq
\mcu = Q_g(\sum_{\alpha_1 = 2}^{k-1} h_{\alpha_1}(x) (v - \tilde v)^{\alpha_1}),  \ \ \tilde v \in \hml^{4s+1}(M). 
\eeq
Finally, we use Lemma \ref{lm51} and Lemma \ref{lm52} again to conclude that $\mcu \in \hml^{4s+2}(M)$. This finishes the proof.
\epf

The point of this lemma is that if the leading singularities of $\mcu_{0}^{(k)}$ is non-vanishing on $\La^g_{q_0}\backslash (\La^{(3), g}\cup \La^{(1), g})$ (see \eqref{eqlak} and \eqref{flowset}), then $\mcu_{0}^{(k)} \notin H^{-4\mu-2}_{\loc}(M)$ and we get
\beq
\mcu_0^{(k)} = \mcu^{(k)} \text{ mod }  H^{-4\mu -2}_{\loc}(M\backslash \mck).
\eeq
In particular, we can separate the leading singularities.

Finally, we compute the principal symbols of $\mcu^{(k)}_0$ and we will see that this involves the Taylor coefficients $h_k(x)$. We first express the principal symbols of $v_1^{k-3}$. By Lemma \ref{lm51}, this is 
\beq
\sigma(v_1^{k-3})  = (2\pi)^{-\ha(k-4)}\sigma(v_1)\ast^{(k-4)}\sigma(v_1) = (2\pi)^{-\ha(k-4)} \underbrace{\sigma(v_1)\ast \sigma(v_1)\ast \cdots \ast \sigma(v_1)}_{k-4 \text{ convolutions }}.
\eeq
Here as in Lemma \ref{lm51}, all the convolutions are over the fiber variables. Before we continue, we make an observation when this is non-vanishing at a given point $(q_0, \zeta_1)$. Suppose $\sigma(v_1)(q_0, \zeta_1) > 0$ and $\sigma(v_1) \geq 0$ on its support. Then the convolution $\sigma(v_1)\ast \sigma(v_1)$ is positive at $(q_0, \zeta_1)$. This argument can be continued so that $\sigma(v_1^{k-3})$ is non-vanishing at $(q_0, \zeta_1)$.

For any $\zeta\in \La_{q_0}\backslash (\La^{(3)}\cup \La^{(1)})$, we can write $\zeta = \sum_{i = 1}^4\zeta_i$ where $\zeta_i\in N^*_{q_0}K_i$. Also, we let $A_i$ be the principal symbols of $v_i, i = 2, 3, 4$. We let $A^{(k-3)}_1 = \sigma(v_1^{k-3})$. Then we have
\beq
\sigma_{\La_{q_0}}(h_k(x)v_1^{k-3}v_2v_3v_4)(q_0, \zeta)  = (2\pi)^{-\ha(k-4)-3} h_k(q_0)A^{(k-3)}_1(q_0, \zeta_1)\prod_{i = 2}^4A_i(q_0, \zeta_i).
\eeq
By Prop.\ 2.1 of \cite{GrU93}, we get that 
 \beq
 \sigma_{\La_{q_0}^g}(\mcu_0^{(k)})(q, \eta)  = (2\pi)^{-\ha(k-4)-3} \sigma_{\La_g}(Q_g)(q, \eta, q_0, \zeta)h_k(q_0)A^{(k-3)}_1(q_0, \zeta_1)\prod_{i = 2}^4A_i(q_0, \zeta_i),
\eeq
 where $(q, \eta)$ is joined with $(q_0, \zeta)$ by bicharacteristics.  The symbols satisfy
\beq
A_i(q_0, \xi_i) = \sigma(Q_{g})(q_0, \xi_i, x_i, \zeta_i)B_i(x_i, \zeta_i), \ \ i = 1, 2, 3, 4,
\eeq
where $x_i\in V$, $(q_0, \xi_i)$  and $(x_i, \zeta_i)$ are joined by bicharacteristics and $B_i$ are the principal symbols of $f_i = \square_g v_i, i = 1,2, 3, 4$. Therefore, we can choose $f_i$ such that the leading singularities in $\mcu^{(k)}$ is non-vanishing at $q_0$.

\subsection{Multiplication of distributions}\label{multising}
In the previous subsection, we analyzed the leading singularities in $\mcu^{(k)}_0, k\geq 4$. In principal, we can analyze every term in $\mcu^{(k)}$ by the same method in Section \ref{singu}. For example, we will use the singularities of $\mcu^{(5)}$ in Section \ref{pfmain}. These terms involve multiplication of Lagrangian distributions whose wave front sets intersects at $\La_1$. We analyze them in this subsection. All the proofs below follow the same ideas as in Section \ref{singu}. We continue using the notations in Section \ref{singu}, especially $\La_{\bullet}$ and $K_{\bullet}$ in \eqref{lakj}.

\begin{lemma}\label{prod51}
Let $u\in  I^{\mu}(\La_1), v \in  I^{\mu'}(\La_{q_0})$ and $\mu, \mu'<-\frac{3}{2}$. For any $\eps > 0$, we can write $w = uv  = w_0 + w_1$ such that $w_0\in I^{\mu'}(\La_{q_0})$ and $w_1 \in \mcd'(M; \La_1(\eps))$. Moreover, the principal symbol of $w$ satisfies
\beq
\sigma_{\La_{q_0}}(w_0)= (2\pi)^{-\ha}\sigma_{\La_{1}}(u)\ast_1 \sigma_{\La_{q_0}}(v),
\eeq
where $\ast_1$ denotes the partial convolution in the fiber variable of $N^*K_1$ i.e.\
\beq
\sigma_{\La_{1}}(u)\ast_1 \sigma_{\La_{q_0}}(v)(x, \xi) = \int_\mbr \sigma_{\La_{1}}(u)(x, \xi-\eta) \cdotp \sigma_{\La_{q_0}}(v)(x, \eta) d\eta, \ \ (x, \xi), (x, \eta)\in N^*K_1.
\eeq 
\end{lemma}
\bpf
Thanks to Lemma \ref{symp}, it suffices to work with the model Lagrangians $\tilde \La_i$. We can write 
\beq
u = \int_{\mbr} e^{ix^1 \xi_1} A(x^2, x^3, x^4; \xi_1)d\xi_1,
\eeq
with $A\in S^{M}(\mbr^3_{x^2, x^3, x^4}; \mbr_{\xi_1})$ being standard symbol and where $M = \mu -\fnf + \ha$. Also, we can write 
\beq
v = \int_{\mbr} e^{ix\cdot \xi} B(x; \xi)d\xi, 
\eeq
with $B\in S^{m}(\mbr^4; \mbr^4)$ being a standard symbol and $m = \mu'-\fnf + 2$. 
Their product is
\beq
\begin{gathered}
 uv = \int_{\mbr^2} e^{i (x^1 \eta_1 + x\cdot \xi)} A(x^2, x^3, x^4, \eta_1)B(x, \xi)d\xi d\eta_1 = \int_{\mbr} e^{ix\cdot \zeta} C(x, \zeta) d\zeta,\\
 \text{ where } C(x, \zeta) = \int_{\mbr} A(x^2, x^3, x^4, \eta_1)B(x, \zeta_1 - \eta_1, \zeta_2, \zeta_3, \zeta_4) d\eta_1 = A\ast_1 B.
 \end{gathered}
\eeq
Here we use $\ast_1$ to denote the partial convolution in the $\eta_1$ variable. In a coordinate invariant way, this is the partial convolution in the bundle variable of $N^*K_1$. Similar to the proof of Lemma \ref{prod3} and \ref{prod41}, we introduce a cut-off function to stay away from $\tilde \La_1$. This proves the lemma.
\epf

\begin{lemma}\label{prod52}
Let $u\in  I^{p, l}(\La_{12}, \La_2), v \in  I^{p, l}(\La_{13}, \La_3)$ with $p +l < -\frac{3}{2}$. For $\eps>0$ sufficiently small, we can write  $w = uv$ as
\beq
\begin{gathered}
w = w_0 + w_1+ w_2, \ \ w_0\in I^{p+ l -\ha}(\La_{123}), \ \ w_1\in \mcd'(M; \La^{(1)}(\eps)) \\
 w_2 \in I^{p, p+1}(\La_{23}, \La_2) + I^{p, p+1}(\La_{23}, \La_3).
 \end{gathered}
\eeq
Moreover, the principal symbol of $w_0$ satisfies
\beq
\sigma_{\La_{123}}(w_0) = (2\pi)^{-\ha}\sigma_{\La_{12}}(u)\ast_1 \sigma_{\La_{13}}(v).
\eeq
\end{lemma}

\bpf
It suffices to consider the distributions on model Lagrangians $\tilde \La_i$, see Lemma \ref{symp}. We  write 
\beq
u = u_0+u_1, \ \ u_0 = \int_{\mbr^2} e^{i(x^1 \xi_1 + x^2\xi_2)} A(x^3, x^4; \xi_1, \xi_2)d\xi_1d\xi_2,\ \ u_1\in I^p(\tilde \La_2),
\eeq
with $A\in S^{M, M'}(\mbr^2_{x^3, x^4}; \mbr_{\xi_2}, \mbr_{\xi_1})$ being a symbol of product type where $M = p-\fnf + 1 + 1, M' = l-1$. Also, $u_1$ is supported away from $K_{123}.$ Next we write 
\beq
v = v_0+v_1, \ \ v_0 = \int_{\mbr^2} e^{i(x^1 \xi_1 + x^3\xi_3)} B(x^2, x^4; \xi_1, \xi_3)d\xi_1d\xi_3, \ \ v_1 \in I^{p'}(\tilde \La_3),
\eeq
with $B\in S^{m, m'}(\mbr^2_{x^2, x^4}; \mbr_{\xi_1}, \mbr_{\xi_3})$ being a symbol of product type where $m = p-\fnf +1 + 1, m' = l-1$. Also, $v_1$ is supported away from $K_{123}$. 
The product of $u_0v_1$ and $u_1v_0$ are smooth if we arrange the supports. We have
\beq
\begin{gathered}
 v_1u_1 \in I^{p, p+1}(\La_{23}, \La_2) + I^{p, p+1}(\La_{23}, \La_3).
 \end{gathered}
\eeq
It remains to find the product
\beq
 u_0v_0 = \int_{\mbr^4} e^{i(x^1 (\xi_1+ \eta_1) + x^2\xi_2 + x^3\xi_3)} A(x^3, x^4; \eta_1, \xi_2)B(x^2, x^4; \xi_1, \xi_3)d\eta_1d\xi_1d\xi_2d\xi_3.
\eeq
Similar to the proof of Lemma \ref{prod3} and \ref{prod41}, we introduce a cut-off function to stay away from $\tilde \La_2, \tilde \La_3$. The rest is similar to the proof of Lemma \ref{prod51}. 
\epf

\begin{lemma}\label{prod53}
Let $u\in  I^{p, l}(\La_{12}, \La_2), v \in  I^{\mu'}(\La_{134})$ with $p, \mu'< -\frac{3}{2}, l<0$. For $\eps>0$ sufficiently small, we can write  $w = uv$ as
\beq
w = w_0 + w_1, \ \ w_0\in I^{p+\mu'+1}(\La_{q_0}), \ \ w_1 \in \mcd'(M; \La_2(\eps)).
\eeq
Moreover, the principal symbol of $w_0$ satisfies
\beq
\sigma_{\La_{q_0}}(w_0) = (2\pi)^{-\ha}\sigma_{\La_{12}}(u) \ast_1 \sigma_{\La_{134}}(v).
\eeq
\end{lemma}

\bpf
It suffices to consider the distributions on model Lagrangians $\tilde \La_i$, see Lemma \ref{symp}. We  write 
\beq
u = u_0+u_1, \ \ u_0 = \int_{\mbr^2} e^{i(x^1 \xi_1 + x^2\xi_2)} A(x^3, x^4; \xi_1, \xi_2)d\xi_1d\xi_2,\ \ u_1\in I^p(\tilde \La_2),
\eeq
with $A\in S^{M, M'}(\mbr^2_{x^3, x^4}; \mbr_{\xi_2}, \mbr_{\xi_1})$ being a symbol of product type where $M = p-\fnf + \ha + 1, M' = l-1$. Also, $u_1$ is supported away from $K_{134}$. Next, we write 
\beq
v = \int_{\mbr^3} e^{i(x^1 \xi_1 + x^3\xi_3+ x^4\xi_4)} B(x^2; \xi_1, \xi_3, \xi_4)d\xi_1d\xi_3d\xi_4,  
\eeq
with $B\in S^{m}(\mbr_{x^2}; \mbr^3_{\xi_1, \xi_3, \xi_4})$ being a standard symbol where $m =\mu' -\fnf + \frac{3}{2}$.
 The product of $u_1v$ is smooth because of the disjoint support. We only need to find the product
\beq
 u_0v = \int_{\mbr^5} e^{i(x^1 (\xi_1+ \eta_1) + x^2\xi_2 + x^3\xi_3 + x^4\xi_4)} A(x^3, x^4; \eta_1, \xi_2)B(x^2; \xi_1, \xi_3, \xi_4) d\eta_1d\xi_1d\xi_2d\xi_3d\xi_4.
\eeq
Again we introduce cut-off functions to separate the singularities. The rest is similar to the proof of Lemma \ref{prod51} as well.
\epf

\section{Proof of the main results}\label{pfmain}
We continue from Section \ref{disto} to finish the proof of our main theorems in the introduction. As we already mentioned, the arguments are very similar and the difference is that we make use of singularities in higher order asymptotic expansions analyzed in Section \ref{singuhigh}. We first show that the conformal class can be determined, which is part of Theorem \ref{maingh}. Moreover, we prove a more general result including linear terms  in the wave operators.
  
\begin{theorem}\label{mainconf}
Let $(M^{(j)}, g^{(j)}), j = 1, 2$ be two $4$-dimensional globally hyperbolic Lorentzian manifolds. Let $\hat \mu^{(j)}(t) \subset M^{(j)}$ be time-like geodesics where $t\in [-1, 1]$ and $V^{(j)}\subset M^{(j)}$ be open relatively compact neighborhood of $\hat \mu^{(j)}([s_-, s_+])$ where $-1<s_-<s_+<1$. Let $M_0^{(j)} = (-\infty, T_0)\times N^{(j)}, T_0 > 0$ such that $V^{(j)}\subset M_0^{(j)}$. Consider the semilinear wave equations with source terms
\beqq\label{eqsem2}
\begin{gathered}
\square_{g^{(j)}} u(x) + Q^{(j)}(x)u(x)  +  H^{(j)}(x, u(x))  = f(x), \text{ on } M^{(j)}_0,\\
u = 0 \text{ in } M_0\backslash J^+_{g^{(j)}}(\supp(f)),
\end{gathered}
\eeqq
where $\supp(f)\subset V^{(j)}$ and $Q^{(j)}(x)$ are smooth functions.  We assume that $H^{(j)}(x, z)$ are genuinely nonlinear on $I(p_-^{(j)}, p_+^{(j)})$ respectively, where $p^{(j)}_\pm = \hat\mu^{(j)}(s_\pm)$. Suppose that there is a diffeomorphism $\Phi: V^{(1)}\rightarrow V^{(2)}$ such that $\Phi(p^{(1)}_\pm) = p^{(2)}_\pm$ and the source-to-solution maps $L^{(i)}$ satisfy
\beq
(\Phi^{-1})^*(L^{(1)}(\Phi^*f)) = L^{(2)}(f)
\eeq
for all $f$  in a small neighborhood of the zero function in $C_0^4(V^{(2)})$. 
Then there exists a diffeomorphism $\Psi: I(p^{(1)}_-, p^{(1)}_+)\rightarrow I(p^{(2)}_-, p^{(2)}_+)$ such that the metric $\Psi^*g^{(2)}$ is conformal to $g^{(1)}$ in $I(p^{(1)}_-, p^{(1)}_+)$. Moreover, $\Psi = \Phi$ on $V^{(1)}\cap I(p^{(1)}_-, p^{(1)}_+)$.
\end{theorem} 
\bpf
As in the proof of Theorem \ref{main11}, it suffices to consider one manifold $(M, g)$. We let $\tilde Q_g$ be the causal inverse of $\square_g + Q$ where $Q$ is a smooth function. Notice that the principal symbol of $\tilde Q_g$ is the same as that of $Q_g = \square_g^{-1}$. Therefore, all of our analysis for $Q_g$ so far works for $\tilde Q_g$, and we shall abuse the notation below by taking $\tilde Q_g = Q_g$.  We divide the proof into three steps.

\textbf{Step 1:} If there are $k\geq 4$ such that $h_k(x)$ is non-vanishing at $q_0\in I(p_-, p_+)$, we know from Section \ref{sechs} that there are sources $f_i$ such that $\mcu^{(k)}$ has non-vanishing singularities at  $\La_{q_0}^g\backslash (\La^{(3), g}\cup \La^{(1), g})$. It follows from the same argument as in Theorem \ref{main11} (consider the set when $s_0\rightarrow 0$) that we can determine $\mce_V(q_0)$ for such $q_0$.

\textbf{Step 2:} If $h_k(x)$ vanishes in a neighborhood of $q_0$ for all $k\geq 4$, we can assume that $H(x, z) = az^2 + bz^3$. Actually, the case when $a$ is non-vanishing was considered already in Theorem \ref{main11}. So we focus on the case when $a$ vanishes but $b$ is non-vanishing i.e.\ $H(x, z) = b(x)z^3$. We've seen that in this case $\mcu^{(4)}$ does not produce new singularities. Although $\mcu^{(3)}$ contains conic singularities, it is not clear if they can propagate back to $V$. The idea below is to consider singularities in $\mcu^{(5)}$ as defined in \eqref{highnon}, which are produced due to non-vanishing $b$. 

We compute the asymptotic term $\mcu^{(5)}$ for $H(x, u(x)) = b(x)u(x)^3$ following the same approach as in Section \ref{semiwave}. We have
\beq
u = v - Q_g(bu^3).
\eeq
Below we use $\mcr$ to denote terms in $H^{4}(M)$ and not of the order $\eps_1^2\eps_2\eps_3\eps_4$. We compute
 \beq
\begin{split}
u^3 &= v^3 - 3v^2Q_g(bu^3) + \mcr = v^3 - 3v^2Q_g(bv^3)+ \mcr.
\end{split}
\eeq
Therefore, we obtain 
\beq
u = v - Q_g(bv^3) + 3Q_g(bv^2Q_g(bv^3)) + \mcr.
\eeq
Our $\mcu^{(5)}$ consists of terms in $3Q_g(bv^2Q_g(bv^3))$. In particular, we get 
\beq
\mcu^{(5)} =  \sum_{i,j,k,l,m}3Q_g(bv_iv_jQ_g(bv_kv_lv_m)),
\eeq
where the summation in $(i, j, k, l, m)$ is over permutations of $(1, 1, 2, 3, 4).$ These terms involve multiplication of $I^*(\La_{1i}, \La_i)$ and $I^*(\La_{1jk})$ where $(i, j, k)$ are permutations of $(2, 3, 4)$. These are analyzed in Lemma \ref{prod52}. The principal symbols of these terms can be calculated as in Section \ref{nonsym}, using the results in Section \ref{multising}.  It suffices to consider three model cases: 
\begin{enumerate}
\item[(a)] $i, j = 1$ and $(k, l, m) = (2, 3, 4)$; 
\item[(b)] $i=1, j = 2$ and $(k, l, m) = (1, 3, 4)$; 
\item[(c)] $i = 2, j = 3$ and $(k, l, m) = (1, 1, 4)$. 
\end{enumerate}
We first find the principal symbols of $\mco_{ijklm} =  bv_iv_jQ_g(bv_kv_lv_m)$ in each case. Then we find and show that the principal symbol of $\mcu^{(5)}$ is non-vanishing. As the order of the distributions is not important for the analysis below, we shall use $*$ to replace the orders.

In case (a), consider $\mco_{11234} =  bv_1^2Q_g(bv_2v_3v_4)$. We know from Prop.\ \ref{inter3} that $Q_g(bv_2v_3v_4)\in I^{*,*}(\La_{234}, \La^g_{234})$ and the principal symbol is a product of $\sigma(v_2), \sigma(v_3), \sigma(v_4)$ with positive coefficient. Also, we know from Lemma \ref{prod51} that $v_1^2 \in I^{*}(\La_1)$ and the principal symbol is a partial convolution of $\sigma(v_1)$. By Lemma \ref{prod42}, we know that $Q_g(\mco_{11234})\in I^*(\La_{q_0}\backslash \Xi)$  with $\Xi$ defined in Prop.\ \ref{inter4}. Now as in Section \ref{sechs}, we take $\zeta \in L_{q_0}^{*, +}M$ and $\zeta = \sum_{i = 1}^4 \zeta_i$ with $\zeta_i \in N^*_{q_0}K_i$. Let $A_i(q_0, \zeta_i), i = 2, 3, 4$ be the principal symbols of $v_2, v_3, v_4$. Following notations in Section \ref{sechs}, we let $A^{(1)}_{1}(q_0, \zeta_1)$ be the principal symbol of $v_1^2$, so $A^{(1)}_1 = \sigma(v)\ast \sigma(v)$ is the convolution in the fiber variable of $\La_1.$ Then the principal symbol can be written as
\beq
\sigma(\mco_{11234})(q_0, \zeta) = (2\pi)^{-\frac{7}{2}}b(q_0)^2|\zeta_2 + \zeta_3 + \zeta_4|_{g^*(q_0)}^{-2} \cdot A^{(1)}_{1}(q_0, \zeta_1) \prod_{i = 2}^4 A_i(q_0, \zeta_i),
\eeq
see the principal symbol of term $\mcy_3$ in the proof of Prop.\ \ref{nonvan}. 

In case (b), consider $\mco_{12134} =  bv_1v_2Q_g(bv_1v_3v_4)$we know from Lemma \ref{prod2} that $v_1v_2\in I^{*,*}(\La_{12}, \La_1) + I^{*,*}(\La_{12}, \La_2)$ with principal symbol the product of $\sigma(v_1), \sigma(v_2)$. Also, Prop.\ \ref{inter3} tells that $Q_g(bv_1v_3v_4)\in I^{*,*}(\La_{134}, \La^g_{134})$ microlocally away from $\La^{(1)}$ with principal the product of $\sigma(v_1), \sigma(v_3), \sigma(v_4)$. Now we use Lemma \ref{prod53} to conclude that $Q_g(\mco_{12134}) \in I^*(\La_{q_0}\backslash \Xi)$. The principal symbol is a partial convolution over the fiber variable of $\La_1$. In particular, we find that 
\beqq\label{sym12134}
\sigma(\mco_{12134})(q_0, \zeta) = (2\pi)^{-\frac{7}{2}} b(q_0)^2 \int \frac{\sigma(v_1)(q_0, \zeta_1 - \eta) \sigma(v_1)(q_0, \eta)}{|\eta + \zeta_3 + \zeta_4|_{g^*(q_0)}^{2}} d\eta \cdot \prod_{i = 2}^4 A_i(q_0, \zeta_i),
\eeqq
where the integration is over the fiber of $\La_1$ at $q_0.$

In case (c),  consider $\mco_{23114} =  bv_2v_3Q_g(bv_1^2 v_4)$. We first use Lemma \ref{prod2} to conclude that $v_2v_3\in I^{*,*}(\La_{23}, \La_2) + I^{*,*}(\La_{23}, \La_3)$. Using Lemma \ref{prod51} and Lemma \ref{prod2}, we know that $v_1^2v_4\in I^{*,*}(\La_{14}, \La_1) + I^{*,*}(\La_{14}, \La_4)$. Now we apply Lemma \ref{prod52} to conclude that  $Q_g(\mco_{23114}) \in I^*(\La_{q_0}\backslash \Xi)$. The principal symbol can be found as for term $\mcy_2$ in the proof of Prop.\ \ref{nonvan} 
\beq
\sigma(\mco_{23114})(q_0, \zeta) = (2\pi)^{-\frac{7}{2}} b(q_0)^2 |\zeta_1 + \zeta_4|_{g^*(q_0)}^{-2} \cdot A^{(1)}_{1}(q_0, \zeta_1) \prod_{i = 2}^4 A_i(q_0, \zeta_i).
\eeq
Finally, we can write down the symbol of $\mcu^{(5)}$ at $(q_0, \zeta)$ as
\beqq\label{symu5}
\begin{split}
\sigma(\mcu^{(5)})(q_0, \zeta)  = &\sum_{i, j, k, l, m} \sigma(\mcu^{(5)}_{ijklm})(q_0, \zeta) \\
 = & 3\cdot 6 \cdot \sigma(\mco_{11234})(q_0, \zeta) + 3\sum_{i, j, k}  \sigma(\mco_{1i1jk})(q_0, \zeta) +  3\sum_{i, j, k}  \sigma(\mco_{ij11k})(q_0, \zeta)\\
  = & (2\pi)^{-\frac{7}{2}} b(q_0)^2 [18 |\zeta_2 + \zeta_3 + \zeta_4|_{g^*(q_0)}^{-2} + 6\sum_{k=2}^4  |\zeta_1 + \zeta_k|_{g^*(q_0)}^{-2}] \\
  &\cdot A^{(1)}_{1}(q_0, \zeta_1) \prod_{n = 2}^4 A_n(q_0, \zeta_n) +  2\sum_{i, j, k}  \sigma(\mco_{1i1jk})(q_0, \zeta),
\end{split}
\eeqq
where the summation in $(i, j, k)$ is over the permutations of $(2, 3, 4)$. The symbols of $\mco_{1i1jk}$ have a similar express as \eqref{sym12134}, which involves the partial convolutions. Now we show  that when the submanifolds $K_i$ intersect
in  a generic way, 
the symbol is non-vanishing on any open set of $\La_{q_0}\backslash \Xi$. We basically follow the same argument as in Prop.\ \ref{nonvan}. In particular, we prove that for the light-like vectors $\zeta_i$ constructed in \eqref{zetas}, the term with $|\zeta_1 + \zeta_4|_{g^*}^{-2}$ dominates as $\rho \rightarrow 0$ hence the symbol is non-vanishing. 

We'll estimate the symbols of $\mco_{1i1jk}$ because they involve the convolution. First we choose the symbol of $v_1$ such that $\sigma(v_1)(q_0, \eta), \eta\in N^*_{q_0}K_1\subset L^{*, +}M$ is supported in $|\eta|_{g^+} > \eps$ for some $\eps >0$ small and $\sigma(v_1)\geq0$. In particular, by changing $\eps$, we can assume $\sigma(v_1)$ is supported on $\eta\in\{ \alpha \zeta_1:  \alpha > \eps\}$. Now consider the symbol of $\mco_{1i1jk}$, which has a similar expression as \eqref{sym12134}. We can estimate the convolution kernel as  
\beq
\begin{gathered}
|\frac{1}{|\eta + \zeta_2 + \zeta_3|_{g^*(q_0)}^{2}}| \leq C \rho^{-2}, \ \ 
|\frac{1}{|\eta + \zeta_2 + \zeta_4|_{g^*(q_0)}^{2}}| \leq C \rho^{-1}, \\
|\frac{1}{|\eta + \zeta_3 + \zeta_4|_{g^*(q_0)}^{2}}| \leq C \rho^{-1},
\end{gathered}
\eeq
using the computations in the proof of Prop.\ \ref{nonvan}. Thus we conclude that 
\beq
|\sigma(\mco_{1i1jk})(q_0, \zeta)| \leq C \rho^{-2} b(q_0)^2   \cdot  A^{(1)}_{1}(q_0, \zeta_1) \prod_{i = 2}^4 A_i(q_0, \zeta_i),
\eeq
for some constant $C>0$. Notice that now the convolution is in $A^{(1)}_1(q_0, \zeta_1)$. By comparing the growth orders of each terms in  \eqref{symu5} as $\rho\rightarrow 0$, we find that the leading term is given by $|\zeta_1 + \zeta_4|_{g^*(q_0)}^{-2} = -\rho^{-10}$.  Hence we obtain that 
\beq
\sigma(\mcu^{(5)})(q_0, \zeta)  = (-3\rho^{-10} + o(\rho^{-10})) (2\pi)^{-\frac{7}{2}} b(q_0)^2  \cdot A^{(1)}_{1}(q_0, \zeta_1)\prod_{i = 2}^4 A_i(q_0, \zeta_i),
\eeq
which is non-vanishing if $A_i, i = 2, 3, 4$ are non-vanishing and $\rho$ sufficiently small. Since the symbol is an analytic function of $\zeta$, we showed that when the submanifolds $K_i$ intersect
in  a generic way, the symbol $\sigma(\mcu^{(5)})(q_0, \zeta)$
 is non-vanishing on any open set of $\La_{q_0}^g\backslash \Xi$ with $\Xi$ defined in Prop.\ \ref{inter4}. Therefore, we can determine $\mce_V(q_0)$ as in Theorem \ref{main11}. 

\textbf{Step 3:} Part (1) and (2) determines the earliest light observation set for a dense subset of $I(p_-, p_+)$ i.e.\ $\cup_{k\geq 4} \{q: h_k(q) \neq 0\} \cup \text{int}(\{q: h_k(q) = 0, k\geq 4\})$. We can finish the proof now as in Theorem \ref{main11}.
\epf

Now we prove the key relation of the conformal factor and the nonlinear terms. The proof follows the same idea in Theorem \ref{main13}.
\begin{theorem}\label{maincoef}
Let $g^{(1)}, g^{(2)}$ be two globally hyperbolic Lorentzian metrics on a $4$-dimensional smooth manifold $M$ such that $g^{(1)} = e^{2\gamma} g^{(2)}$ where $\gamma \in C^\infty(M)$ and $\gamma = 0$ on an open relatively compact set $V \subset M_0 = (-\infty, T_0)\times N, T_0 > 0$.  In particular, $g^{(1)}$ and $g^{(2)}$ are isometric on $V$.  Let $\hat \mu([-1, 1]) \subset V$ be a time-like geodesic, and $p_\pm = \hat \mu(s_\pm)$ where $-1<s_-<s_+<1$. Consider the semilinear wave equations with source terms
\beq 
\begin{gathered}
\square_{g^{(i)}} u(x) + Q^{(i)}(x)u(x) + H^{(i)}(x, u(x))  = f(x), \text{ on } M_0,\\
u = 0 \text{ in } M_0\backslash J^+_{g^{(i)}}(\supp(f)),
\end{gathered}
\eeq
where $\supp(f)\subset V$ and $Q^{(i)}$ are smooth functions.  Assume that $H^{(i)}(x, z), i = 1, 2$ are genuinely nonlinear and have Taylor expansions
\beq
H^{(i)}(x, z) = \sum_{k = 2}^N h_k^{(i)}(x) z^k + O(z^{N+1}).
\eeq
 Let $L^{(1)}, L^{(2)}$ be the source-to-solution map with respect to the metrics $g^{(1)}, g^{(2)}$ and they satisfy 
$L^{(1)}f = L^{(2)}f$ for all $f$ in a small neighborhood of the zero function in $C_0^4(V)$. Then for $x\in I(p_-, p_+)$, we have 
\beq
h^{(1)}_k(x) = e^{(k-3)\gamma(x)}h_k^{(2)}(x),  \ \ \forall k\geq 4. 
\eeq
Furthermore, if $Q^{(j)} = 0, j = 1, 2$, we have 
\begin{enumerate}
\item $h^{(1)}_2(x)h^{(1)}_3(x) = e^{-\gamma(x)}h^{(2)}_2(x) h^{(2)}_3(x)$;
\item $h^{(1)}_2(x) = e^{-\gamma(x)}h^{(2)}_2(x)$ if $h_3^{(1)} = 0$.
\end{enumerate}
\end{theorem} 

\bpf
We divide the proof into two steps.
 
\textbf{Step 1:} Consider the determination of $h_k, k\geq 4$. We assume $f = \sum_{i = 1}^4\eps_i f_i$ constructed as in Section \ref{sec42} and denote
\beq
\mcu^{(k), i} = \p^{k-3}_{\eps_1}\p_{\eps_2}\p_{\eps_3}\p_{\eps_4}L^{(i)}(f)|_{\{\eps_1=\eps_2=\eps_3=\eps_4 =0\}}, \ \ i = 1, 2.   
\eeq
We consider when $h_k^{(1)}  \neq 0$. For any $q_0 \in I(p_-, p_+)$, we will compare the principal symbols of $\mcu^{(k), 1}_0$ and $\mcu^{(k), 2}_0, k\geq 5$ on $\mce_{V}(q_0)\backslash (\mck^{(3)}\cup K^{(1)})$ for $s_0\rightarrow 0$ using the computations  in Section \ref{sechs}.  By abuse notations, we let 
\beq
Q_{g^{(j)}} = (\square_{g^{(j)}} + Q^{(j)})^{-1}, \ \ j = 1, 2.
\eeq
Notice that the principal symbols of $Q_{g^{(j)}}$ are the same as the principal symbols of $\square^{-1}_{g^{(j)}}$ and in particular, we can still apply Prop.\ \ref{symconf} for $Q_{g^{(j)}}$. For two conformal metrics $g^{(1)} =e^{2\gamma} g^{(2)}$ and $\gamma = 0$ on $V\subset M_0$, we have
\beq
\begin{gathered}
\sigma(Q_{g^{(1)}})(q, \eta, q_0, \xi) = \sigma(Q_{g^{(2)}})(q, \eta, q_0, \xi)e^{3\gamma(q_0)},\\
\sigma(Q_{g^{(1)}})(q_0, \xi_i, x_i, \zeta_i) = e^{-\gamma(q_0)}\sigma(Q_{g^{(2)}})(q_0, \xi_i, x_i, \zeta_i),
\end{gathered}
\eeq
where $(x_i, \zeta_i)\in T_V^*M$. With the calculation in Section \ref{sechs}, we obtain 
\beq
\begin{gathered}
A^{(1)}_i(q_0, \xi_i) = e^{-\gamma(q_0)}A^{(2)}_i(q_0, \xi_i) , \ \ i = 1, 2, 3, 4,\\
A^{(k-3), (1)}_1(q_0, \xi_1) = e^{-(k-3)\gamma(q_0)}A^{(k-3), (2)}_1(q_0, \xi_1).
\end{gathered}
\eeq
Notice that when computing the principal symbol $\sigma(v_1^{k-3})$, the partial convolution does not affect the conformal factor $e^{\gamma}$. Therefore, we obtain the following relation
\beq
\sigma(\mcu^{(k), 1}_0)(q, \eta) h_k^{(2)}(q_0) = h_k^{(1)}(q_0) e^{(-k + 3)\gamma(q_0)}\sigma(\mcu^{(k), 2}_0)(q, \eta).
\eeq
From $L^{(1)}f = L^{(2)}f$, we know the principal symbols of $\mcu_0^{(k), i}$ should be the same. Hence we get the relation of $h_{k}^{(i)}$.  
This finishes the proof of Step 1.

\textbf{Step 2:} It remains to determine $h_2, h_3$. For convenience, we change notations by setting $h_2(x) = a(x), h_3(x) = b(x)$. Recall that in Prop.\ \ref{inter5}, we found the leading singularities of $\mcu^{(4)}$. Below, we make use of the leading singularity of 
\beq
\mcu^{(4)}_1 \doteq \mcu^{(4)} - 4Q_g(h_4 v_1v_2v_3v_4) \in I^{\mu_0}(\La_{q_0}^g\backslash \La_{q_0}).
\eeq
Here, observe that since we determined $h_4$ from Step 1 and we know the manifold $(M, g)$, if there is no potential term, the term $4Q_g(h_4 v_1v_2v_3v_4)$ is already determined.  Notice that this is not just the leading singularity but the whole term. Therefore, if follows from Prop.\ \ref{inter5} the conclusions below 
\begin{enumerate}
\item[(i)] if $\mu_0 = 4\mu-\frac{5}{2}$, then $a\neq 0, b\neq 0$ at $q_0$; 
\item[(ii)] if $\mu_0 = 4\mu - \frac{9}{2}$, then $a\neq 0, b = 0$ at $q_0$; 
\end{enumerate}
Here we assume that Prop.\ \ref{inter5} is valid near $q_0$. In the end, we will use the density argument (as in the proof of Theorem \ref{main11}) to complete the proof. 

We start with case (ii). As seen in Section \ref{singu}, when $b$ vanishes near $q_0$, we have
\beq
\mcu^{(4)}_1 = -\sum_{(i, j, k, l)} [4Q_g(av_iQ_g(av_jQ_g(av_kv_l))) + Q_g(aQ_g(av_iv_j)Q_g(av_kv_l))],
\eeq 
where the summation is over the permutations of $(1, 2, 3, 4)$. 
This is part (3) of Theorem \ref{main13}, so we have $a^{(1)} = e^{-\gamma}a^{(2)}$ at $q_0$. This finishes the proof for part (3).

Now we consider case (i) which is part (2). In this case, we need the singularities in $\mcu^{(4)}_1$ 
\beq
\mcu^{(4)}_1 = \sum_{(i, j, k, l)} [2Q_g(av_iQ_g(bv_jv_kv_l)) + 3Q_g(bv_iv_jQ_g(av_kv_l))]
\eeq
with summation over permutations of $(1, 2, 3, 4)$. This is indeed part (2) of Theorem \ref{main13},  so we have 
\beq
a^{(1)}b^{(1)} =  e^{-\gamma}a^{(2)}b^{(2)}.
\eeq
Notice that this is also true for $b^{(1)} = 0$ (or $b^{(2)} = 0$). This finishes the proof.  
\epf

Finally, we finish the proof of all the main results stated in the introduction. 

\bpf[Proof of Theorem \ref{maingh}]
Using Theorem \ref{mainconf} and Remark 3.1 of \cite{KLU}, we know that $e^{\gamma} = 1$ on $V^{(1)}\cap I(p^{(1)}_-, p^{(1)}_+)$,  see also the proof of Theorem \ref{main11}, part (2). Now we can apply Theorem \ref{maincoef} to complete the proof. 
\epf

\bpf[Proof of Theorem \ref{main2}]
This is a direct consequence of Theorem \ref{maingh}.
\epf

\bpf[Proof of Theorem \ref{maingauge}] 
When the measurement domain $V$ and the map $L_{V;{g},H}$  are given, we can use Theorem \ref{maingh} (or Theorem \ref{mainconf}) to construct the conformal type of $g$. After that, let us choose some metric $g_0$ that is conformal to $g$. After fixing $g_0$, we can consider equation 
$$
\Y_{g_0} u(x) + H_0(x,u(x)) = f(x),
$$
that is a gauge transformation of the original equation. Since the zeroth order term $-\frac {n-4}{4(n-2)}R_{g_0}$ in the operator $\Y_{g_0}$ does not change the principal symbol of the parametrix $\Y_{g_0}^{-1}=(\square_g -\frac {n-4}{4(n-2)}R_{g_0})^{-1}$ of the Yamabe operator $\Y_{g_0}$ and we know the metric $g_0$, it follows from Theorem \ref{main2} that we can construct the non-linear term $H_0$ as
\beq
H_0(x, z) = \sum_{k = 4}^\infty \p_z^k H_0(x, 0)z^k
\eeq
under the assumption of the theorem i.e.\ $H_0(x, z)$ is real analytic in $z$ and $H_0(x, z) = O(z^4)$ as $z\rightarrow 0$. Thus we can construct the representative $\Y_{g_0}+H_0(x,\,\cdotp)$  of the gauge-equivalence class of the original operator $\Y_g+H(x,\,\cdotp)$. 
\epf

\bpf[Proof of Theorem \ref{main1}]
We proved the determination of the conformal class in Theorem \ref{mainconf}. To determine the conformal factor, the case when the Ricci curvatures are zero follows from the same geometric argument in Corollary 1.3 of \cite{KLU}. For the other  case when $H^{(i)}(x, z)$ are independent of $x$, we claim that $\p^k_z H^{(1)}(0) = \p^k_z H^{(2)}(0), k\geq 2$. Then by applying Theorem \ref{maingh}, we can determine the conformal factor in any case except when the nonlinear terms are purely cubic.

Actually, we prove a little bit more general statement than what we claimed. We'll show that when $H = H(x,z)$ depends also on $x$, the measurements in the set V determine the derivatives $\p_z^k H(x,z)|_{z=0}$ for all $x\in V$ and $k\geq 2$. 

Assume that we are given the open set $V$ and  the operator $f\mapsto L_{V,\square_g,H}(f)$ in a neighborhood $W\subset C^2_0(V)$ of the zero-function. As seen in the proof of Theorem \ref{main13} (also Remark 3.1 of \cite{KLU}), the derivative of the non-linear map  $f\mapsto L_{V,\square_g,H}(f)$ determines the source-to-field map of the linearized wave equation and it determines the metric $g|_V$ in $V$. Let 
$$
W_0=\{f\in W;\ \hbox{support of $L_{V,\square_g,H}(f)$  is compact subset of $V$}\}.
$$
Then the solutions $u^f$ corresponding to sources $f\in W_0$ vanish in a neighborhood of $\partial V$
and thus those vanish also outside $V$. Choosing all sources $f\in W$ such that $L(f)=u^f|_V$ is $C^4$-smooth in $V$ and compactly
supported in $V$, we can determine the set 
$$
\{u^f|_V\in C^4(V);\ f\in W_0\}=\{v\in C^4_0(V);\ \square_g v\in C_0^2(V)\cap W\}
$$
Now, there is $\delta_0>0$ such that
$$
S(\delta_0)=\{v\in C^4_0(V);\ \|v\|_{C^4(\overline V)}<\delta_0\}\subset \{u^f|_V\in C^4(M_0);\ f\in W_0\}=\{L(f) \in C^4_0(V);\ f\in W_0\}.
$$
Using the fact that the  set $V$ and  the operator $f\mapsto L_{V,\square_g,H}(f)$ are given
we can first determine $\delta_0>0$ such that the above is valid and then determine the inverse of the map $L|_{W_0}:f\to L(f)$, defined in the set $S(\delta_0)$.
That is, we can find the map $(L|_{W_0})^{-1}:S(\delta_0)\to W_0$.
Then, using the equation $\square_g v+H(x,v)=f$ we see that 
$$
\{(v,f-\square_g v);\ v\in S(\delta_0),\ f=(L|_{W_0})^{-1}v\}=\{(v,H(\,\cdotp,v));\ v\in S(\delta_0)\}.
$$
In other words, we can find the pairs $(v,H(\,\cdotp ,v))\in C^4_0(V)\times C^2_0(V)$  for all sufficiently small $v\in C^4_0(V)$, i.e.\,
we find the graph of the map $v\mapsto H(\,\cdotp,v)$ close to zero function. Hence we can determine the derivatives of the function $z\mapsto H(x,z)$ for $x\in V$.  
\epf

%
%

\section{Acknowledgements}
The authors would like to thank Peter Hintz for his careful reading of a preliminary version of the paper and for very valuable comments. They also thank Andr\'as Vasy for helpful discussions. G.U.\ was partially supported by NSF and Academy of Finland, FiDiPro-project. M.L.\ was partially supported by the Academy of Finland.


\end{document}